\documentclass[12pt]{amsart}

\usepackage{amsmath,amsfonts,amssymb}

\theoremstyle{plain}
\newtheorem{lemma}{Lemma}[section]
\newtheorem{theorem}[lemma]{Theorem}

\newtheorem{proposition}[lemma]{Proposition}

\theoremstyle{definition}
\newtheorem{remark}[lemma]{Remark}

\numberwithin{equation}{section}
\DeclareMathOperator{\diam}{diam}
\DeclareMathOperator{\Co}{Co}
\DeclareMathOperator{\Var}{Var}

\DeclareMathOperator{\Int}{Int}
\DeclareMathOperator{\Ext}{Ext}
\DeclareMathOperator{\HSkel}{HSkel}
\DeclareMathOperator{\HPath}{HPath}

\DeclareMathOperator{\Open}{Open}

\DeclareMathOperator{\main}{main}
\DeclareMathOperator{\Mid}{Mid}

\begin{document}

\newcommand{\ZZ}{\mathbb{Z}^{2}}
\newcommand{\RR}{\mathbb{R}^{2}}
\newcommand{\tL}{\tau_{\mathcal{L}}}
\newcommand{\lra}{\leftrightarrow}
\newcommand{\omB}{\overline{\mathcal{B}}}
\newcommand{\kt}{\kappa_{\tau}}
\newcommand{\ep}{\epsilon}


\title[Droplets in Random Cluster Models]
{Cube--Root \\
Boundary Fluctuations for Droplets \\
in Random Cluster Models}
\author{Kenneth S. Alexander}
\address{Department of Mathematics DRB 155\\
University of Southern California\\
Los Angeles, CA  90089-1113}
\email{alexandr@math.usc.edu}
\thanks{Research supported by NSF grant DMS-9802368.}

\date{\today}

\begin{abstract}
For a family of bond percolation models on $\ZZ$
that includes the Fortuin-Kasteleyn
random cluster model, we consider properties of the ``droplet''
that results, in the percolating regime, from conditioning on the 
existence of an open dual circuit surrounding the origin and enclosing
at least (or exactly) a given large area $A$.  
This droplet is a close surrogate 
for the one obtained by Dobrushin, Koteck\'y and Shlosman by 
conditioning the Ising model; it approximates an area-$A$
Wulff shape.  The local part of the deviation from the Wulff shape
(the ``local roughness'') is the inward deviation of 
the droplet boundary from the boundary of its own convex hull; the
remaining part of the deviation, that 
of the convex hull of the droplet from the Wulff shape,
is inherently long-range.  We show that the local roughness is described
by at most the exponent 1/3 predicted by nonrigorous theory; this same
prediction has been made for a wide
class of interfaces in two dimensions.  Specifically, the average
of the local roughness over the droplet surface is shown to be
$O(l^{1/3}(\log l)^{2/3})$ in
probability, where $l = \sqrt{A}$ is the linear scale of 
the droplet.  We also bound the maximum of the local roughness over
the droplet surface and bound the long-range part of the deviation from
a Wulff shape, and we establish the absense of ``bottlenecks,'' which
are a form of self-approach by the droplet boundary, down to scale
$\log l$.  Finally, if we condition instead on the event that the total
area of all large droplets inside a finite box exceeds $A$, we show
that with probability near 1 for large $A$, 
only a single large droplet is present.
\end{abstract}
\maketitle

\section{Introduction} \label{S:intro}
Consider an Ising model in a finite box
(or other ``nice'' region) $\Lambda$ in $\mathbb{Z}^{2}$
at a supercritical inverse temperature $\beta$, with minus boundary condition.
There is a positive magnetization $m(\beta)$, and if $\Lambda$ is
large, the expected and actual fraction of plus
spins observed in $\Lambda$ will
each be approximately $(1 - m(\beta))/2$, with high probability.  If,
however, one conditions on the observed number of plus spins
 being sufficiently
greater than the expected number, then, as first explicated by Dobrushin,
Koteck\'y and Shlosman \cite{DKS}, the typical configuration contains a
single macroscopic droplet of the plus phase, that is, a droplet in which
the usual
proportions of plus and minus spins are reversed. Further, the droplet
will have a characteristic equilibrium crystal shape, the solution of an
isoperimetric problem, given by the Wulff construction.  If, for example, the
temperature is such that the usual proportion of plus to minus
spins is 20/80, and
one conditions the fraction of plus spins to be 30\%, then the
typical configuration
will show a 20/80 mix except inside a Wulff droplet in which the proportion
is approximately 80/20.  This droplet will cover approximately
1/6 of the box $\Lambda$,
so as to account for nearly all of the excess plus spins.

This in an example of the general phenomenon of phase separation.  The
work of Dobrushin, Koteck\'y and Shlosman
in \cite{DKS} provides the first rigorous derivation of phase
separation beginning from a local interaction.

In the joint construction \cite{ES} of the Ising model and the
corresponding Fortuin-Kasteleyn random cluster model
(\cite{FK}; see \cite{Gr95}), abbreviated ``FK model,'' the droplet
boundary appears as a circuit of open dual bonds.  
If a  particular site, say the
origin, is inside the droplet, one expects that  
the outermost open dual circuit
$\Gamma_{0}$ containing the origin will closely 
approximate the droplet boundary.  Since the Ising droplet has
approximately a fixed area, we can gain information
useful in studying the droplet by studying the FK model conditioned
to have $\Gamma_{0}$ enclose at least, or exactly, a given area $A$.
This is our main aim in this paper.

Since it is of interest to study phase separation beyond the context 
of the Ising model, we establish our results not just for the FK model 
but for general percolation models having properties known, or 
reasonably expected, to hold quite widely
in the percolating regime.  These include the FKG 
property, a special case of the Markov property, exponential decay of
dual connectivity
and certain mixing properties.

In \cite{DKS} for very low temperatures, and in \cite{IS97} for all
subcritical temperatures, bounds are given
for the boundary fluctuations of the Ising droplet, that is, 
for the deviation of the
boundary of the  observed droplet from 
the boundary of an appropriately translated and rescaled Wulff (that is, 
equilibrium crystal) 
shape.  Let $N$ denote the linear scale of the box.  For a 
droplet also of linear scale $N$, the boundary 
fluctuations are shown to be of order at most $N^{3/4}\sqrt{\log N}$, and 
for a droplet of linear scale $l_{N} = N^{\alpha}$, with 
$2/3 < \alpha < 1$, the boundary fluctuations are shown to be 
of order less than $\sqrt{N^{2/3}l_{N}}$.
There is typically never a 
droplet of linear 
scale $N^{\alpha}$ for $0 < \alpha < 2/3$; if the excess number of plusses 
is of 
order $N^{2\alpha}$ with $0 < \alpha < 2/3$, then these plusses are 
dispersed 
throughout the minus phase without the formation of a 
large droplet.  Also, when 
phase separation does occur, other than the single large droplet there are 
typically no droplets of linear scale greater than log $N$.

Heuristics suggest that the boundary fluctuation bounds of 
\cite{DKS} and \cite{IS97} are not sharp.  To see what the correct
fluctuation size should be, one must refine the analysis by considering
three separate types of fluctuations.  The first type is 
\emph{shrinkage}---the 
actual droplet may be smaller than an ideal ``full size''
Wulff shape large enough to account for
all excess plus spins, since some of the excess may be 
dispersed in the surrounding minus
phase.   So one should actually consider fluctuations about a 
shrunken Wulff shape
enclosing the same area as the actual droplet.  The second type is 
\emph{local roughness}, defined as 
inward deviations of the droplet boundary from the boundary of
its own convex hull.  The third type is \emph{long-wave fluctuations},
defined as deviations of the convex hull of the droplet from the
shrunken Wulff shape.  (More precise definitions will be made below.)

The local roughness is of particular interest, and is our main focus in
this paper, because it is subject to the same type of 
interface-roughness heuristics as a wide variety of other dynamic and
equilibrium systems, 
including first-passage percolation \cite{NP}, \cite{LNP}, various
deposition models \cite{KS91}, polymers in random environments \cite{Pi},
asymmetric exclusion processes \cite{KS91}
and longest increasing subsequences of random permutations \cite{BDJ},
only the last of which is now well-understood
rigorously.  For a two-dimensional object of linear scale $l$, 
these heuristics predict fluctuations of order $l^{1/3}$ and a transverse
correlation length of order $l^{2/3}$.  For the local roughness, this 
transverse correlation length should appear as the typical separation 
between adjacent extreme points, where the 
droplet boundary touches the boundary of
its convex hull (see Section \ref{S:main}.)  
This is what makes local roughness local---one expects distinct inward 
excursions of the droplet boundary from the convex hull boundary to interact
only minimally.
The main result of this paper is that in
the random cluster model context, with probability 
approaching 1 as $l \to \infty$, 
the average local roughness is $O(l^{1/3}(\log l)^{2/3})$.

Results of Dobrushin and Hryniv \cite{DH} and Hryniv \cite{Hr} (at very low
temperatures) strongly
suggest that the the fluctuations of the droplet boundary about the shrunken 
Wulff shape should be Gaussian, heuristically resembling roughly a rescaled 
Brownian bridge added radially to the Wulff shape.  In particular, the
long-wave fluctuations should be of order $l^{1/2}$.  We are only able to
bound these by $l^{2/3}(\log l)^{1/3}$, however, 
in the random cluster model context.

\section{Definitions, Heuristics and Statement of Main Results} \label{S:main}
The results in this paper make use of only a few basic properties
of the FK or other percolation model, so we will state 
our results for general bond percolation models satisfying these properties.
A \emph{bond}, denoted $\langle xy \rangle$, is an unordered pair of 
nearest neighbor sites of $\ZZ$.  When convenient 
we view bonds as being open line segments in the plane; this should be
clear from the context.  In particular for $R \subset 
\mathbb{R}^{2}, \ \mathcal{B}(R)$
denotes the set of all bonds for which the corresponding open line segments 
are contained in $R$, and when we refer to distances between sets of
bonds, we mean distances between the corresponding sets of line
segments.  The exception is for $\Lambda \subset
\ZZ$,  for which we set
$\mathcal{B}(\Lambda) = \{\langle xy \rangle: x, y \in \Lambda\}$.
(Again, this should be clear from the context.)
For a set $\mathcal{D}$ of bonds we let 
$V(\mathcal{D})$ denote the set of all endpoints of bonds in 
$\mathcal{D}$, and 
\[
  \partial \mathcal{D} = \{\langle xy \rangle: x \in V(\mathcal{D}),
  y \notin V(\mathcal{D})\}, \quad \overline{\mathcal{D}} = 
  \mathcal{D} \cup \partial \mathcal{D}.
\]
We write $\omB(\Lambda)$ for 
$\overline{\mathcal{B}(\Lambda)}$.
A \emph{bond configuration} is an element $\omega \in 
\{0,1\}^{\mathcal{B}(\mathbb{Z}^{2})}$.

The \emph{dual lattice} is the translation of the integer lattice 
by (1/2,1/2); we write $x^{*}$ for $x + (1/2,1/2)$.
To each (regular)
bond $e$ of the lattice there corresponds a \emph{dual bond}
$e^{*}$ which is its perpendicular bisector; the dual bond is 
defined to be open in a configuration $\omega$ precisely 
when the regular bond is closed, and the corresponding configuration
of dual bonds is denoted $\omega^{*}$.  We write
$(\ZZ)^{*}$ for $\{x^{*}: x \in \ZZ\}$.
A \emph{cluster} in a given configuration is a connected 
component of the graph with site set $\ZZ$ and all open
bonds; \emph{dual clusters} are defined analogously for open dual bonds.
(In contexts where there is a boundary condition consisting of 
a configuration on the complement
$\mathcal{D}^{c}$ for some set $\mathcal{D}$ of bonds,
a cluster may include bonds in $\mathcal{D}^{c}$.)  
$C_{x}$ and $C_{x^{*}}$ denote the regular and dual clusters 
containing sites $x$ and $x^{*}$, respectively.
Given a set $\mathcal{D}$ of bonds, we write $\mathcal{D}^{*}$ for 
$\{e^{*}: e \in \mathcal{D}\}$.  The set of all endpoints of bonds in
$\mathcal{D}^{*}$ is denoted
$V^{*}(\mathcal{D})$ or $V^{*}(\mathcal{D}^{*})$. 

For $\Lambda \subset \ZZ$ or $\Lambda \subset (\ZZ)^{*}$ we define
\[
  \partial \Lambda = \{x \notin \Lambda: x \text{ adjacent to } 
  \Lambda\}, \quad \partial_{in}\Lambda = 
  \{x \in \Lambda: x \text{ adjacent to } 
  \Lambda^{c}\}
\]
where adjacency is in the appropriate lattice $\ZZ$ or $(\ZZ)^{*}$.

A \emph{(dual) path} is a sequence $\gamma = 
(x_{0},\langle x_{0}x_{1} \rangle, x_{1},\ldots
x_{n-1},\langle x_{n-1}x_{n} \rangle, x_{n})$ of alternating (dual)
sites and bonds.  $\gamma$ is \emph{self-avoiding} if all sites are
distinct.  We write $x \lra y$ (in $\omega$) 
if there is a path of open bonds (or open dual 
bonds, if $x$ and $y$ are dual sites) from $x$ to $y$ in $\omega$.
A \emph{circuit} is a path with $x_{n} = x_{0}$
which has all bonds distinct and which does not ``cross itself'' (in 
the obvious sense.)  Note we do allow $x_{i} = x_{j}$ for any $i \neq j$
here, i.e. a circuit may touch itself without
crossing.  A path or circuit is \emph{open} in a bond configuration
$\omega$ if all its bonds are open.
The \emph{exterior} of a circuit
$\gamma$, denoted $\Ext(\gamma)$, is the
unique unbounded component of the complement of $\gamma$ in 
$\mathbb{R}^{2}$, and the \emph{interior} $\Int(\gamma)$ is the
union of the bounded components.  An open circuit $\gamma$ is called an 
\emph{exterior circuit} in a configuration $\omega$
if $\gamma \cup \Int(\gamma)$ is
maximal among all open circuits in $\omega$.  (These definitions
differ slightly from what is common in the literature.)  Similar
definitions apply to dual circuits.
A site $x$ is surrounded by
at most one exterior circuit; when this circuit exists we denote it
$\Gamma_{x}$.  For $u, v$ points
in a path or circuit $\zeta$, let $\zeta^{[u,v]}$
and $\zeta^{(u,v)}$ denote the closed and open segments,
respectively, of $\zeta$ from $u$ to $v$ (in the direction of positive 
orientation, for circuits.)  $| \cdot |$ denotes 
the Euclidean norm for vectors, Euclidean length for curves,
cardinality for finite sets, and Lebesgue measure
for regions in $\mathbb{R}^{2}$
(which one should be clear from the context.)  Euclidean distance
is denoted $d(\cdot, \cdot)$.
Define $d(A,B)
= \inf\{d(x,y): x \in A, y \in b\}$ for $A, B \subset \mathbb{R}^{2}$
and $d(x,A) = d(\{x\},A)$.  
We define the \emph{average local roughness}
of a circuit $\gamma$ by
\[
  ALR(\gamma) = \frac{|\Co(\gamma) \backslash \Int(\gamma)|}
  {|\partial \Co(\gamma)|},
\]
where $\Co(\cdot)$ denotes the convex hull.  The \emph{maximum local
roughness} is
\[
  MLR(\gamma) = sup\{d(x,\partial\Co(\gamma)): x \in \gamma\}.
\]

By a \emph{bond percolation model} we mean a probability 
measure $P$ on 
$\{0,1\}^{\mathcal{B}(\mathbb{Z}^{2})}$.  
The conditional
distributions for the model $P$ are
\[
  P_{\mathcal{D},\rho} = P(\cdot \mid \omega_{e}
  = \rho_{e} \text{ for all } e \in \mathcal{D}^{c}),
\]
where $\mathcal{D} \subset \mathcal{B}(\ZZ)$.  
We say a bond percolation
model $P$ has \emph{bounded energy} if there exists $p_{0} > 0$ such that
\begin{equation} \label{E:boundener}
  p_{0} < P(\omega_{e} = 1 \mid \omega_{b}, b \neq e) < 1 - p_{0} \quad 
  \text{for all } \{\omega_{b}, b \neq e\}.
\end{equation}
From \cite{BuK}, bounded energy and translation invariance imply that 
there is at most one infinite cluster $P$-a.s.
Write $\omega_{\mathcal{D}}$ for $\{\omega_{e}: e \in \mathcal{D}\}$
and let $\mathcal{G}_{\mathcal{D}}$ denote the 
$\sigma$-algebra generated by $\omega_{\mathcal{D}}$.
$P$ has the \emph{weak mixing property} if for some $C, \lambda > 0$, for all
finite sets $\mathcal{D},\mathcal{E}$ with $\mathcal{D} \subset \mathcal{E}$,
\begin{align} 
  \sup \{\Var(&P_{\mathcal{E},\rho}(\omega_{\mathcal{D}} \in \cdot),
    P_{\mathcal{E},\rho^{\prime}}(\omega_{\mathcal{D}} \in \cdot)): 
    \rho, \rho^{\prime} \in \{0,1\}^{\mathcal{E}^{c}}\} \notag \\
  &\leq C \sum_{x \in V(\mathcal{D}),y \in V(\mathcal{E}^{c})} 
    e^{-\lambda |x - y|},
\notag
\end{align}
where $\Var(\cdot,\cdot)$ denotes total variation distance between measures.  
Roughly, the influence of the boundary condition on a finite region decays 
exponentially with distance from that region.  Equivalently, for some 
$C, \lambda > 0$, for all sets $\mathcal{D}, \mathcal{F} \subset 
\mathcal{B}(\ZZ)$,
\begin{align} \label{E:weakmix}
  \sup \{|&P(E \mid F) - P(E)|:  E \in \mathcal{G}_{\mathcal{D}}, F \in
    \mathcal{G}_{\mathcal{F}}, P(F) > 0\} \\
  &\leq C \sum_{x \in V(\mathcal{D}),y \in V(\mathcal{F})} 
    e^{-\lambda |x - y|}. \notag
\end{align}
$P$ has the \emph{ratio weak mixing property} if for some $C, \lambda > 0$,
for all sets $\mathcal{D}, \mathcal{F} \subset 
\mathcal{B}(\ZZ)$,
\begin{align} \label{E:rweakmix}
  \sup &\left\{ \left| \frac{P(E \cap F)}{P(E)P(F)} - 1 \right| : E \in 
    \mathcal{G}_{\mathcal{D}}, F \in
    \mathcal{G}_{\mathcal{F}}, P(E)P(F) > 0 \right\} \\
  &\leq C \sum_{x \in V(\mathcal{D}),y \in V(\mathcal{F})}
    e^{-\lambda |x - y|}, \notag
\end{align}
whenever the right side of (\ref{E:rweakmix}) is less than 1.

Let $\Open(\mathcal{D})$ denote the event 
that all bonds in $\mathcal{D}$ are open.

The FK model \cite{FK} with parameters $(p,q),\ p \in [0,1],q>0$
on a finite $\mathcal{D} \subset \mathcal{B}(\ZZ)$ is described
by a weight attached to each bond configuration $\omega \in 
\{0,1\}^{\mathcal{B}(\mathcal{D})}$, which is 
\[
  W(\omega) = p^{|\omega|}(1-p)^{|\mathcal{D}| - |\omega|}
  q^{C(\omega)},
\]
where $|\omega|$ denotes the number of open bonds in $\omega$ and
$C(\omega)$ denotes the number of open clusters in $\omega$, 
counted in accordance with the boundary condition, if any; see
\cite{Gr95} for details and further 
information.  For integer $q \geq 1$ the FK model is a random
cluster representation of the $q$-state Potts model at inverse 
temperature $\beta$ given by
$p = 1 - e^{-\beta}$.  For the study of phase separation involving
more than two species, for example in the Potts model, it is useful 
to be able to ``tilt'' the distribution with one or more external
fields before calculating various probabilities, 
as well as quantities such 
as surface tension and magnetization.  For the 
$q$-state Potts model with external fields
$h_{i}$ on species $i$, $i = 0,1,\ldots,q-1$, 
we need only consider $0 = h_{0} \geq h_{1} \geq \ldots \geq h_{q-1}$
and then the factor $q^{C(\omega)}$ in the weight $W(\omega)$ is
replaced by 
\[
  \prod_{C \in \mathcal{C}(\omega)} \left(1 + (1-p)^{h_{1}|C|} + \ldots
  + (1-p)^{h_{q-1}|C|}\right),
\]
where $\mathcal{C}(\omega)$ is the set of clusters in the
configuration $\omega$ and
$|C|$ denotes the number of sites in the cluster $C$.  We call
species $i$ \emph{viably dominant} if $h_{i}$ is maximal,
i.e. $h_{i} = h_{0}$.  For each 
species $i$ and for finite $\Lambda
\subset \ZZ$, corresponding to the species-$i$ boundary
condition for the Potts model on $\Lambda$ there is
the $i$-\emph{wired} boundary
condition on $\omB(\Lambda)$ for the FK model,
in which sites in $\Lambda$
connected to $\partial \Lambda$ are considered a 
single cluster $C_{\partial}$ and assigned weight
\[
  (1-p)^{h_{i}|C_{\partial}|}.
\]
Given a circuit $\gamma$ and a 
configuration $\rho \in \{0,1\}^{\mathcal{B}(\Ext(\gamma))}$, 
conditioning on $\rho$ and on $\Open(\gamma)$ induces a 
boundary condition on $\mathcal{B}(\Int(\gamma))$ which is a 
mixture over $i$ of the
different $i$-wired boundary conditions.  The weight assigned
to $i$-wiring in this mixture is proportional
to $(1-p)^{h_{i}N(\rho)}$, where 
$N(\rho)$ is the number of sites in $\gamma$ plus the number
of sites in $\Ext(\gamma)$ connected to $\gamma$ in $\rho$. In
the absense of external fields, $i$-wiring is the same for all 
$i$ and the choice of $\rho$ does not affect the boundary condition
induced on $\mathcal{B}(\Int(\gamma))$, which is a form of
Markov property, but if an external field 
is present this property fails.  However, the 
weight assigned to $i$-wiring in the mixture 
for non-viably-dominant $i$ is 
exponentially small in $|\gamma|$, uniformly in $\rho$, so for 
large $\gamma$ the effect of $\rho$ on the boundary condition
is uniformly small.

Motivated by the preceding, we say a bond percolation model
$P$ has the \emph{Markov property for open circuits} if for
every circuit $\gamma$ (of regular bonds), the bond 
configurations inside and 
outside $\gamma$ are independent given the event $\Open(\gamma)$.  
We have seen that the FK model has this 
property if and only if there are no external fields.  If $P$ is the
infinite-volume 
$k$-wired FK model for some viably dominant $k$, then
letting $\omega_{int}$
and $\omega_{ext}$ denote the bond configurations inside and outside
$\gamma$, respectively, we have from the preceding 
discussion for some $C, a > 0$
\begin{align} \label{E:Markov}
  \sup &\left\{ \left|\frac{P(\omega_{int} \in A \mid \Open(\gamma), \omega_{ext} 
    \in B)}{P(\omega_{int} \in A \mid \Open(\gamma))} - 1\right|: 
    A \in \mathcal{G}_{\mathcal{B}(\Int(\gamma))},
    B \in \mathcal{G}_{\mathcal{B}(\Ext(\gamma))}\right\} \\
  &\qquad \leq Ce^{-a|\gamma|}
    \qquad \text{for all } \gamma. \notag 
\end{align}
When (\ref{E:Markov}) holds we say $P$ has the \emph{near-Markov
property for open circuits}.
It is easy to see that one can interchange the roles of interior and exterior 
in (\ref{E:Markov}).  Further, if $\gamma_{1},..,\gamma_{k}$ are
circuits with disjoint interiors, $B_{i} \in 
\mathcal{G}_{\mathcal{B}(\Int(\gamma_{i}))}, A \in
\mathcal{G}_{\mathcal{B}(\cap_{i} 
\Ext(\gamma_{i}))}$, then
by easy induction on $k$,
\begin{equation} \label{E:Markov2}
  \frac{P(A \mid \Open(\gamma_{i}) \cap B_{i} \text{ for all } 
  i \leq k)}{ P(A 
  \mid \Open(\gamma_{i}) \text{ for all } i \leq k)} \leq \prod_{i \leq k}
  \frac{1 + Ce^{-a|\gamma_{i}|}}{1 - Ce^{-a|\gamma_{i}|}}
\end{equation}
and
\begin{equation} \label{E:Markov3}
  \frac{P(A \mid \Open(\gamma_{i}) \text{ for all } i \leq k)}
  {P(A \mid \Open(\gamma_{i}) \cap B_{i} \text{ for all } 
  i \leq k)} \leq \prod_{i \leq k}
  \frac{1 + Ce^{-a|\gamma_{i}|}}{1 - Ce^{-a|\gamma_{i}|}}.
\end{equation}

An event $A \subset \mathcal{B}(\ZZ)$ is called \emph{increasing}
if $\omega \in A$ and $\omega \leq \omega^{\prime}$ imply 
$\omega^{\prime} \in A$.  Here $\omega \leq \omega^{\prime}$
refers to the natural coordinatewise partial ordering.  A bond percolation
model $P$ has the \emph{FKG property} if $A, B$ increasing implies
$P(A \cap B) \geq P(A)P(B)$.

Throughout the paper, $\epsilon_{1}, \epsilon_{2},\ldots$, 
$c_{1}, c_{2},\ldots$ and $K_{1}, K_{2},\ldots$ are constants which
depend only on $P$.  We reserve $\epsilon_{i}$ for constants which 
are ``sufficiently small,'' $K_{i}$ for constants which are
``sufficiently large,'' and $c_{i}$ for those which fall in 
neither category.

Our basic assumptions will be that
\begin{align} \label{E:assump}
  &P \text{ is translation-invariant,
    invariant under $90^{\circ}$ rotation, and has the FKG} \\
  &\text{property, bounded
    energy and exponential decay of dual connectivity, and } \notag \\
  &P_{\Lambda,\rho} \text{ has the FKG property for all }
    \Lambda, \rho. \notag
\end{align}
When necessary we will also assume weak mixing, ratio weak mixing
and/or the near-Markov property for open circuits.

Since $P$ has the FKG property,
$-\log P(0^{*} \leftrightarrow x^{*})$ is a subadditive function
of $x$, and therefore the limit
\begin{equation} \label{E:surftens}
  \tau(x) = 
  \lim_{n \to \infty} -\frac{1}{n}\log
  P(0^{*} \leftrightarrow (nx)^{*}),
\end{equation}
exists for $x \in \mathbb{Q}^{2}$, provided we take the limit
through values of $n$ for which $nx \in \ZZ$.  This definition 
extends to $\mathbb{R}^{2}$ by continuity (see \cite{Al97}); the 
resulting $\tau$ is a norm on $\mathbb{R}^{2}$, when the dual
connectivity decays exponentially (i.e. $\tau(x)$ is positive for
all $x \neq 0$, or equivalently by lattice symmetry, $\tau(x)$ is
positive for some $x \neq 0$; we abbreviate this by saying
$\tau$ \emph{is positive}.)
By standard subadditivity results,
\begin{equation} \label{E:conupr}
  P(0^{*} \leftrightarrow x^{*}) \leq e^{-\tau(x)}
  \quad \text{for all} \ x.
\end{equation}
In the opposite direction, it is known \cite{Al97pwr}
that if $\tau$ is positive, 
ratio weak mixing holds and some milder assumptions 
hold then for some $\ep_{1}$ and $K_{1}$,
\begin{equation} \label{E:conlwr}
  P(0^{*} \leftrightarrow x^{*}) \geq
  \ep_{1}|x|^{-K_{1}}e^{-\tau(x)}
  \quad \text{for all} \ x \ne 0.
\end{equation}
It follows from the fact that the surface tension $\tau$
is a norm on $\mathcal{R}^{2}$ with axis symmetry that,
letting $e_{i}$ denote the $i$th unit coordinate vector, for
$\kt = \tau(e_{1})$ we have
\begin{equation} \label{E:equivnorm}
  \frac{1}{\sqrt{2}}\kt \leq
  \frac{\tau(x)}{|x|} \leq \sqrt{2}\kt \quad
  \text{ for all } \ x \ne 0.
\end{equation}
For a curve $\gamma$ tracing the boundary of a convex region
we define the $\tau$-length of $\gamma$ as the line integral
\[
  \mathcal{W}(\gamma) = \int_{\gamma} \tau(v_{x})\ dx,
\]
where $v_{x}$ is the unit forward tangent vector at $x$ and
$dx$ is arc length.
The \emph{Wulff shape} is the convex set $\mathcal{K}_{1} =
\mathcal{K}_{1}(\tau)$ which minimizes $\mathcal{W}(\partial V)$
subject to the constraint $|V| = 1$.  (We also refer to multiples 
of $\mathcal{K}_{1}$ as Wulff shapes, when confusion is 
unlikely.)  The Wulff shape actually minimizes $\mathcal{W}$
over a much larger class of $\gamma$ than just boundaries of 
convex sets (\cite{Ta1},\cite{Ta2}) but that fact will not concern us here.

Let $d_{\tau}(\cdot , \cdot)$ denote $\tau$-distance;
diam($\cdot$) and 
$\diam_{\tau}(\cdot)$ denote
Euclidean diameter and $\tau$-diameter, respectively.
$B(x,r)$ and $B_{\tau}(x,r)$
denote the closed Euclidean and $\tau$-balls,
respectively, of radius
$r$ about $x$.
We write $x + A$ for the translation of the set $A$ by the vector
$x$.  $d_{H}$ denotes Hausdorff distance.  

The deviation of a closed curve $\gamma$ from the boundary of an
area-$A$ Wulff shape is given by 
\[
  \Delta_{A}(\gamma) = \inf_{z} d_{H}(\gamma,z + 
    \partial(\sqrt{A}\mathcal{K}_{1})).
\]

As a convention, whenever we refer to the object in a finite 
class which
maximizes or minimizes something, we implicitly assume there is a 
deterministic algorithm for breaking ties.

Our description of heuristics for the local roughness is 
nonrigorous, so we permit ourselves the following partly vague
assumptions:
\begin{enumerate}
  \item[(i)] The Wulff shape boundary has curvature bounded away
    from 0 and $\infty$.
  \item[(ii)] For a droplet of any linear scale $l$, there is a
    characteristic length scale $\xi = \xi(l)$ representing the typical
    spacing between adjacent extreme points where the droplet 
    touches the boundary of its convex hull.
  \item[(iii)] On any length scale $n \leq \xi$ the fluctuations
    of the droplet boundary are of order $\sqrt{n}$.
\end{enumerate}
Here (iii) is reasonable because within each inward excursion 
between extreme points, the droplet boundary is nearly 
unconstrained, except by surface tension.
(i) is known for the Ising case from the exact solution (see
\cite{AVSZ}, \cite{MW}.)
Under (i), an arc of $\partial(l\mathcal{K}_{1})$ of length $n$
deviates from the corresponding secant line by 
a distance of order $n^{2}/l$, so 
we call $n^{2}/l$ the \emph{curvature deviation} (on scale $n$.)

On the characteristic scale $\xi$ the fluctuations and the
curvature deviation should be of the same order, that is,
\begin{equation} \label{E:sameorder}
  \frac{\xi^{2}}{l} \approx \sqrt{\xi}.
\end{equation}
To see this, consider two adjacent extreme points $x$ and $y$
of the droplet boundary separated by a distance of order $\xi$.  If
the curvature deviation
$\xi^{2}/l \gg \sqrt{\xi}$ this means the boundary between $x$ and $y$ 
is following the straight segment $\overline{xy}$ much more closely
than it follows the arc of $\partial(l\mathcal{K}_{1})$
from $x$ to $y$, that is, the droplet has an approximate facet from
$x$ to $y$.  But such facets are isoperimetrically disadvantageous 
since the Wulff shape lacks them under (i), 
so this is not a probable picture. 
Therefore we expect $\xi^{2}/l \leq \sqrt{\xi}$.  On the other hand, 
if the curvature deviation $\xi^{2}/l \ll \sqrt{\xi}$ then even on 
length scales $n \gg \xi$ an arc of the Wulff shape boundary looks 
nearly flat compared to the droplet boundary fluctuations, so along 
such an arc the roughly $n/\xi$ extreme points appear as a large 
number of local maxima of the droplet boundary above the Wulff
shape boundary, all approximately collinear.  This too is an unlikely 
picture, so we expect $\xi^{2}/l \geq \sqrt{\xi}$.

From (\ref{E:sameorder}) we get $\xi \approx l^{2/3}$, and then from 
(iii), we expect local roughness of order $l^{1/3}$.  The same relation 
(\ref{E:sameorder}) occurs in the assorted systems 
mentioned in the introduction.  

For $r > q > 0$, an $(q,r)$-\emph{bottleneck} 
in an exterior dual circuit 
$\gamma$ is an ordered pair $(u,v)$ of sites in $\gamma$ such that
there exists a path of length at most $q$ from $u$ to $v$ in
$\Int(\gamma)$, and the 
segments $\gamma^{[u,v]}$ and $\gamma^{[v,u]}$ each have 
diameter at least $r$.  When $r$ is not very large (as in our main
theorem, where $r$ can be of order $\log l$)
the absense of $(q,r)$-bottlenecks reflects a high degree of regularity 
in the structure of the boundary.

Our main theorem is the following.

\begin{theorem} \label{T:main}
  Let $P$ be a percolation model on $\mathcal{B}(\ZZ)$ satisfying
  (\ref{E:assump}), the near-Markov property for open circuits,
  and the ratio weak mixing property.  There exist
  $K_{i}, \epsilon_{i}$ such that for $A > K_{2}$ and $l = \sqrt{A}$,
  under the measure $P( \cdot \mid |\Int(\Gamma_{0})| \geq A)$ with 
  probability approaching 1 as $A \to \infty$ we have
  \begin{equation} \label{E:ALRmax}
    ALR(\Gamma_{0}) \leq K_{3}l^{1/3}(\log l)^{2/3},
  \end{equation}
  \begin{equation} \label{E:Hausmax}
    \Delta_{A}(\partial\Co(\Gamma_{0})) \leq K_{4}l^{2/3}(\log l)^{1/3},
  \end{equation}
  \begin{equation} \label{E:MLRmax}
    MLR(\Gamma_{0}) \leq K_{5}l^{2/3}(\log l)^{1/3},
  \end{equation}
  and, for $\epsilon_{2}A \geq r \geq 15q \geq K_{6}\log A$,
  \begin{equation} \label{E:bnkfree}
    \Gamma_{0} \text{ is } (q,r)-\text{bottleneck-free}.
  \end{equation} 
\end{theorem}

It is easy to see that if $MLR(\gamma)$ and 
$\Delta_{A}(\partial\Co(\gamma))$ are each a sufficiently
small fraction of $\sqrt{A}$, then $MLR(\gamma) = 
d_{H}(\gamma,\partial\Co(\gamma))$ and hence
\begin{equation} \label{E:DeltaA}
  \Delta_{A}(\gamma) \leq \Delta_{A}(\partial\Co(\gamma))
  + MLR(\gamma).
\end{equation}
(Here ``sufficiently small'' does not depend on $A$.)
Hence provided $A$ is large, (\ref{E:Hausmax}) and
(\ref{E:MLRmax}) imply 
\begin{equation} \label{E:Hausmax1}
  \Delta_{A}(\Gamma_{0}) \leq (K_{4} + 
  K_{5})l^{2/3}(\log l)^{1/3}.
\end{equation}

Theorem \ref{T:newmain} below shows that one
may condition on $|\Int(\Gamma_{0})| = A$ instead of on
$|\Int(\Gamma_{0})| \geq A$ in Theorem \ref{T:main}.

The basic strategy of the proof of 
Theorem \ref{T:main} is like that of \cite{ACC}, \cite{DKS} 
and \cite{IS97}:  one establishes a lower bound for 
$P(|\Int(\Gamma_{0})| 
\geq A)$ and upper bounds for the (unconditional) probability that 
$|\Int(\Gamma_{0})| \geq A$ but
$\Gamma_{0}$ does not have the desired behavior, these upper bounds 
being much smaller than the lower bound.  Both the upper and lower 
bounds involve coarse-graining to create a skeleton 
$(y_{0},..,y_{n},y_{0})$ for $\Gamma_{0}$, and both require showing
that the various connections $y_{j} \lra y_{j+1}$ occur 
approximately independently, that is,
\[
  P(y_{0} \lra ... \lra y_{n} \lra y_{0})
\]
can be approximated in an appropriate sense by
\[
  \prod_{j \leq n} P(y_{j} \lra y_{j+1})
\]
and hence by
\[
  \exp(- \sum_{j \leq n} \tau(y_{j+1} - y_{j}))
\]
(setting $y_{n+1} = y_{0}$.)
For the lower bound this is a straightforward application of the
FKG inequality, but for the upper bound the mixing properties 
established in Section \ref{S:prelim} are 
needed, and our methods must handle ``pathological'' forms of 
$\Gamma_{0}$ having many near-self-intersections.  The difficulties are
of two types.  First,
near-independence generally requires 
large spatial separation, or, under the 
near-Markov property, separation by a circuit of open bonds, neither of
which need be present in our context.  Second, direct
application of standard mixing properties such as 
weak mixing requires specifying in
advance some \emph{deterministic} spatially 
separated regions on which the near-independent
events will occur, but in our context
one does not know \emph{a priori} where the paths
$y_{j} \lra y_{j+1}$ may go.  This is where Lemma \ref{L:decouple}
below is important.

We consider now the special case of the FK model on $\mathcal{B}(\ZZ)$.
For each $(p,q)$ there is a value $p^{*}$ dual to $p$ at level $q$ 
given by 
\[
  \frac{(1 - p^{*})}{p^{*}} = \frac{p}{q(1 - p)};
\]
the dual configuration to the infinite-volume wired-boundary FK model at
$(p,q)$ is the infinite-volume free-boundary FK model at $(p^{*},q)$
(see \cite{Gr95}.)  The model has a percolation critical point
$p_{c}(q)$ which for $q = 1, q = 2$ and $q \geq 25.72$ is known to 
coincide with the self-dual point $p_{sd}(q) = \sqrt{q}/(1+\sqrt{q})$
\cite{LMR}; positivity
of $\tau$ is known to hold for $p > p_{sd}(q)$ for 
these same values of $q$.  For
$2 < q < 25.72$, it is known that positivity of $\tau$ holds for 
$p > p_{sd}(q-1)^{*}$, where the $*$ refers to duality at level $q$
\cite{Al97arc}.  The FK model without external field has the Markov 
property for open circuits; assuming positivity of $\tau$ it satisfies 
(\ref{E:assump}) 
(see \cite{Gr95}) and has the ratio weak mixing property \cite{Al97mix}.
In a forthcoming paper we will show that positivity of $\tau$
also implies ratio weak 
mixing for the FK model with external fields.  For now, we can conclude the
following from Theorem \ref{T:main}.

\begin{theorem} \label{T:FKcase}
  Let $P$ be the FK model at $(p,q)$ on $\mathbb{B}(\ZZ)$
  with $q \geq 1$ (without
  external fields) and suppose the surface tension
  $\tau$ is positive.
  There exists
  $K_{2}$ such that for $A > K_{2}$ and $l = \sqrt{A}$,
  under the measure $P( \cdot \mid |\Int(\Gamma_{0})| \geq A)$ with 
  probability approaching 1 as $A \to \infty$, (\ref{E:ALRmax})
  -- (\ref{E:bnkfree}) hold.
\end{theorem}

Following \cite{IS97} we call a dual circuit $\gamma \ 
r$-\emph{large} if $\diam_{\tau}(\gamma) > r$ and 
$r$-\emph{small} if $\diam_{\tau}(\gamma) \leq r$.

When $\tau$ is positive, in the box $\Lambda_{N} = [-N,N]^{2}$ 
the largest
open dual circuit typically has diameter of order $\log N$.  
Let $\mathfrak{C}_{N}$ denote the collection 
of all $(K \log N)$-large
exterior open dual circuits contained in $\Lambda_{N}$; here and
throughout this section, $K$ is a fixed ``sufficiently large''
constant.  To avoid any 
ambiguity, we impose a wired boundary condition on $\Lambda_{N}$.
Rather than conditioning on
$|\Int(\Gamma_{0})| \geq A$ as in Theorem \ref{T:main}, it is 
sometimes of interest to condition on the event
\[
  \sum_{\gamma \in \mathfrak{C}_{N}} |\Int(\gamma)| \geq A.
\]
In particular, it is natural 
to ask whether under such conditioning one has $|\mathfrak{C}_{N}|
= 1$ with high probability.  The alternative most difficult to 
rule out, as suggested by the form of the error term in Theorem
\ref{T:lowerbound}, is the presence of one or more ``small''
$(K \log N)$-large open dual circuits, of diameter of order 
$l^{1/3}(\log l)^{2/3}$ or less, outside a single large open 
dual circuit enclosing nearly area $A$.
In the context of phase separation in 
the Ising model, single-droplet theorems have been proved in
\cite{DKS} and \cite{IS97}.  However, the method there for ruling out 
``small'' $(K \log N)$-large droplets uses not purely surface 
tension but rather an argument that it is energentically 
preferable to dissolve such ``small'' droplets by spreading
the spins they contain throughout the bulk of the system, which
is allowed because there is no constraint on the total volume
contained in $(K \log N)$-large droplets.  Thus what we seek 
here is different, essentially a single-droplet 
theorem based purely on
surface tension considerations.  Let $P_{N,w}$ denote the
measure P conditioned on all bonds outside $\Int(\Lambda_{N})$ 
being open, that is, the measure under a wired 
boundary condition on $\mathcal{B}(\Int(\Lambda_{N}))$.

\begin{theorem} \label{T:single}
  Let $P$ be a percolation model on $\mathcal{B}(\ZZ)$ satisfying
  (\ref{E:assump}), the near-Markov property for open circuits,
  and the ratio weak mixing property.  There exist
  $\epsilon_{i},K_{i}$ such that for $N \geq 1, K_{7}(\log N)^{2} \leq A
  \leq c_{2}N^{2}$ and $l = \sqrt{A}$,
  under the measure $P_{N,w}( \cdot \mid 
  \sum_{\gamma \in \mathfrak{C}_{N}} |\Int(\gamma)| \geq A)$, with 
  probability approaching 1 (uniformly in $A$) 
  as $N \to \infty$ we have
  \begin{equation} \label{E:single}
    |\mathfrak{C}_{N}| = 1,
  \end{equation}
  and, for the unique open dual circuit $\gamma$ in 
  $\mathfrak{C}_{N}$,
  \begin{equation} \label{E:ALRmax2}
    ALR(\gamma) \leq K_{3}l^{1/3}(\log l)^{2/3},
  \end{equation}
  \begin{equation} \label{E:MLRmax2}
    MLR(\gamma) \leq K_{8}l^{2/3}(\log l)^{1/3},
  \end{equation}
  \begin{equation} \label{E:Hausmax2}
    \Delta_{A}(\partial\Co(\gamma)) \leq K_{4}l^{2/3}(\log l)^{1/3}
  \end{equation}
  and, for $\epsilon_{3}A \geq r \geq 15q \geq K_{9}\log A$,
  \begin{equation} \label{E:bnfree}
    \gamma \text{ is } (q,r)-\text{bottleneck-free}.
  \end{equation}
  Here $c_{2}$ is any constant less than 
  \[
    \sup\{c > 0: \sqrt{c}\mathcal{K}_{1} \subset [-1,1]^{2}\},
  \]
  and $K_{3},K_{4}$ are from Theorem \ref{T:main}.
\end{theorem}

As noted after Theorem \ref{T:main}, the FK model satisfies the
assumptions of Theorem \ref{T:single}, provided that $\tau$ is 
positive.

\section{Preliminaries---Coarse Graining and Mixing Properties} 
  \label{S:prelim}
We first define our coarse-graining concepts.
Our definition of the \emph{s-hull skeleton} follows
\cite{Al92}, with some added refinements.  For a contour
$\gamma$ let $E_{\gamma}$ denote the set of extreme
points of $\Co(\gamma)$ and let
$\gamma_{co}: [0,1] \to \mathbb{R}^{2}$
be a curve which traces $\partial\Co(\gamma)$ in the
direction of positive orientation, beginning at the
leftmost lattice site $u_{0}$ having minimal second
coordinate.  When confusion is unlikely we also use
$\gamma_{co}$ to denote the image of this curve.
Note that
\[
  u_{0} \in E_{\gamma} \subset \gamma_{co} \cap
  \gamma \cap
  (\mathbb{Z}^{2})^{*}.
\]
To define the $s$-hull skeleton we require that the $\tau$-diameter 
of $\gamma$ be at least $2s$.  We 
traverse $\gamma$ in the direction of positive orientation,
beginning at, say, the leftmost lattice site $u_{0}$ having minimal
second coordinate.
Given $u_{0},..,u_{j} \in E_{\gamma}$, go forward from
$u_{j}$ along $\gamma_{co}$ until either $u_{0}$ or
$\partial B_{\tau}(u_{j},s)$
is reached, at some point $u_{j+1}^{\prime}$.  Then
backtrack along $\gamma_{co}$ until a point of
$E_{\gamma}$ is reached (possibly $u_{j+1}^{\prime}$,
meaning we backtrack zero distance.)  If this backtracking
does not require
going all the way back to $u_{j}$, then this new point of
$E_{\gamma}$ is labeled $u_{j+1}$.  If instead the
backtracking does require going all the way back to
$u_{j}$, then from $u_{j+1}^{\prime}$ continue forward
along $\gamma_{co}$, necessarily in a straight line,
to the next point of $E_{\gamma}$, which is then
labeled $u_{j+1}$.  Stop the process when $u_{m+1} =
u_{0}$ for some $m$. The \emph{s-hull pre-skeleton} is
then $(u_{0},..,u_{m+1})$.  (A similar definition,
under the name ``$s$-hull skeleton,'' may be found in
\cite{Al92}.) The sites $u_{0},..,u_{m+1}$ are 
sites of $E_{\gamma}$ which
appear in order in $\gamma_{co}$ (and in
$\gamma$.)  Therefore
$\Co(\{u_{0},..,u_{m+1}\})$ is a convex polygon bounded by the
polygonal path $u_{0} \to .. \to u_{m+1}$.  Now
\[
  \mathcal{W}(\gamma_{co}) \geq
  \sum_{j=0}^{m} \tau(u_{j+1} - u_{j}),
\]
and as noted in \cite{Al92}, we have
\[
  \tau(u_{j+2} - u_{j}) > s \quad \text{for all}
  \quad 0 \leq j \leq m-2.
\]
Therefore
\begin{equation} \label{E:maxvert}
  m \leq 1 + 2\mathcal{W}(\gamma_{co})/s.
\end{equation}

To obtain the $s$-hull skeleton we refine the $s$-hull
pre-skeleton.  This is necessary because if $\gamma_{co}$ 
has a sharp corner, then the polygonal path
$u_{0} \to .. \to u_{m+1}$ may clip this corner excessively,
meaning part of $\gamma$ may be too far
outside the polygon for our needs.
We must add vertices so that, whenever possible,
the angular change between successive segments of the
polygonal path does not exceed $s/\text{diam}(\gamma)$,
which is of the order of the angular change we would obtain if
$\gamma$ were a circle.
By convexity, for each $x \in \gamma_{co}$ there exists a forward
tangent vector $v_{x}$ and a corresponding
forward tangent line; for $y \neq x \in \gamma_{co}$ let 
$\alpha(x,y)$ denote the angle measured counterclockwise from 
$v_{x}$ to $v_{y}$.  Fix $0 \leq j \leq m$ and let $u_{j0} = u_{j}$.
Note that $\alpha(u_{j},\cdot)$ is a nondecreasing function 
as one traces $\gamma_{co}$ from $u_{j}$
to $u_{j+1}$.
Having defined $u_{j0},..,
u_{jk} \in E_{\gamma} \cap
\gamma_{co}^{[u_{j-1},u_{j+1}]}$, let $u_{j,k+1}$ be the first
point of $\gamma_{co}$ after $u_{j,k}$ for which
$\alpha(u_{j,k},u_{j,k+1}) \geq s/\diam_{\gamma}$.  If there is no
such point, then set $u_{j,k+1} = u_{j+1}$ and stop the process.
Necessarily
$u_{j,k+1}$ is a lattice site in $E_{\gamma}$.

We call the sites $u_{jk}$ strictly between $u_{j}$
and $u_{j+1}$ \emph{refinement sites}.
Let $(w_{0},..,w_{n},$ $w_{n+1})$ be a relabeling of
all sites $u_{jk}, 0 \leq j \leq m, k \geq 1$, in order on
$\gamma_{co}$, with $w_{0} = w_{n+1} = u_{0}$;
we call $(w_{0},..,w_{n},w_{n+1})$ the
\emph{s-hull skeleton}
of $\gamma$ and denote it $\HSkel_{s}(\gamma)$.
The polygonal path $w_{0} \to .. \to w_{n+1}$ is denoted
$\HPath_{s}(\gamma)$.
It is easy to see
that the angle between $u_{jk} - u_{j,k-1}$ and
$u_{j,k+2} - u_{jk}$ cannot be less than
$s/\diam(\gamma)$, for $k \geq 1$.  It follows
that the number of refinement sites satisfies
\[
  n - m \leq 4\pi \text{diam}(\gamma)/s
\]
With (\ref{E:maxvert})
this shows that the number of sites in the $s$-hull
skeleton satisfies
\begin{equation} \label{E:maxvert2}
  n + 1 \leq K_{10}\text{diam}(\gamma)/s.
\end{equation}
As with the pre-skeleton, the sites $w_{0},..,w_{n+1}$
appear in order in $\gamma$ and in $\gamma_{co}$,
and $\Co(\{w_{0},..,w_{n+1}\})$ is a convex polygon bounded by
$\HPath_{s}(\gamma)$.

A key property of the $s$-hull skeleton involves the
extent to which $\gamma$ can go outside
$\Co(\{w_{0},..,w_{n+1}\})$. Let
$T_{j}$ denote the triangle formed by the segment
$\overline{w_{j}w_{j+1}}$,
the forward tangent line at $w_{j}$
and the backward tangent line at $w_{j+1}$.  The angle between 
these two tangent lines is at most $s/\diam(\gamma)$.  Now
every point of $\gamma$ outside $\Co(\{w_{0},..,w_{n+1}\})$ is in
some $T_{j}$, and $\gamma_{co} \cap T_{j} \subset
\gamma^{[w_{j},w_{j+1}]}$. Let
\[
  J = \{j \leq n: \tau(w_{j+1} - w_{j}) \leq 2s\}.
\]
For distinct points $x, y \in \mathbb{R}^{2}$ let $H_{xy}$
($\overline{H}_{xy}$)
denote the open (closed)
halfspace which is to the right of the line from
$x$ to $y$.
From the construction of the $s$-hull pre-skeleton,
if $\gamma^{(u_{j},u_{j+1})} \not\subset B_{\tau}(u_{j},s)$
for some $j$, then $\gamma \cap H_{u_{j}u_{j+1}} = \phi$. It follows
that if $\tau(w_{j+1} - w_{j}) > 2s$ for some $j$, then $w_{j}$
and $w_{j+1}$ are sites of the $s$-hull pre-skeleton and
$\gamma \cap H_{w_{j}w_{j+1}} = \phi$.  Thus we have
\begin{equation} \label{E:outside}
  \Int(\gamma) \backslash \Int(\HPath_{s}(\gamma)) \subset
  \cup_{j \in J} T_{j}.
\end{equation}
For $j \in J$ we have $|T_{j}| \leq K_{11}s^{3}/
\text{diam}(\gamma)$ and
\[
  d\bigl(x, \Int(\HPath_{s}(\gamma))\bigr) \leq
  K_{12}s^{2}/\text{diam}(\gamma)
  \quad \text{for all } x \in T_{j}.
\]
With (\ref{E:outside}) this shows that
\begin{equation} \label{E:outside2}
  |\Int(\gamma) \backslash \Int(\HPath_{s}(\gamma))| \leq
  K_{13}s^{2}
\end{equation}
and
\begin{equation} \label{E:outside3}
  \sup_{x \in \Co(\gamma)} d\bigl(x,\Int(\HPath_{s}(\gamma))\bigr) \leq
  K_{12}s^{2}/\text{diam}(\gamma).
\end{equation}
From (\ref{E:outside3}) and convexity it follows that
\begin{equation} \label{E:extralength}
  \mathcal{W}(\gamma_{co}) \leq \mathcal{W}(\HPath_{s}(\gamma))
  + K_{14}s^{2}/\diam(\gamma).
\end{equation}

We turn now to mixing properties.
The following is an immediate consequence of the definition of
ratio weak mixing.

\begin{lemma} \label{L:ratiowm}
(\cite{Al97pwr})
Suppose $P$ has the ratio weak mixing property.
There exists a constant $K_{15}$ as follows.
Suppose $r > 3$ and $\mathcal{D}, \mathcal{E} \subset \mathbb{Z}^{2}$
with $\diam(\mathcal{E}) \leq r$ and $d(\mathcal{D},\mathcal{E}) \geq
K_{15}\log r$.  Then for all $A \in \mathcal{G}_{\mathcal{D}}$ and
$B \in \mathcal{G}_{\mathcal{E}}$, we have
\[
  \frac{1}{2}P(A)P(B) \leq P(A \cap B) \leq 2P(A)P(B).
\]
\end{lemma}

A weakness of Lemma \ref{L:ratiowm} is that the locations
$\mathcal{D}, \mathcal{E}$ of the two events must be deterministic.  
The next lemma applies only to a limited class of events
but allows the locations to be partially random.
For $\mathcal{C} \subset \mathcal{D}
\subset \mathcal{B}(\ZZ)$ we say an event $A 
\subset \{0,1\}^{\mathcal{D}}$ 
\emph{occurs on} $\mathcal{C}$ (or on $\mathcal{C}^{*}$)
in $\omega \in \{0,1\}^{\mathcal{D}}$ 
if $\omega^{\prime} \in A$  for every
$\omega^{\prime}\in \{0,1\}^{\mathcal{D}}$  satisfying 
$\omega^{\prime}_{e} = \omega_{e}$
for all $e \in \mathcal{C}$.  
For a possibly random set $\mathcal{F}(\omega)$ we
say $A$ \emph{occurs only on} $\mathcal{F}$ (or equivalently,
on $\mathcal{F}^{*}$) if $\omega \in A$ implies
$A$ occurs on $\mathcal{F}(\omega)$ in $\omega$.
We say events $A$ and $B$ 
\emph{occur at separation} $r$ in $\omega$ if there exist 
$\mathcal{C}, \mathcal{E} \subset \mathcal{D}$ with 
$d(\mathcal{C},\mathcal{E}) \geq r$ such
that $A$ occurs on $\mathcal{C}$ and $B$ occurs on
$\mathcal{E}$ in $\omega$.  
Let $A \circ_{r} B$ denote the event that $A$ and $B$ 
occur at separation $r$.
Let $\mathcal{D}^{r} = \{e \in \mathcal{B}(\ZZ): 
d(e,\mathcal{D}) \leq r\}$.

In the next lemma we give two alternate hypotheses, in the
interest of wider applicability, though either hypothesis alone
suffices for our purposes in this paper.

\begin{lemma} \label{L:decouple}
  Assume (\ref{E:assump}) and either (i) the ratio weak mixing property
  or (ii) both the weak mixing property and the near-Markov property
  for open circuits.  
  There exist constants $K_{i},\epsilon_{i}$ as follows.
  Let $\mathcal{D} \subset \mathcal{B}(\ZZ), x^{*} \in 
  (\ZZ)^{*}$ and 
  $r > K_{16} \log |\mathcal{D}|$, and let
  $A, B$ be
  events such that $A$ occurs only on $C_{x^{*}}$ and
  $B \in \mathcal{G}_{\mathcal{D}}$.  Then 
  \begin{equation} \label{E:decouple}
    P(A \circ_{r} B) \leq
    (1 + K_{17}e^{-\epsilon_{5}r}) P(A)P(B).
  \end{equation}
\end{lemma}
\begin{proof}  
First suppose $P$ has the ratio weak mixing property.
For $y^{*} \in V^{*}(\mathcal{D}^{r})$, let $B^{y^{*}} = 
B(y^{*},r/12),\ B^{y^{*}}_{\mathbb{Z}} = B^{y^{*}} \cap (\ZZ)^{*}$ and
$D^{y^{*}} = B(y^{*},r/6),\ D^{y^{*}}_{\mathbb{Z}} = 
D^{y^{*}} \cap (\ZZ)^{*}$.
Define events 
\[
  G_{y^{*}} = [y^{*} \leftrightarrow \partial_{in} B^{y^{*}}_{\mathbb{Z}} 
  \text{  by an open dual path}],
\]
\[
  Z_{y^{*}} = [A \text{ and } B 
  \text{ occur at separation } r \text{ on }
  \omB(D^{y^{*}})^{c}],
\]
\[
  Q = \cup_{y^{*} \in V^{*}(\mathcal{D}^{*})} 
  (Z_{y^{*}} \cap G_{y^{*}}).
\]
Then by (\ref{E:conupr}),
\begin{equation} \label{E:Gbound}
  P(G_{y^{*}})  \leq K_{18}e^{-\epsilon_{6}r}.
\end{equation}
Let $C$ and $\lambda$ be as in (\ref{E:weakmix}).  Then for some $K_{i},
\epsilon_{i}$
depending on $C, \lambda$, provided $K_{16}$ is chosen large enough 
\begin{align} \label{E:Qsmall}
  P&((A \circ_{r} B) \cap Q) \\
  &\leq \sum_{y^{*} \in V^{*}(\mathcal{D}^{r})} P(Z_{y^{*}})
    P(G_{y^{*}} \mid Z_{y^{*}}) \notag \\
  &\leq \sum_{y^{*} \in V^{*}(\mathcal{D}^{r})} P(Z_{y^{*}})
    (P(G_{y^{*}}) + K_{19}r^{2}e^{-\epsilon_{7}r}) \notag \\
  &\leq K_{20}r^{2}|\mathcal{D}| P(A \circ_{r} B) 
    (K_{18}e^{-\epsilon_{6}r} + K_{19}r^{2}e^{-\epsilon_{7}r}) \notag \\
  &\leq K_{21}e^{-\epsilon_{8}r} P(A \circ_{r} B). \notag
\end{align}
Therefore
\begin{equation} \label{E:ridQ}
  P(A \circ_{r} B) \leq (1 -  K_{21}e^{-\epsilon_{8}r})^{-1}
  P((A \circ_{r} B) \cap Q^{c}).
\end{equation}

Suppose $\omega \in (A \circ_{r} B) \cap Q^{c}$, and consider
$\mathcal{C}, \mathcal{E}$ with $\mathcal{C}^{*} \subset 
C_{x^{*}}(\omega),\ \mathcal{E} \subset \mathcal{D}$  
and $d(\mathcal{C},\mathcal{E}) \geq r$
for which $A$ occurs on $\mathcal{C}^{*}$ and $B$ occurs on 
$\mathcal{E}$ in $\omega$.  If $r/2 < d(y^{*},
\mathcal{C}) \leq r/2 + 1$ for 
some $y^{*} \in V^{*}(\mathcal{D}^{r})$, then $y^{*} \notin 
C_{x^{*}}(\omega)$
since $\omega \notin Z_{y^{*}} \cap G_{y^{*}}$.  Therefore
$C_{x^{*}}(\omega) \subset \mathcal{C}^{r/2}$ and hence
$d(\mathcal{E},C_{x^{*}}(\omega)) > r/2$.  
For $\mathcal{F} \subset \mathcal{B}(\ZZ)$ let $B_{\mathcal{F},r}$ 
be the event that $B$ occurs on some set $\mathcal{E} \subset
\mathcal{D}$ with
$d(\mathcal{E},\mathcal{F}) > r/2$.  We call $\mathcal{F} \ 
A$-\emph{sufficient} if $\omega_{e} = 0$
for all $e \in \mathcal{F}$ implies $\omega \in A$.
Let $\mathfrak{A}$ denote the set of all finite $A$-sufficient 
subsets of $\mathcal{B}(\ZZ)$.  Since
the event $[C_{x^{*}} = \mathcal{F}^{*}]$ occurs on 
$\omB(V(\mathcal{F}))$ for all $\mathcal{F} \subset \mathcal{B}(\ZZ)$, we have 
by ratio weak mixing
\begin{align} \label{E:noQ}
  P&((A \circ_{r} B) \cap Q^{c}) \\
  &\leq \sum_{\mathcal{F} \in \mathfrak{A}} 
    P([C_{x^{*}} = \mathcal{F}^{*}] \cap B_{\mathcal{F},r})
    \notag \\
  &\leq \sum_{\mathcal{F} \in \mathfrak{A}}
    (1 + K_{22}|\mathcal{D}|e^{-\epsilon_{9}r})
    P(C_{x^{*}} = \mathcal{F}^{*})P(B_{\mathcal{F},r})  \notag \\
  &\leq (1 + K_{23}e^{-\epsilon_{10}r})P(A)P(B). \notag
\end{align}
Combining (\ref{E:ridQ}) and (\ref{E:noQ}) proves
(\ref{E:decouple}), provided 
$K_{16}$ is sufficiently large.

Under hypothesis (ii) of 
weak mixing and the near-Markov property for open circuits,
the proof through the first inequality of (\ref{E:noQ}) is still valid, but
we need to modify the rest of
(\ref{E:noQ}) as follows.  Fix $\mathcal{F} \in 
\mathfrak{A}$ and let $F
= [C_{x^{*}} = \mathcal{F}^{*}]$,
$\mathcal{D}^{\prime} = \mathcal{D} \backslash \mathcal{F}^{r}$,
$\mathcal{A} =
\{e \in \mathcal{B}(\mathbb{Z}^{2}): r/6 \leq
d(e,\mathcal{D}^{\prime}) \leq r/3\}$.  Let $\mathfrak{C}$ be the set
of all circuits (of regular bonds) 
in $\mathcal{A}$ which surround
$\mathcal{D}^{\prime}$ and let $O_{\mathcal{A}}$ be 
the event that some circuit $\alpha \in \mathfrak{C}$ is open; 
for $\omega \in O_{\mathcal{A}}$  there is a
unique outermost open circuit in $\mathfrak{C}$, which 
we denote $\Gamma = \Gamma(\omega)$.  Note that 
\begin{equation} \label{E:notOA}
  O_{\mathcal{A}}^{c} \subset \cup_{y^{*} \in V^{*}(\mathcal{A})}
  \ G_{y^{*}},
\end{equation}
\begin{equation} \label{E:alphamin}
  |\alpha| \geq r \quad \text{for all } \alpha \in \mathfrak{C},
\end{equation}
and 
\begin{equation} \label{E:outermost}
  [\Gamma = \alpha] = \Open(\alpha) \cap G_{\alpha} \quad
  \text{for some event } G_{\alpha} \in 
  \mathcal{G}_{\mathcal{B}(\Ext(\alpha))}.
\end{equation}
Let
\[
  p_{\alpha} = P(\Gamma = \alpha),
    \qquad p_{F\alpha} = P(\Gamma
    = \alpha \mid F), \qquad 
    p_{FB\alpha} 
    = P(\Gamma = \alpha 
    \mid F \cap B_{\mathcal{F},r}).
\]
By weak mixing we have for $\delta = K_{24}e^{-\epsilon_{11}r}$:
\begin{equation} \label{E:totalvar}
  \sum_{\alpha \in \mathfrak{C}} |p_{F\alpha} - p_{\alpha}| 
  \leq \delta, \qquad
  \sum_{\alpha \in \mathfrak{C}} |p_{FB\alpha } - p_{\alpha}| 
  \leq \delta.
\end{equation}
Define the set of ``good'' circuits
\[
  \mathfrak{R}
  = \{\alpha \in \mathfrak{C}: p_{F\alpha} \leq 
    p_{\alpha}(1 + \sqrt{\delta})\}
\]
and let
\[
  h(\alpha) = \left(\frac{p_{F\alpha}}{p_{\alpha}} - 1\right) 1_{[\alpha \in 
  \mathfrak{C} \backslash \mathfrak{R}]}.
\]
From (\ref{E:alphamin}), (\ref{E:outermost}) and the near-Markov property 
for open circuits, if $r$ is sufficiently
large we obtain
\begin{equation} \label{E:Bchange}
  P(B_{\mathcal{F},r} \mid F \cap [\Gamma = \alpha]) \leq
    (1 + \delta^{\prime})P(B_{\mathcal{F},r} \mid \Gamma = \alpha),
\end{equation}
where $\delta^{\prime} = 3Ce^{-ar} < 1$ for $C, a$ as in 
(\ref{E:Markov}).

We need to bound $P(F \cap B_{\mathcal{F},r})$.  To do this, we 
decompose $F \cap B_{\mathcal{F},r}$ into 3 pieces by intersecting it with 
$[\Gamma \in \mathfrak{R}], [\Gamma \in \mathfrak{C} \backslash 
\mathfrak{R}]$ and $O_{\mathcal{A}}^{c}$ (the latter meaning there
is no $\Gamma$.)  We then show that the first piece is approximately 
bounded by $P(F)P(B_{\mathcal{F},r})$, and
the other two pieces are negligible relative 
to the size of the full event.  Specifically, from (\ref{E:Bchange}),
\begin{align} \label{E:goodalpha}
  P(&F \cap B_{\mathcal{F},r} \cap [\Gamma \in \mathfrak{R}]) \\
  &= \sum_{\alpha \in \mathfrak{R}} p_{F\alpha}P(F)P(B_{\mathcal{F},r} 
    \mid F \cap [\Gamma = \alpha]) \notag \\
  &\leq (1 + \sqrt{\delta})(1 + \delta^{\prime})P(F) \sum_{\alpha \in 
    \mathfrak{R}}
    p_{\alpha}P(B_{\mathcal{F},r} \mid \Gamma = \alpha) \notag \\
  &\leq (1 + 2\sqrt{\delta} + \delta^{\prime})P(F) P(B_{\mathcal{F},r}).
    \notag
\end{align}
Next, one application of (\ref{E:totalvar}) yields
\[ 
  E(h(\Gamma)1_{O_{\mathcal{A}}}) = \sum_{\alpha 
    \in \mathfrak{C} \backslash \mathfrak{R}} 
    h(\alpha)p_{\alpha} \leq \delta;
\]
this and a second application of (\ref{E:totalvar}), with Markov's 
inequality, yield
\begin{align} \label{E:badalpha}
  P(&F \cap B_{\mathcal{F},r} \cap 
    [\Gamma \in \mathfrak{C} \backslash \mathfrak{R}]) \\
  &= P(F \cap B_{\mathcal{F},r}) 
    P(h(\Gamma) > \sqrt{\delta} \mid 
    F \cap B_{\mathcal{F},r}) \notag \\
  &\leq P(F \cap B_{\mathcal{F},r}) 
    (P(h(\Gamma) > \sqrt{\delta}) 
    + \delta)  \notag \\
  &\leq P(F \cap B_{\mathcal{F},r}) (\sqrt{\delta}
    + \delta).  \notag 
\end{align}
Finally, similarly to (\ref{E:Qsmall}) we have using (\ref{E:notOA})
\begin{align} \label{E:noalpha}
 P(F \cap B_{\mathcal{F},r} \cap O_{\mathcal{A}}^{c})
    &\leq \sum_{y^{*} \in V^{*}(\mathcal{A})}
      P(F \cap B_{\mathcal{F},r} \cap G_{y^{*}}) \\
  &\leq \delta^{\prime \prime}P(F \cap B_{\mathcal{F},r}),
    \notag
\end{align}
where $\delta^{\prime \prime} = K_{25}e^{-\epsilon_{12}r}$.
Combining (\ref{E:goodalpha}), (\ref{E:badalpha}) and (\ref{E:noalpha})
we obtain
\begin{align} \label{E:FBbound}
  P(F \cap B_{\mathcal{F},r}) &\leq 
    (1 - \delta - \sqrt{\delta} - \delta^{\prime \prime})^{-1}
    P(F \cap B_{\mathcal{F},r} \cap [\Gamma \in \mathfrak{R}]) \\
  &\leq (1 + K_{26}e^{-\epsilon_{13}r})P(F)P(B_{\mathcal{F},r}). \notag
\end{align}
Summing over $\mathcal{F}$ yields
\[
  \sum_{\mathcal{F} \in \mathfrak{A}} 
  P([C_{x^{*}} = \mathcal{F}^{*}] \cap B_{\mathcal{F},r} )
  \leq (1 + K_{26}e^{-\epsilon_{13}r}) P(A)P(B) 
\]
which substitutes for (\ref{E:noQ}).
\end{proof}

\section{Lower Bounds for Open Dual Circuit Probabilities} \label{S:lower}
In this section we prove the following result.

\begin{theorem} \label{T:lowerbound}
  Let $P$ be a percolation model on $\mathcal{B}(\ZZ)$ satisfying
  (\ref{E:assump}), the near-Markov property for open circuits,
  positivity of $\tau$
  and the ratio weak mixing property.  There exist
  $K_{i}$ such that for $A > K_{27}$ and $l = \sqrt{A}$,
  \[
    P(|\Int(\Gamma_{0})| \geq A) \geq \exp(-w_{1}\sqrt{A} - K_{28}l^{1/3}
    (\log l)^{2/3}).
  \]
\end{theorem}

The size of the error term $K_{28}l^{1/3}(\log l)^{2/3}$ in 
Theorem \ref{T:lowerbound} is important because it determines what 
``bad'' behaviors can be ruled out as unlikely---in particular,
those which have
probability at most $\exp(-w_{1}l - cl^{1/3}(\log l)^{2/3}))$ 
for some $c > K_{28}$.  Though our error term is likely not 
optimal---according to \cite{Hr} the optimal error term may be
of order $\log l$---it is enough of an improvement over the 
corresponding results in \cite{DKS} and \cite {IS97} to enable us to 
establish an apparently near-optimal bound on the local roughness.

The proof of our Theorem \ref{T:lowerbound}
relies on the following result, the halfspace version of (\ref{E:conlwr}).

\begin{theorem} \label{T:halfspace}
  (\cite{Al97pwr})
  Let $P$ be a percolation model on $\mathcal{B}(\ZZ)$ satisfying
  (\ref{E:assump}), positivity of
  $\tau$ and the ratio weak mixing property.  There exist
  $\epsilon_{14}, K_{29}$ such that for all $x \neq y \in \mathbb{R}^{2}$
  and all dual sites $u, v \in H_{xy}$, 
  \[
    P(u \lra v \text{ via an open dual path in } 
    H_{xy}) \geq \frac{\epsilon_{14}}{|x|^{K_{29}}}e^{-\tau(v-u)}
  \]
\end{theorem}

\begin{proof}[Proof of Theorem \ref{T:lowerbound}]
Let $s = l^{2/3}(\log l)^{1/3}$ and $\delta = 
K_{30}s^{2}/l$, with $K_{30}$
to be  specified.  Let $(y_{0},..,y_{n},y_{0})$ be the $s$-hull
skeleton of $\partial(l + \delta)\mathcal{K}_{1}$. 
For each $i$ let $y_{i}^{\prime}$ be a dual site with
$y_{i}^{\prime} \in H_{y_{i-1}y_{i}} \cap H_{y_{i}y_{i+1}}$ and
$|y_{i}^{\prime} - y_{i}| \leq 2\sqrt{2}$. 
By (\ref{E:outside3}), provided $K_{30}$ is large enough we have
\[
  \Co(\{y_{0},..,y_{n}\}) \supset l\mathcal{K}_{1} \quad
  \text{and hence} \quad |\Co(\{y_{0},..,y_{n}\})| \geq A.
\]
Further,
\begin{equation} \label{E:lengthhull}
  \sum_{j=0}^{n} \tau(y_{j+1}^{\prime} - y_{j}^{\prime}) \leq 
  \mathcal{W}(\partial(l + \delta)\mathcal{K}_{1}) + 4\kt n
  \leq w_{1}l + K_{31}l^{1/3}(\log l)^{2/3}
\end{equation}
(with $y_{n+1} = y_{0}$.)
Therefore using the FKG property, Theorem \ref{T:halfspace},
(\ref{E:lengthhull}) and (\ref{E:maxvert2}),
\begin{align} \label{E:lowerbd}
  P(|\Int(\Gamma_{0})| \geq A) &\geq P(\Gamma_{0}
    \text{ encloses } \Co(\{y_{0},..,y_{n}\})) \\
  &\geq P(y_{j}^{\prime} \lra y_{j+1}^{\prime} \text{ via a path in }
    H_{y_{j}y_{j+1}} \text{ for all } j \leq n) \notag \\
  &\geq \prod_{j=0}^{n} P(y_{j}^{\prime} \lra y_{j+1}^{\prime} \text{ via a path in }
    H_{y_{j}y_{j+1}}) \notag \\
  &\geq \left(\frac{\epsilon_{14}}{l^{K_{29}}}\right)^{(n+1)}
    \exp\left(-\sum_{j=0}^{n} \tau(y_{j+1}^{\prime} 
    - y_{j}^{\prime})\right) \notag \\
  &\geq \exp(-w_{1}l - K_{32}l^{1/3}(\log l)^{2/3}). \notag
\end{align}
\end{proof}

\section{Upper Bounds for Open Dual Circuit Probabilities} \label{S:upper}  
We need to develop a method of cutting a dual circuit
across a bottleneck, modifying the bond configuration to create two dual 
circuits.  The cutting procedure is simplified if the bottleneck is clean, in
the following sense.
The \emph{canonical path} from dual site $u = (x_{1},y_{1})$ 
to dual site $v = (x_{2},y_{2})$ is the path, denoted $\zeta_{uv}$, 
which goes horizontally from 
$(x_{1},y_{1})$ to $(x_{2},y_{1})$, then vertically to 
$(x_{2},y_{2})$.  We call a bottleneck $(u,v)$ \emph{clean} if 
$\zeta_{uv} \subset \Int(\gamma)$ (except for the endpoints $u, v$.)
The next lemma will enable us to restrict our cutting to clean bottlenecks.

\begin{lemma} \label{L:bottleneck}
  If a dual circuit $\gamma$ contains a $(q,r)$-bottleneck for some
  $r > 3q > 0$, then $\gamma$ contains a clean $(q,r/3)$-bottleneck.
\end{lemma}
\begin{proof}
Suppose $\gamma$ contains a $(q,r)$-bottleneck $(u,v)$.
We have two disjoint paths from $u$ to $v$: $\gamma^{[u,v]}$ 
and $\gamma^{[v,u]}$ (traversed backwards.)
Each of these may intersect $\zeta_{uv}$ a number of times.
Accordingly, $\zeta_{uv}$ contains a finite sequence of sites
$u = x_{0},x_{1},..,x_{m} = v$ such that the segment $\zeta_{i}$ of
$\zeta_{uv}$ between $x_{i-1}$ and $x_{i}$ satisfies $\zeta_{i}
\subset \Int(\gamma)$ for all $i \in I$ and 
$\zeta_{i} \subset \Ext(\gamma) \cup \gamma$  for all 
$i \notin I$, where $I$ consists
either of all odd $i$ or of all even $i$.
For $i \in I$, we call the segment of $\zeta_{uv}$
with endpoints $x_{i-1}$ and
$x_{i}$ an \emph{interior gap}.
Let $\psi$ be a dual path from $u$ to $v$ in $\Int(\gamma)$ with
$|\psi| \leq q$.  We can extend $\psi$ to a doubly infinite path
$\psi^{+}$ by adding on (possibly non-lattice) paths $\psi_{1}$ from $v$
to $\infty$ and $\psi_{2}$ from $\infty$ to $u$, both in $\Ext(\gamma)$.
The path $\psi^{+}$ divides the plane into two regions,
$A_{L} \supset \gamma^{[v,u]}$ to the left of
$\psi^{+}$ and $A_{R} \supset \gamma^{[u,v]}$ to the right.
Replacing $\psi$ with $\zeta_{uv}$ in the definition of $\psi^{+}$, 
we obtain another doubly 
infinite path $\zeta^{+}$.  The path $\zeta^{+}$ is
not necessarily self-avoiding, but $\mathbb{R}^{2} \backslash
\zeta^{+}$ has exactly two unbounded components $B_{L}$ and
$B_{R}$, to the left and right of $\zeta^{+}$, respectively.
Since $\diam(\zeta_{uv}) \leq q$, there exist sites $z_{1} \in \gamma^{[u,v]}$
and $z_{2} \in \gamma^{[v,u]}$ for which $d(z_{j},\zeta_{uv}) \geq
(r-q)/2 > q \ (j = 1,2)$ and hence $z_{1} \in B_{R}, z_{2} \in B_{L}$.
Let $\theta$ be a (possibly non-lattice)
path from $z_{1}$ to $z_{2}$ in $\Int(\gamma)$.
Then $\theta$ must intersect $\zeta^{+}$, and hence must intersect 
$\zeta_{uv}$,
necessarily in some interior gap.  Thus every $\theta$ from $z_{1}$
to $z_{2}$  in $\Int(\gamma)$ must cross at least one interior gap, so 
there exists an interior gap $\zeta_{i}$ which separates $z_{1}$ 
and $z_{2}$, that is, exactly one of $z_{1}, z_{2}$ is in
$\gamma^{[x_{i-1},x_{i}]}$.  It follows 
that $(x_{i-1},x_{i})$ is a clean $(q,(r-q)/2)$-bottleneck.  Since 
$(r-q)/2 > r/3$, the proof is complete.
\end{proof}

Define $R_{x} = x + [-1/2,1/2]^{2}$ and
$R_{x}^{+} = x + [-1,1]^{2}$.  Let $Q_{1}(u,v) = \cup_{x \in \zeta_{uv}}
R_{x}$ and $Q_{2}(u,v) = \cup_{x \in \zeta_{uv}}
R_{x}^{+}$ .  Note that
\[
  |\partial Q_{2}(u,v)| \leq 4q + 8.
\]
Let $J_{1}(u,v),J_{2}(u,v),...$ be an enumeration of the subsets of
$\partial Q_{2}(u,v)$.  We say a clean $(q,r)$-bottleneck
$(u,v)$ in a dual circuit
$\gamma$ is of \emph{type} $\eta$
if the set of bonds in $\partial Q_{2}(u,v)$ which are contained (except
possibly for endpoints) in $\Int(\gamma)$ is precisely $J_{\eta}(u,v)$.

We assume we have a fixed but arbitrary algorithm for choosing a
particular $(q,r)$-bottleneck, which we then call \emph{primary},
from any circuit containing one or more $(q,r)$-bottlenecks.
When a configuration $\omega$ includes an exterior dual circuit $\gamma$
for which $(u,v)$ is
a primary $(q,r)$-bottleneck of type $\eta$, we can apply a procedure, which
we term \emph{bottleneck surgery} (on $\gamma$, at
$(u,v)$) to create a new configuration denoted
$Y_{uv\eta}(\omega)$.  Bottleneck surgery consists
of replacing the configuration $\omega$ with the configuration
given by
\begin{equation} \label{E:Yomega}
  Y_{uv\eta}(\omega)_{e} =
  \begin{cases}
    1, &\text{if $e \in \partial Q_{1}(u,v)$;} \\
    0, &\text{if $e^{*} \in J_{\eta}(u,v)$;} \\
     \omega_{e}, &\text{otherwise,}
  \end{cases}
\end{equation}
for each bond $e$.
The configuration $Y_{uv\eta}(\omega)$ 
then contains two or more disjoint
open dual circuits $\alpha_{i}$, each consisting of some 
dual bonds of $\gamma$ and 
some dual bonds of $J_{\eta}(u,v)$, with no open dual path
connecting $\alpha_{i}$ to $\alpha_{j}$ for $i \neq j$, and with
\begin{equation} \label{E:newcircuits}
  \cup_{i} \Int(\alpha_{i}) = \Int(\gamma) \backslash Q_{2}(u,v).
\end{equation}
Further,
\begin{equation} \label{E:smalloff}
  |Q_{2}(u,v)| +
  \sum_{i: \alpha_{i} (\kt r/3)-\text{small}}
  |\Int(\alpha_{i})| \leq K_{33}r^{2},
\end{equation}
and, since $\gamma$ is exterior,
there is no open dual path connecting $\alpha_{i}$ to 
$\alpha_{j}$ for $i \neq j$.
We call each $\alpha_{i}$ an $(q,r)$-\emph{offspring} 
or an $(q,r)$-\emph{descendant} of $\gamma$.  A $(q,r)$-offspring
of an $(q,r)$-descendant is also an $(q,r)$-descendant, iteratively.  We 
may perform bottleneck surgery on each $(q,r)$-offspring of $\gamma$ 
which contains a clean $(q,r)$-bottleneck, and iterate this process until 
no descendant of $\gamma$ contains such a clean $(q,r)$-bottleneck
(necessarily after a finite
number of surgeries.)  The
bottleneck-free $(q,r)$-descendants are called \emph{final} 
$(q,r)$-\emph{descendants}.
Among final $(q,r)$-descendants, the one enclosing maximal area
is called the
\emph{maximal} $(q,r)$-\emph{descendant} of $\gamma$
and denoted $\alpha_{\max,\gamma}$. The set of all $(\kt r/3)$-large 
final $(q,r)$-descendants of $\gamma$ is denoted
$\mathfrak{F}_{(q,r)}(\gamma)$; the non-maximal among these
form the set $\mathfrak{F}_{(q,r)}^{\prime}(\gamma) = 
\mathfrak{F}_{(q,r)}(\gamma) \backslash 
\{\alpha_{\max,\gamma}\}$.  Note that since $\gamma$ is 
exterior, so is each offspring of $\gamma$.

It is useful to note the following general fact about norms on 
$\mathbb{R}^{2}$, which can be verified by a simple geometric
argument.  Let $C$ be a convex set; then
\begin{equation} \label{E:bdrydiam}
  \mathcal{W}(\partial C) \leq 6 \diam_{\tau}(C).
\end{equation}

Define
\[
   u(c,A) = \max(w_{1}A^{1/2} - cA^{1/6},0)
\]
and
\[
  D_{(q,r)}(\gamma) = 
  \sum_{\alpha \in \mathfrak{F}_{(q,r)}(\gamma)}
  \diam_{\tau}(\alpha), \qquad
  D_{(q,r)}^{\prime}(\gamma) = 
  \sum_{\alpha \in \mathfrak{F}_{(q,r)}^{\prime}(\gamma)}
  \diam_{\tau}(\alpha).
\]
The following lemma is related to (\ref{E:bdrydiam}).

\begin{lemma} \label{L:diamsum}
  Let $\gamma$ be a circuit, let $A = |\Int(\gamma)|$, 
  and let $q \geq 1, r \geq 15q$.  Then
  \begin{equation} \label{E:diamsum}
    D_{(q,r)}(\gamma) \geq \frac{1}{6}w_{1}\sqrt{A}.
  \end{equation}
\end{lemma}
\begin{proof}
We may assume $\gamma$ contains a clean $(q,r)$-bottleneck
$(u,v)$,
for otherwise (\ref{E:diamsum}) is immediate from (\ref{E:bdrydiam}).
We have 
\[
  q + 2\sqrt{2} \leq \frac{1}{4}r.
\]
Let $S$ denote the union of $Q_{2}(u,v)$ and all 
($\kt r/3$)-small offspring of $\gamma$, and let
\[
  R = \{z \in \mathbb{R}^{2}:  d_{\tau}(z,Q_{2}(u,v)) \leq 
  \kt r/3\}.
\]
Then 
\[
   S \subset R, \qquad \diam(R) \leq q + 2\sqrt{2} + \frac{2\sqrt{2}r}{3}
\]
and
\begin{equation} \label{E:Rbound}
  |R| \leq \left|\left\{z \in \mathbb{R}^{2}:  d(z,Q_{2}(u,v)) \leq 
  \frac{\sqrt{2}r}{3}\right\}\right| \leq \pi \left(\frac{q + 2\sqrt{2} + 
  \tfrac{2\sqrt{2}r}{3}}{2}\right)^{2}
  \leq r^{2}.
\end{equation}
Note that the set $\{\alpha_{1},\alpha_{2},...\}$ of
($\kt r/3$)-large offspring of $\gamma$ can be divided into two 
disjoint classes:  \emph{right offspring}, which intersect
$\gamma^{[u,v]}$, and \emph{left offspring}, which intersect
$\gamma^{[v,u]}$.  Also, every point of $\gamma$ is either in a left 
offspring, in a right offspring, or in $S$.
The diameter of $Q_{2}(u,v)$ is at most 
$q + 2\sqrt{2} \leq r/6$, while the diameters of $\gamma^{[u,v]}$
and $\gamma^{[v,u]}$ are at least $r$,
so the right and left classes each include at least one 
($5\kt r/6\sqrt{2}$)-large
offspring,.  Further, if
\begin{equation} \label{E:Ddiam}
  D_{(q,r)}(\gamma) \geq \diam_{\tau}(\gamma)
\end{equation}
then (\ref{E:diamsum}) follows from (\ref{E:bdrydiam}).  
Let $w$ and $x$ be sites of $\gamma$ with $d_{\tau}(w,x) = 
\diam_{\tau}(\gamma)$.
At least one of these points is not in $Q_{2}(u,v)$, so we 
may assume $w$ is in some $\alpha_{i}$.  There are now
four possibilities.  First, if also $x \in 
\alpha_{i}$, then (\ref{E:Ddiam}) holds.  Second, if 
instead $x \in S$, then there exists a
($5\kt r/6\sqrt{2}$)-large offspring 
$\alpha_{j} \neq \alpha_{i}$, and we have
\begin{align} 
  D_{(q,r)}(\gamma) &\geq \diam_{\tau}(\alpha_{i})
    +  \diam_{\tau}(\alpha_{j}) \notag \\
  &\geq \diam_{\tau}(\gamma) - \sqrt{2}\kt(q + 2\sqrt{2})
    - \frac{\kt r}{3} + \frac{5\kt r}{6\sqrt{2}} \notag \\
  &\geq \diam_{\tau}(\gamma) \notag
\end{align}
and again (\ref{E:Ddiam}) holds.  Third, suppose
$x \in \alpha_{k}$ for some $k \neq i$ and there exists a third 
($\kt r/3$)-large offspring $\alpha_{l}$ with $l \neq i, k$.
Let $d_{m} = 
\diam_{\tau}(\alpha_{m})$.  Then 
\[
  d_{i} + d_{k} \geq \diam_{\tau}(\gamma) - \sqrt{2}\kt(q + 2\sqrt{2})
  \geq \diam_{\tau}(\gamma) - d_{l}
\]
so once more, (\ref{E:Ddiam}) holds.

The fourth possibility is that $x \in \alpha_{k}$ for some
$k \neq i$ and $\alpha_{i}, \alpha_{k}$ are the only 
($\kt r/3$)-large offspring; each is necessarily actually
($5\kt r/6\sqrt{2}$)-large. 
From (\ref{E:bdrydiam}) we have
\[
  A \leq |R| + |\Int(\alpha_{i})| + |\Int(\alpha_{k})|
  \leq r^{2} +  \frac{36}{w_{1}^{2}}(d_{i}^{2} + d_{k}^{2}).
\]
Using this and the fact that $w_{1} \leq 4\kt$ (since the unit square
encloses unit area) we obtain
\[
  \frac{w_{1}^{2}}{36}A \leq \frac{4}{9}\kt^{2}r^{2}
  + d_{i}^{2} + d_{k}^{2} \leq 2d_{i}d_{k} + 
  d_{i}^{2} + d_{k}^{2}.
\]
Taking square roots yields (\ref{E:diamsum}).
\end{proof}

For $k \geq 0$ define the events
\begin{align}
  M_{y}(k,q,r,A,A^{\prime},d^{\prime},t) = 
    [|\mathfrak{F}_{(q,r)}^{\prime}&(\Gamma_{0})| = k] \cap
    [|\Int(\Gamma_{0})| = A] \cap [|\Int(\alpha_{\max,\Gamma_{0}})| 
    = A^{\prime}] \notag \\
  &\cap [D_{(q,r)}^{\prime}(\Gamma_{0}) \in 
    [d^{\prime},d^{\prime}+1)] \cap [\mathcal{W}(\partial
    \Co(\alpha_{\max,\Gamma_{0}})) \geq t]. \notag
\end{align}
We first consider $k = 0$, which means $\alpha_{\max,\gamma} = 
\gamma$ and $D_{(q,r)}(\gamma) = 0$; larger 
values will be handled later by induction.

\begin{proposition} \label{P:nobnbound}
  Assume (\ref{E:assump}) and either (i) the ratio weak mixing property
  or (ii) both the weak mixing property and the near-Markov property for 
  open circuits.
  Then there exist constants $\epsilon_{i},K_{i}$ as
  follows.  Let $A \geq 1, 
  t_{+} \geq 0, t = w_{1}\sqrt{A} + t_{+} \geq 2$, and
  $\epsilon_{16}t > r > 
  15q > K_{34}\log t$.  Then
  \begin{equation} \label{E:nobnbound}
    P(M_{0}(0,q,r,A,A,0,t)) \leq e^{-u(K_{35}r^{2/3},A) - 
    \frac{1}{2}t_{+}}.
  \end{equation}
\end{proposition}
\begin{proof}
From the definition of $w_{1}$ we may assume $t_{+} \geq 0$.
It follows easily 
from (\ref{E:conupr}) that for some $K_{36}, K_{37}$,
\begin{equation} \label{E:biggamma}
  P(\diam_{\tau}(\Gamma_{0}) \geq t) \leq K_{36}t^{4}e^{-t}
    \leq e^{-u(K_{37}r^{2/3},A) -  t_{+}}
\end{equation}
so by (\ref{E:bdrydiam}) it suffices to consider
\[
  M_{0}^{\prime}(q,r,A,t) = M_{0}(0,q,r,A,A,0,t)
    \cap [t/6 \leq \diam_{\tau}(\Gamma_{0}) \leq t].
\]
Suppose $\omega \in M_{0}^{\prime}(q,r,A,t)$.
Fix $\alpha > 0$ to be
specified, 
let $s = \alpha t^{2/3}r^{1/3}$ and suppose
$\HSkel_{s}(\Gamma_{0}) = (y_{0},..,y_{m+1})$.  
By (\ref{E:maxvert}),
\begin{equation} \label{E:maxvert3}
  m \leq K_{38}\alpha^{-1}t^{1/3}r^{-1/3}.
\end{equation}
Let $B_{i} = B(y_{i},
4r) \cap (\ZZ)^{*}$.  Let $I = \{i \leq m: |y_{i+1} - y_{i}| > 8r\}$.
For each $i \in I$ there is a segment $\Gamma_{0}^{[w_{i},x_{i}]}
\subset \Gamma_{0}^{[y_{i},y_{i+1}]}$ entirely outside $B_{i}
\cup B_{i+1}$, with
$w_{i} \in \partial B_{i}$ and $x_{i} \in \partial B_{i+1}$.  We next
show that
\begin{equation} \label{E:separated}
  d(\Gamma_{0}^{[w_{i},x_{i}]},\Gamma_{0}^{[w_{j},x_{j}]}) >
  q/2 \quad \text{for all } i \neq j \text{ in } I.
\end{equation}
If not, there exist $u \in \Gamma_{0}^{[w_{i},x_{i}]}$,
$v \in \Gamma_{0}^{[w_{j},x_{j}]}$ and a dual path $\psi$ from $u$ to
$v$ in $\Co(\Gamma_{0})$ with $|\psi| \leq q$.  Let $a$ be the last site of 
$\psi$ in $\Gamma_{0}^{[y_{i},y_{i+1}]}$, and $b$ the first site 
of $\psi$ after $a$ which is in some segment 
$\Gamma_{0}^{[y_{k},y_{k+1}]}$ with $k \neq i$.  Since all sites
$y_{l}$ are extreme points, we must have
$\psi^{(a,b)} \subset \Int(\Gamma_{0})$. 
We claim that $(a,b)$ is a $(q,3r)$-bottleneck.  By
Lemma \ref{L:bottleneck} this is a contradiction, so our claim will 
establish (\ref{E:separated}). Suppose $i < k$; the proof
if $i > k$ is similar.  We have $\psi \subset B(u,q)$ and $u \notin
B_{i+1}$ so $\psi \cap B(y_{i+1},3r+1) = \phi$.  Therefore 
$\Gamma_{0}^{[a,b]}$ contains a segment in $\omB(B_{i+1})$ which 
includes $y_{i+1}$ and has diameter at least $3r$.  Similarly 
since $v \notin B_{i}, \Gamma_{0}^{[b,a]}$ contains a segment 
in $\mathcal{B}(B_{i})$ which 
includes $y_{i}$ and has diameter at least $3r$.  This proves the claim
and thus (\ref{E:separated}).

From (\ref{E:extralength}) and
(\ref{E:maxvert3}) we have
\begin{align} \label{E:shortertau}
  \sum_{i \in I} \tau(x_{i} - w_{i}) &\geq 
    \mathcal{W}(\HPath_{s}(\Gamma_{0})) - K_{39}mr \\
  &\geq \mathcal{W}(\partial \Co(\Gamma_{0})) - K_{40}
    \alpha^{2} t^{1/3}r^{2/3} -     
    K_{39}mr \notag \\
  &\geq t - K_{41}(\alpha^{2} + \alpha^{-1})t^{1/3}r^{2/3}. \notag
\end{align}
Equation (\ref{E:shortertau}) shows that it is optimal to take $\alpha$ of 
order 1 in our choice of $s$, so we now set $\alpha = 1$.

For $w_{i} \in \partial(B_{i} \cap \ZZ)$ and $x_{i} \in \partial 
(B_{i+1} \cap \ZZ)$ for 
each $i \leq m$,
let $A(w_{0},x_{0},..,w_{m},x_{m})$ be the event 
that for each $i \in I$ there
is an open dual path $\alpha_{i}$ from $w_{i}$ to $x_{i}$
in $B_{\tau}(w_{i},t)$, with $d(\alpha_{i},\alpha_{j}) > q/2$
for all $i \neq j$.
Then we have shown
\begin{align} \label{E:L0bound1}
  P&(M_{0}^{\prime}(q,r,A,t)) \cap
    [\HSkel_{s}(\Gamma_{0}) = (y_{0},..,y_{m+1})]) \\
  &\leq \sum_{w_{0}} \sum_{x_{0}} \cdots \sum_{w_{m}}
    \sum_{x_{m}} P(A(w_{0},x_{0},..,w_{m},x_{m})) \notag \\
  &\leq (K_{42}r^{2})^{m+1} \max_ {w_{0},x_{0},..,w_{m},x_{m}}
    P(A(w_{0},x_{0},..,w_{m},x_{m})). \notag 
\end{align}
Provided $K_{34}$ is sufficiently large, Lemma \ref{L:decouple},
(\ref{E:shortertau}) and induction give
\begin{align} \label{E:L0bound2}
  P(A(w_{0},x_{0},..,w_{m},x_{m})) &\leq 2^{m} \prod_{i \in I}
    P(w_{i} \lra x_{i}) \\
  &\leq 2^{m} e^{-t + 2K_{41}t^{1/3}r^{2/3}} \notag
\end{align}
which with (\ref{E:maxvert3}) and (\ref{E:L0bound1}) yields
\begin{align} \label{E:L0bound3}
  P&(M_{0}^{\prime}(q,r,A,t) \cap
    [\HSkel_{s}(\Gamma_{0}) = (y_{0},..,y_{m+1})]) \\
  &\leq e^{-t + K_{43}t^{1/3}r^{2/3}}. \notag
\end{align}
The number of possible $(y_{0},..,y_{m+1})$ in 
(\ref{E:L0bound3}) is at most $(K_{44}t^{2})^{m+1}$, which with 
(\ref{E:L0bound3}) yields 
\begin{equation} \label{E:mprimebd}
  P(M_{0}^{\prime}(q,r,A,t) \leq 
  e^{-t + K_{45}t^{1/3}r^{2/3}}
\end{equation}
provided $K_{34}$,
and hence $r$, is large enough.  It is easily verified that, 
provided $\epsilon_{16}$ is small enough,
\[
  K_{45}(t_{+} + w_{1}A^{1/2})^{1/3}r^{2/3} \leq
  2K_{45}(w_{1}A^{1/2})^{1/3}r^{2/3} + \frac{1}{2}t_{+},
\]
by considering two cases according to which of $t_{+}$ and
$w_{1}A^{1/2}$ is larger.  This and (\ref{E:mprimebd}) 
establish (\ref{E:nobnbound}) for $M_{0}^{\prime}$;
as we have noted, this and 
(\ref{E:conupr}) establish (\ref{E:nobnbound}) as given.
\end{proof}

\begin{remark} \label{R:increvent}
Let $I$ be any increasing event.  Since the event on the left side of
(\ref{E:L0bound2}) is decreasing, its probability is not increased by
conditioning on $I$.  It follows easily that Proposition 
\ref{P:nobnbound} remains true if the probability on the left side
of (\ref{E:nobnbound}) is conditioned on $I$, even though 
$M_{0}(0,q,r,A,A,0,t)$ is not itself a decreasing event.
\end{remark}

Under (\ref{E:assump}), open dual bonds do not percolate, so for
every bounded set $A$ there is a.s. an innermost open circuit 
surrounding $A$; we denote this circuit by $\Theta(A)$.

An \emph{enclosure event} is an event of form
\[
  \cap_{i \leq n} (\Open(\alpha_{i}) \cap [\alpha_{i} \leftrightarrow
  \infty])
\]
where $\alpha_{1},..,\alpha_{n}$ are circuits (of regular bonds.)
This includes the degenerate case of the full space
$\{0,1\}^{\mathcal{B}(\ZZ)}$.
Clearly any such 
event is increasing.

\begin{proposition} \label{P:withbnbound}
  Assume (\ref{E:assump}), the weak mixing property and the
  near-Markov property for open circuits.  There exist
  constants $K_{i}, \epsilon_{i}$ as follows.
  Let $A \geq A^{\prime} \geq 3, k \geq 0, t_{+} \geq 0,
  t = w_{1}\sqrt{A^{\prime}} + t_{+}, d^{\prime} \geq 0,$ and
  $\epsilon_{17}(w_{1}\sqrt{A} + d^{\prime} + t_{+}) \geq r \geq
  15q \geq K_{46}\log A$.  
  Then
  \begin{equation} \label{E:mkbound2}
    P(M_{0}(k,q,r,A,A^{\prime},d^{\prime},t)) \leq
    \exp\left(- u(K_{47}r^{2/3},A)
    - \frac{1}{60}t_{+} - \frac{1}{10}d^{\prime}\right)
  \end{equation}
  and
  \begin{equation} \label{E:mkbound3}
    P(M_{0}(k,q,r,A,A^{\prime},d^{\prime},t)) 
    \leq \exp\left(-\frac{1}{2}d^{\prime}\right)
    P(M_{0}(0,q,r,A^{\prime},A^{\prime},0,t)). 
  \end{equation}
  \end{proposition}
\begin{proof}
We will refer to the requirement $\epsilon_{17}(w_{1}\sqrt{A} 
+ d^{\prime} + t_{+}) \geq r$ as the
\emph{size condition}, and to all other assumptions of the
Proposition collectively as the \emph{basic assumptions}.

We first prove (\ref{E:mkbound2}).
We proceed by induction on $k$, using Proposition
\ref{P:nobnbound} for $k = 0$.
Fix $q, r$ and define 
\[
  L_{y}(k,A,A^{\prime},d,d^{\prime},t) 
  = M_{y}(k,q,r,A,A^{\prime},d^{\prime},t) \cap
  [d \leq D_{(q,r)}(\Gamma_{y}) < d + 1], 
\]
where 
\begin{equation} \label{E:drange}
  \frac{1}{6}w_{1}\sqrt{A} - 1 \leq d \leq K_{48}A, \quad
  d \geq \frac{1}{6}t_{+} - 1.
\end{equation}
If $K_{48}$ is large enough then, from Lemma \ref{L:diamsum}
and the lattice nature of $\Gamma_{y}$,
$L_{y}(k,A,A^{\prime},$ $d,d^{\prime},t)$ is empty if any of the
inequalities in (\ref{E:drange}) fails.
Note that for some $K_{49}$,
\[
  M_{y}(k,q,r,A,A^{\prime},d^{\prime},t) \subset
  [\diam(\Gamma_{y}) \leq K_{49}A].
\]
Our induction hypothesis is that
for some constants $c_{i}$, for all $j < k$, all $A, A^{\prime},
t,d^{\prime}$ satisfying the basic assumptions, and all
$d$ satisfying (\ref{E:drange}), for any
enclosure event $E$, we have
\begin{equation} \label{E:induchyp}
  P(L_{0}(j,A,A^{\prime},d,d^{\prime},t) \mid E)
  \leq \exp(- \frac{9}{10}d + \frac{\kt}{40}jr);
\end{equation}
if the size condition is
also satisfied, then in addition
\begin{align} \label{E:induchyp2}
  P(L_{0}&(j,A,A^{\prime},d,d^{\prime},t) \mid E) \\
  &\leq \exp(- u(K_{35}r^{2/3},A)
    - \frac{1}{60}t_{+} - \frac{1}{10}d^{\prime}),
    \notag
\end{align}
with $K_{35}$ from (\ref{E:nobnbound}).
We wish to verify these hypotheses for $j = k$.

For $j = 0$ it suffices to consider $d^{\prime} = 0$ and 
(\ref{E:induchyp2}) is Proposition \ref{P:nobnbound} 
(together with Remark \ref{R:increvent}), while (\ref{E:induchyp}) follows 
easily from the first inequality in
(\ref{E:biggamma}), if $K_{46}$ is large.
Hence we may assume $k \geq 1$ and fix $A,A^{\prime},d,d^{\prime}$.
Let $\Psi(u,v,\eta)$ denote the event that
$L_{0}(k,A,A^{\prime},d,d^{\prime},t)$ occurs with $(u,v)$ a primary
$(q,r)$-bottleneck
in $\Gamma_{0}$ of type $\eta$, and let $\Phi(u,v,\eta) =
\{Y_{uv\eta}(\omega): \omega \in \Psi(u,v,\eta)\}$.
Let $E$ be an enclosure event; it is easy to see that
bottleneck surgery cannot destroy $E$, that is,
\[
  \Psi(u,v,\eta) \cap E \subset \Phi(u,v,\eta) \cap E.
\]
(This is the reason for considering only enclosure events, not general
increasing events.)
Since $|\partial Q_{1}(u,v)|
+ |J_{\eta}(u,v)| \leq K_{48}q$, we then have from the bounded
energy property:
\begin{equation} \label{E:psiphi}
  P(\Psi(u,v,\eta) \mid E) \leq e^{K_{49}q}P(\Phi(u,v,\eta) \mid E).
\end{equation}

Fix $u,v,\eta$ and for $y_{1},..,y_{l} \in \mathbb{Z}^{2}$ and
$l,m,A_{i},k_{i},d_{i},
d_{i}^{\prime} \geq 0$ in $\mathbb{Z}$, let
\[
  Z = Z(l,m,y_{1},..,y_{l},A_{1},..,A_{l}, d_{1},..,d_{l},
  d_{1}^{\prime},..,d_{l}^{\prime},k_{1},..,k_{l})
\]
denote the event that there exist disjoint exterior open dual circuits
$\alpha_{1},..,\alpha_{l}$ such that:
\begin{enumerate}
  \item[(i)] $\alpha_{i}$ surrounds
    $y_{i}$, $|\Int(\alpha_{i})| = A_{i}, \diam_{\tau}(\alpha_{i}) \geq 
    \kt r/3,
    |\mathfrak{F}^{\prime}_{(q,r)}(\alpha_{i})| = k_{i},
    d_{i}$ $\leq D_{(q,r)}(\alpha_{i})$ $< d_{i} + 1$ and $ d_{i}^{\prime}
    \leq D_{(q,r)}^{\prime}(\alpha_{i}) < d_{i}^{\prime} + 1$ for all 
    $i \leq l$;
  \item[(ii)] letting $\alpha_{\main}$
    denote the open dual circuit enclosing maximal area among all
    offspring of all $\alpha_{i}$, we have
    $\alpha_{\main}$
    an offspring of $\alpha_{m}$ satisfying
    $|\Int(\alpha_{\main})| = A^{\prime}$ and
    $\mathcal{W}(\partial \Co(\alpha_{\main})) \geq w_{1}\sqrt{A^{\prime}}
    + t_{+}$ ;
  \item[(iii)] there is no open dual path connecting
    $\alpha_{i}$ to $\alpha_{n}$ for $i \neq n$.
\end{enumerate}
We supress the parameters in $Z$ when
confusion is unlikely.
Then considering only
$(\kt r/3)$-large offspring we see that
\begin{equation} \label{E:Zunion}
  \Phi(u,v,\eta) \subset \cup \
  Z(l,m,y_{1},..,y_{l},A_{1},..,A_{l}, d_{1},..,d_{l},
  d_{1}^{\prime},..,d_{l}^{\prime},k_{1},..,k_{l})
\end{equation}
where the union is over all parameters satisfying
\[
  2 \leq l \leq \min(4q,k+1), \quad m \leq l, \quad y_{i}
  \in (\ZZ)^{*} \text{ with }
  d(y_{i},\zeta_{uv}) \leq 2, \quad A_{i} \geq \frac{\kt r}{6},
\]
\begin{equation} \label{E:ksum}
  A_{m} \geq A^{\prime}, \quad
  A - K_{33}r^{2} \leq \sum_{i \leq l} A_{i} \leq A, \quad
  \sum_{i \leq l} k_{i} = k + 1 - l,
\end{equation}
\begin{equation} \label{E:dirange}
  \frac{1}{6}w_{1}\sqrt{A_{i}} - 1 \leq d_{i} \leq K_{48}A_{i}, \quad
  d_{i}^{\prime} \leq d_{i}, 
\end{equation}
\begin{equation} \label{E:disums}
  d^{\prime} - l \leq d_{m}^{\prime} + \sum_{i \neq m} d_{i} \leq
  d^{\prime}, \quad
  d - l \leq \sum_{i\leq l} d_{i} \leq d,
\end{equation}
\begin{equation} \label{E:parameters}
  d_{m} - d_{m}^{\prime}
  + 1 \geq \frac{1}{6}(w_{1}\sqrt{A^{\prime}} + t_{+}).
\end{equation}
Here (\ref{E:parameters}) and the first inequality in
(\ref{E:dirange}) follow from (ii) above and
(\ref{E:bdrydiam}),and $K_{33}$ is from
(\ref{E:smalloff}). Temporarily fix such a set of parameters and
let $\nu_{1},..,\nu_{l}$ be circuits with
\[
  \diam_{\tau}(\nu_{i}) \geq \frac{\kt r}{3} \quad \text{for all $i$,   and }
  \Int(\nu_{i}) \cap \Int(\nu_{j}) = \phi \text{ for } i \neq j.
\]
Define events
\[
  \tilde{L}_{i} = \cup_{3 \leq B \leq A_{i}}
  L_{y_{i}}(k_{i},A_{i},B,d_{i},
  d_{i}^{\prime},w_{1}\sqrt{B}), \quad i \neq m,
\]
\[
  \tilde{L}_{m} = L_{y_{m}}(k_{m},A_{i},A^{\prime},d_{m},
  d_{m}^{\prime},t),
\]
\[
  T_{i} = [\Theta(\Gamma_{y_{i}}) = \nu_{i}], \quad
  Y_{i} = \tilde{L}_{i} \cap T_{i}
  \quad \text{for } i \leq m, \quad T = \cap_{i \leq l} T_{i}.
\]
Then
\begin{equation} \label{E:Zsubset}
  Z \cap T \subset \cap_{i \leq l} Y_{i}.
\end{equation}
Note that as in (\ref{E:outermost}),
\begin{equation} \label{E:theta}
  Y_{i} = \Open(\nu_{i}) \cap [\nu_{i} \leftrightarrow \infty] \cap G_{i}
  \text{    for some } G_{i} \in \mathcal{G}_{\mathcal{B}
  (\Int(\nu_{i}))}, \quad \text{for every } i \leq m.
\end{equation}
Define events
\[
  R = \cup_{i \leq l-1} \Int(\nu_{i}), \qquad F = \cap_{i \leq l-1}
  \Ext(\nu_{i}),
\]
\[
  \tilde{G} = \cap_{i \leq l-1} G_{i}, \qquad H = \cap_{i \leq l-1}
  \Open(\nu_{i}), \qquad N = \cap_{i \leq l-1} [\nu_{i} \leftrightarrow
  \infty],
\]
and let $L_{l}$ denote the event that $\tilde{L}_{l}$ occurs
on $\mathcal{B}(F)$.  Then
\begin{equation} \label{E:events1}
  N \cap E \cap H = E_{R} \cap E_{F} \cap H \quad \text{for some }
    E_{R} \in \mathcal{G}_{\mathcal{B}(R)}, \
    E_{F} \in \mathcal{G}_{\mathcal{B}(F)},
\end{equation}
\begin{equation} \label{E:events2}
  \cap_{i \leq l-1} Y_{i} = H \cap N \cap \tilde{G}
\end{equation}
and
\begin{equation} \label{E:events3}
  \cap_{i \leq l} Y_{i} \subset L_{l}.
\end{equation}
The relation between area $A_{i}$ and diameter
$d_{i}$ tells us roughly whether the circuit $\alpha_{i}$ (or its
collection of descendants) is regular or irregular;
we thereby subdivide the circuits into ``large regular,'' ``small regular''
and ``irregular'' categories as follows:
\[
  I_{1} = \{i \leq l: d_{i} < 4w_{1}\sqrt{A_{i}}, A_{i} \geq
  c_{1}r^{2}\},
\]
\[
  I_{2} = \{i \leq l: d_{i} < 4w_{1}\sqrt{A_{i}}, A_{i} <
  c_{1}r^{2}\},
\]
\[
  I_{3} = \{i \leq l: d_{i} \geq 4w_{1}\sqrt{A_{i}} \}
\]
where $c_{1} = \max(1/\epsilon_{17}^{2},(3K_{35}/w_{1})^{3})$
is chosen so that
\begin{equation} \label{E:J1bound1}
  u(K_{35}r^{2/3},A_{i}) \geq \frac{2}{3}w_{1}\sqrt{A_{i}},
  \quad i \in I_{1}.
\end{equation}
Let
\[
  \mu_{i} = \max\left[u(K_{35}r^{2/3},A_{i})
    + \frac{1}{10}d_{i},
    \frac{9}{10}d_{i} - \frac{\kt}{40}k_{i}r\right], \quad i \in I_{1}
    \backslash \{m\},
\]
\[
  \mu_{i} = \frac{9}{10}d_{i} - \frac{\kt}{40}k_{i}r, \quad
    i \in I_{2} \cup I_{3},
\]
\[
  \mu_{m} = \max\left[u(K_{35}r^{2/3},A_{m})
    + \frac{1}{60}t_{+} + \frac{1}{10}d_{m}^{\prime},
    \frac{9}{10}d_{m} - \frac{\kt}{40}k_{m}r\right] \quad
  \text{if } m \in I_{1} 
\]
(cf. (\ref{E:induchyp}).)
Now $H \cap E_{F}$ is an enclosure event so by the induction
hypotheses (\ref{E:induchyp}) and (\ref{E:induchyp2}), summing
over $B \leq A_{l}$ gives
\begin{equation} \label{E:Llbound}
  P(L_{l} \mid  H \cap E_{F}) \leq A_{l}e^{-\mu_{l}}.
\end{equation}
(Note that the size condition can be used here for $i \in I_{1}$.)
Since $|\nu_{i}| \geq \epsilon_{18}r$ for all $i$,
from (\ref{E:Markov2}) and (\ref{E:Markov3}), provided $K_{46}$ is
sufficiently large we get
\begin{equation} \label{E:usemkv1}
  P(L_{l} \cap  E_{F} \mid H \cap \tilde{G} \cap E_{R})
  \leq (1 + e^{-\epsilon_{18}r/2})P(L_{l} \cap  E_{F} \mid H)
\end{equation}
and
\begin{equation} \label{E:usemkv2}
  P( E_{F} \mid H) \leq (1 + e^{-\epsilon_{18}r/2})P( E_{F}
  \mid H \cap \tilde{G} \cap E_{R}).
\end{equation}
Combining (\ref{E:theta}) -- (\ref{E:usemkv2}) we obtain
\begin{align} \label{E:splitoff}
  P(&(\cap_{i \leq l} Y_{i}) \cap E) \\
  &\leq P(L_{l} \cap (\cap_{i \leq l-1} Y_{i}) \cap E)
    P(T_{l} \mid L_{l} \cap (\cap_{i \leq l-1} Y_{i}) \cap E) \notag \\
  &= P(L_{l} \cap H \cap  \tilde{G} \cap E_{R} \cap E_{F})
    P(T_{l} \mid L_{l} \cap (\cap_{i \leq l-1} Y_{i}) \cap E) \notag \\
  &\leq (1 + e^{- \epsilon_{18}r/2})P(L_{l} \cap  E_{F} \mid H)
    P(H \cap \tilde{G} \cap E_{R}) \notag \\
  &\qquad \cdot P(T_{l} \mid L_{l} \cap (\cap_{i \leq l-1} Y_{i}) \cap E)
    \notag \\
  &= (1 + e^{- \epsilon_{18}r/2})P(L_{l} \mid  E_{F} \cap H)
    P( E_{F} \mid H)
    P(H \cap \tilde{G} \cap E_{R}) \notag \\
  &\qquad \cdot P(T_{l} \mid L_{l} \cap (\cap_{i \leq l-1} Y_{i}) \cap E)
    \notag \\
  &\leq (1 + e^{- \epsilon_{18}r/2})^{2} A_{l}e^{-\mu_{l}} P( E_{F} \mid
    H \cap \tilde{G} \cap E_{R}) P(H \cap \tilde{G} \cap E_{R}) \notag \\
  &\qquad \cdot P(T_{l} \mid L_{l} \cap (\cap_{i \leq l-1} Y_{i}) \cap E)
    \notag \\
  &= (1 + e^{- \epsilon_{18}r/2})^{2} A_{l}e^{-\mu_{l}}
    P((\cap_{i \leq l-1} Y_{i}) \cap E)
    P(T_{l} \mid L_{l} \cap (\cap_{i \leq l-1} Y_{i}) \cap E). \notag
\end{align}
Summing over $\nu_{l}$ (which appears via $T_{l}$), dividing
by $P(I)$ and iterating this (taking $H$ and $N$ to be the full space
and $R = \phi, F = \mathbb{R}^{2}$, at the
last iteration step) we obtain using  (\ref{E:Zsubset})
\begin{equation} \label{E:ZIbound}
  P(Z \mid E) \leq (1 + e^{-\epsilon_{18}r/4})^{2l} \prod_{i \leq l} 
  A_{i}e^{-\mu_{i}}
  \leq 2A^{l} \exp(- \sum_{i \leq l} \mu_{i}).
\end{equation}

We now want to sum (\ref{E:ZIbound}) over all parameters of $Z$
allowed in (\ref{E:Zunion}).  We first view $l$ as fixed and allow the
other parameters to vary.  Note that the number of parameter choices
is at most $(K_{50}A)^{5l+1}$, and the number of possible $(u,v,\eta)$
is at most $(K_{50}A)^{3}$, for some $K_{50}$.  Suppose we can
show that, under the basic assumptions,
\begin{equation} \label{E:mubound1}
  \sum_{i \leq l} \mu_{i} \geq \frac{9}{10}d - \frac{\kt}{40}kr
  + \frac{\kt}{100}(l-1)r
\end{equation}
and, if the size condition is also satisfied,
\begin{equation} \label{E:mubound2}
  \sum_{i \leq l} \mu_{i} \geq u(K_{35}r^{2/3},A)
   + \frac{1}{10}d^{\prime}
   + \frac{1}{60}t_{+} + \frac{\kt}{100}(l-1)r.
\end{equation}
Then from (\ref{E:Zunion}) and (\ref{E:ZIbound}),
\[
  P(\Phi(u,v,\eta) \mid I) \leq 2(K_{50}A)^{5l+2}
    \exp\left(-\frac{9}{10}d + \frac{\kt}{40}kr - 
    \frac{\kt}{40}(l-1)r\right), \notag
\]
and if the size condition is satisfied,
\begin{align}
  P(\Phi(u,v,\eta) \mid I) \leq (K_{50}A)^{5l+2}
    \exp \biggl(-\max \biggl[u(K_{35}r^{2/3},&A)
   + \frac{1}{10}d^{\prime}
   + \frac{1}{60}t_{+}, \notag \\
  &\frac{9}{10}d - \frac{\kt}{40}kr\biggr] - \frac{\kt}{100}(l-1)r
    \biggr). \notag
\end{align}
Thus, summing over $u,v,\eta$, then over $l$,
and using $r > K_{46} \log A$
and (\ref{E:psiphi}),
we obtain (\ref{E:induchyp}) and (\ref{E:induchyp2}) for $j = k$.

Now (\ref{E:mubound1}) is a direct consequence of (\ref{E:ksum})
and (\ref{E:disums}),
so we turn to (\ref{E:mubound2}) and assume that the size condition holds.
Let
\begin{equation}
  \delta_{i} =
  \begin{cases}
    \frac{1}{10}d_{i}, &i \neq m \\
    \frac{1}{10}d_{m}^{\prime} + \frac{1}{60}t_{+}, &i = m
  \end{cases}
  \notag
\end{equation}
and set $\beta_{1} = \beta_{3} = 1/2, \beta_{2} = 1/10$.  We
claim that
\begin{equation} \label{E:claim}
  \mu_{i} \geq \beta_{j}w_{1}\sqrt{A_{i}} + \delta_{i} +
  \frac{1}{10}d_{i} - 1 \quad \text{for every } i \in I_{j}, j = 1,2,3.
\end{equation}
Observe that
\begin{equation} \label{E:dbounds}
  d_{i} \geq (k_{i}+1)\frac{\kt r}{3} - 1 \geq (k_{i}+1)\frac{\kt r}{4},
  \quad d_{i}^{\prime} \geq k_{i}\frac{\kt r}{3} - 1 
  \geq k_{i} \frac{\kt r}{4}
\end{equation}
and hence
\begin{equation} \label{E:mudbound}
  \mu_{i} \geq \frac{4}{5}d_{i}, \quad i \leq l.
\end{equation}
This yields
\begin{equation} \label{E:J3bound}
  \mu_{i} \geq \frac{4}{5}d_{i} \geq w_{1}\sqrt{A_{i}} +
  \frac{11}{20}d_{i}, \quad i \in I_{3},
\end{equation}
which proves (\ref{E:claim}) for $i \in I_{3} \backslash \{m\}$.
For $i \in I_{1}$ we can use (\ref{E:J1bound1}), 
(\ref{E:dbounds}) and a convex
combination of the lower bounds (\ref{E:mudbound}) and
$u(K_{35}r^{2/3},A_{i})$ for $\mu_{i}$ to obtain
\begin{equation} \label{E:J1bound2}
  \mu_{i} \geq \frac{3}{4}u(K_{35}r^{2/3},A_{i}) +
  \frac{1}{5}d_{i} 
  \geq \frac{1}{2}w_{1}\sqrt{A_{i}} + \frac{1}{5}d_{i},
  \quad i \in I_{1},
\end{equation}
which proves (\ref{E:claim}) for $i \in I_{1} \backslash \{m\}$.
From (\ref{E:dirange}),
\[
  d_{i} \geq \frac{1}{8}w_{1}\sqrt{A_{i}} + \frac{1}{4}d_{i} -
  \frac{3}{4}, \quad i \leq l,
\]
and hence by (\ref{E:mudbound})
\begin{equation} \label{E:J2bound2}
  \mu_{i} \geq \frac{4}{5} d_{i} \geq \frac{1}{10}w_{1}\sqrt{A_{i}}
    + \frac{1}{5}d_{i} - \frac{3}{5}, \quad i \leq l,
\end{equation}
which proves (\ref{E:claim}) for $i \in I_{2} \backslash \{m\}$.

We need slightly different estimates for $i = m$.
If $m \in I_{3}$ then using (\ref{E:J3bound}) and 
(\ref{E:parameters}) we obtain
\begin{align} \label{E:minJ3}
  \mu_{m} &\geq w_{1}\sqrt{A_{m}} +
    \frac{1}{2}d_{m} \\
  &\geq w_{1}\sqrt{A_{m}} + \frac{1}{4}d_{m} +
    \frac{1}{4}(d_{m}^{\prime} - 1) + \frac{1}{24}t_{+}, \notag
\end{align}
which proves (\ref{E:claim}) for $i = m$.  If $m \in I_{1}$ then
by (\ref{E:J1bound2}) and (\ref{E:parameters}),
\begin{align} \label{E:minJ1}
  \mu_{m} &\geq \frac{1}{2}w_{1}\sqrt{A_{m}} +
    \frac{1}{5}d_{m} \\
  &\geq \frac{1}{2}w_{1}\sqrt{A_{m}} + \frac{1}{10}d_{m} +
    \frac{1}{10}(d_{m}^{\prime} - 1) + \frac{1}{60}t_{+}, \notag
\end{align}
which again proves (\ref{E:claim}) for $i = m$.  Finally if
$m \in I_{2}$ then by (\ref{E:J2bound2}), (\ref{E:minJ1}) 
remains valid with $1/10$ in
place of $1/2$, and 3/5 subtracted from the right side,
once again proving (\ref{E:claim}) for
$i = m$. Thus (\ref{E:claim}) holds in all cases.

The next step is to sum (\ref{E:claim}) over $i$.  
There are 2 cases.

\emph{Case 1}.  $d \geq 20w_{1}\sqrt{A}$. Then using
(\ref{E:claim}), (\ref{E:disums}) and (\ref{E:dbounds}),
\begin{align} \label{E:firstcase}
  \sum_{i \leq l} \mu_{i} &\geq \sum_{i \leq l} (\delta_{i} +
    \frac{1}{10}d_{i} + 1) \\
  &\geq \frac{1}{10}(d^{\prime} - l) + \frac{1}{60}t_{+}
    + \frac{1}{20}(d - l)
    + \frac{1}{20}\sum_{i \leq l} d_{i} + l\notag \\
  &\geq \frac{1}{10}d^{\prime} 
    + \frac{1}{60}t_{+} + w_{1}\sqrt{A}
    + \frac{\kt}{80}lr \notag
\end{align}
which proves (\ref{E:mubound2}).

\emph{Case 2}.  $d < 20w_{1}\sqrt{A}$. 
By (\ref{E:disums}), (\ref{E:parameters}) and (\ref{E:dbounds}),
\[
  6d^{\prime} + t_{+} \leq 6(d + l + 1) \leq 7d \leq 140w_{1}\sqrt{A}.
\]
This and the size condition imply
\[
  r < 141\epsilon_{17}w_{1}\sqrt{A} \leq 
  \frac{\kt \sqrt{A}}{80K_{33}w_{1}}
\]
if $\epsilon_{17}$ is small enough, with $K_{33}$ as in 
(\ref{E:smalloff}), and then
\begin{equation} \label{E:areabound}
  \sqrt{A - K_{33}r^{2}} \geq \sqrt{A}\left(1 - 
  \frac{K_{33}r^{2}}{A}\right) \geq \sqrt{A} - \frac{\kt r}{80w_{1}}.
\end{equation}
Let
\[
  S_{j} = \sum_{i \in I_{j}} A_{i}, \quad j = 1,2,3, \quad \text{and}
  \quad S = S_{1} + S_{2} + S_{3}.
\]
Provided $\epsilon_{17}$ is small enough, we have by (\ref{E:smalloff}) 
and (\ref{E:areabound})
\begin{equation} \label{E:SAbound}
  w_{1}\sqrt{S} \geq w_{1}\sqrt{A} - \frac{\kt r}{80}.
\end{equation}
It is easily checked that
\begin{equation} \label{E:sqroot}
  \sqrt{a+b} \leq
  \sqrt{a} + \theta \sqrt{b} \quad \text{for } 0 \leq b \leq 4\theta^{2}a.
\end{equation}
Choose $i_{13},i_{2}$ satisfying
\[
  A_{i_{13}} = \max_{i \in I_{1} \cup I_{3}} A_{i}, \qquad
  A_{i_{2}} = \max_{i \in I_{2}} A_{i}.
\]
(If $I_{1} \cup I_{3}$ or $I_{2}$ is empty then the corresponding
$i_{13}$ or $i_{2}$ is undefined.)
We now consider two subcases.

\emph{Case 2a}.  $I_{2} = \phi$ or
$A_{i_{2}} \leq \tfrac{1}{25}(S_{1} + S_{3})$.
From the definition of $\mu_{i}$ if $i_{13} \in I_{1}$, and
from (\ref{E:claim}) if $i_{13} \in I_{3}$, we have
\begin{equation} \label{mui13}
  \mu_{i_{13}} \geq w_{1}\sqrt{A_{i_{13}}} - \lambda +
  \delta_{i_{13}},
\end{equation}
where $\lambda = \min(K_{35}r^{2/3}A_{i_{13}}^{1/6},
w_{1}\sqrt{A_{i_{13}}})$.  Therefore by (\ref{E:sqroot}) with
$\theta = 1/2$,
\begin{equation} \label{E:I13sum}
  \sum_{i \in I_{1} \cup I_{3}} \mu_{i} \geq w_{1}\sqrt{S_{1} +
  S_{3}} - \lambda + \sum_{i \in I_{1} \cup I_{3}} \delta_{i}
  + \sum_{i \in I_{1} \cup I_{3}, i \neq i_{13}} \frac{1}{10}d_{i}
  - l.
\end{equation}
Hence by (\ref{E:sqroot}) with $\theta = 1/25$, and
(\ref{E:SAbound}), (\ref{E:dbounds}),
(\ref{E:I13sum}) and (\ref{E:claim}),
\begin{align} \label{E:case1}
  \sum_{i \leq l} \mu_{i} &\geq w_{1}\sqrt{S} - \lambda +
    \sum_{i \leq l} \delta_{i} + \frac{\kt}{40}(l - 1)r - l \\
 &\geq w_{1}\sqrt{A} - \frac{\kt r}{80} - \lambda + 
    \frac{1}{10}(d^{\prime} - l)
    + \frac{1}{60}t_{+} + \frac{\kt}{40}(l - 1)r - l
    \notag \\
  &\geq u(K_{35}r^{2/3},A)  + \frac{1}{10}d^{\prime} +
    \frac{1}{60}t_{+} + \frac{\kt}{100}(l-1)r \notag
\end{align}
which gives (\ref{E:mubound2}).

\emph{Case 2b}.  $A_{i_{2}} > \tfrac{1}{25}(S_{1} + S_{3})$.
Let us relabel $(A_{i}, i \in J_{2})$ as $B_{1} \geq .. \geq B_{n}$.
We have
\[
  S_{1} + S_{3} \leq 25c_{1}r^{2},
\]
while from (\ref{E:smalloff}), provided $\epsilon_{17}$ is small enough,
\[
  S \geq A - K_{33}r^{2} \geq 32,525c_{1}r^{2},
\]
so
\[
  S_{2} \geq 32,500c_{1}r^{2} \geq 325\sum_{m = 1}^{100} B_{i},
\]
which implies
\[
  \frac{1}{20}\left(\sum_{m=101}^{n} B_{m}\right)^{1/2} \geq \frac{9}{10}
  \left(\sum_{m=1}^{100} B_{m}\right)^{1/2}.
\]
Using this, and using (\ref{E:sqroot}) twice (with $\theta = 1$ and
with $\theta = 1/20$), we get
\begin{align} \label{E:Bicompar}
  \frac{1}{10}&\sum_{m=1}^{n} \sqrt{B_{i}} \\
  &\geq \frac{1}{10}\left(\sum_{m=1}^{100} B_{m}\right)^{1/2}
    + \frac{1}{10}\sum_{m=101}^{n} \sqrt{B_{i}} \notag \\
  &\geq \left(\sum_{m=1}^{100} B_{m}\right)^{1/2} +
    \frac{1}{20}\sum_{m=101}^{n}
    \sqrt{B_{i}} + \frac{1}{20}
    \left(\sum_{m=101}^{n} B_{i}\right)^{1/2}
    - \frac{9}{10} \left(\sum_{m=1}^{100} B_{m}\right)^{1/2}
    \notag \\
  &\geq \sqrt{S_{2}}. \notag
\end{align}
If $I_{1} \cup I_{3} \neq \phi$ then (\ref{E:I13sum}) remains valid.
This, with (\ref{E:claim}), (\ref{E:dbounds}) and (\ref{E:Bicompar}),
shows that, whether $I_{1} \cup I_{3} = \phi$ or not, (\ref{E:case1})
(with $\lambda = 0$ if $I_{1} \cup I_{3} = \phi$) still holds.

The proof of (\ref{E:mubound2}), and thus of (\ref{E:induchyp2}),
is now complete.  Taking $y = 0$
and $E$ the full configuration space
in (\ref{E:induchyp2}), and summing over $d$ satisfying 
(\ref{E:drange}) shows that
\begin{equation} \label{E:smalldiam}
  P(M_{0}(k,q,r,A,A^{\prime},d^{\prime},t)) 
  \leq K_{48}A\exp\left(- u(K_{35}r^{2/3},A) - \frac{1}{60}t_{+}
    - \frac{1}{10}d^{\prime}\right). 
\end{equation}
This proves (\ref{E:mkbound2}),
with $K_{47} = 2K_{35}$.

It remains to prove (\ref{E:mkbound3}).  This is similar to the
proof of (\ref{E:induchyp2}), so we will only describe the 
changes.  Again fix $q,r$.  We make the same
induction hypothesis, except that (\ref{E:induchyp2}) is replaced by
\begin{equation} \label{E:induchyp3}
  P(L_{0}(j,A,A^{\prime},d,d^{\prime},t) \mid E) 
  \leq \exp\left(-\frac{4}{5}d^{\prime}\right)
  P(M_{0}(0,q,r,A^{\prime},A^{\prime},
  0,t) \mid E),
\end{equation}
and the requirement that the size condition be satisfied is removed.
This hypothesis is true for $j = 0$, where 
only $d^{\prime} = 0$ is relevant;
hence we fix $k \geq 1$ and $A, A^{\prime}, d, d^{\prime}$.  Let
\[
  f(A^{\prime},t) = - \log P(M_{0}(0,q,r,A^{\prime},A^{\prime},
  0,t) \mid E).
\]
In place of $\mu_{i}$ we use
\[
  \hat{\mu}_{i} = \frac{9}{10}d_{i} - \frac{\kt}{40}k_{i}r, \quad 
  i \neq m,
\]
\[
  \hat{\mu}_{m} = \max\left[\frac{4}{5}d_{m}^{\prime} + f(A^{\prime},t),
  \frac{9}{10}d_{m} - \frac{\kt}{40}k_{i}r\right].
\]
In place of (\ref{E:mubound1}), (\ref{E:mubound2}) and their
multi-case proofs, we have simply, using the first half of
(\ref{E:dbounds}),
\[
  \sum_{i \neq m} \left(\frac{1}{10}d_{i} - \frac{\kt}{40}k_{i}r\right) 
    \geq \frac{\kt r}{40}(l - 1) \geq l + \frac{\kt r}{50}(l - 1)
\]
and hence using (\ref{E:disums})
\begin{align}
  \sum_{i \leq l} \hat{\mu}_{i} &\geq f(A^{\prime},t) + 
    \frac{4}{5}(d^{\prime} - l) + \sum_{i \neq m} \left(\frac{1}{10}
    d_{i} - \frac{\kt}{40}k_{i}r\right) \notag \\
  &\geq f(A^{\prime},t) + \frac{4}{5}d^{\prime} 
    + \frac{\kt r}{50}(l - 1). \notag
\end{align}
This leads directly to (\ref{E:induchyp3}), as in the proof of
(\ref{E:induchyp2}).  In place of (\ref{E:smalldiam}) we have
\[
  P(M_{0}(k,q,r,A,A^{\prime},d^{\prime},t)) \leq 
    K_{48}A\exp\left(-\frac{4}{5}d^{\prime}\right)
    P(M_{0}(0,q,r,A^{\prime},A^{\prime},0,t)).
\]
For $k \geq 1$ we have $d^{\prime} \geq \kt r/3$, so 
provided $K_{46}$ is large
enough, this proves (\ref{E:mkbound3}).
\end{proof}

\begin{lemma} \label{L:hullsize}
  Let $q \geq 1, r \geq 15q$, let
  $\gamma$ be a dual circuit and let $\alpha_{\max,\gamma}$ be
  its maximal $(q,r)$-descendant.  Then
  \[
    \mathcal{W}(\partial \Co(\gamma))
    \leq \mathcal{W}(\partial \Co(\alpha_{\max,\gamma}))
    +  19D^{\prime}_{(q,r)}(\gamma).
  \]
\end{lemma}
\begin{proof}
We may assume $\gamma$ has at least one bottleneck.
If $(u,v)$ is a
primary bottleneck in $\gamma$, and the 
$(\kt r/3)$-large offspring of
$\gamma$ are $\alpha_{1},..,\alpha_{k}$, then
\[
  \Int(\gamma) \subset B_{\tau}(u,\frac{\kt r}{3} + 2\kappa_{\tau}q)
  \cup \cup_{i \leq k} \Int(\alpha_{i})
\]
and therefore from (\ref{E:bdrydiam}),
$\gamma$ can be surrounded by a (non-lattice)
loop of $\tau$-length at most
\[
  12\left(\frac{\kt r}{3} + 2\kappa_{\tau}q\right) + \sum_{i \leq k}
  \mathcal{W}(\partial \Co(\alpha_{i})).
\]
Since $\partial \Co(\gamma)$ minimizes the
$\tau$-length over all such loops, it follows that
\[
  \mathcal{W}(\partial \Co(\gamma)) \leq
  6\kt r +  \sum_{i \leq k} \mathcal{W}(\partial \Co(\alpha_{i})).
\]
Iterating this, and using (\ref{E:bdrydiam}) and
$D^{\prime}_{(q,r)}(\gamma) \geq
\tfrac{\kt r}{3}|\mathfrak{F}^{\prime}_{(q,r)}(\gamma)|$, we obtain
\begin{align}
  \mathcal{W}(\partial \Co(\gamma)) &\leq
    6\kt r|\mathfrak{F}^{\prime}_{(q,r)}(\gamma)|  +
    \sum_{\alpha \in \mathfrak{F}_{(q,r)}(\gamma)}
    \mathcal{W}(\partial \Co(\alpha)) \notag \\
  &\leq \mathcal{W}(\partial \Co(\alpha_{\max,\gamma}))
    + 19D^{\prime}_{(q,r)}(\gamma). \notag
\end{align}
\end{proof}

The next theorem, together with Theorem \ref{T:lowerbound},
shows roughly that for a droplet of size $A$, there is a cost for
the convex hull boundary $\tau$-length exceeding the minimum 
$w_{1}\sqrt{A}$ by an amount $s_{+}$, this cost being exponential
in $s_{+}$, and there is an exponential cost for positive 
$D_{(q,r)}^{\prime}(\Gamma_{0})$.

\begin{theorem} \label{T:onebound}
  Assume (\ref{E:assump}) and either (i) the ratio weak mixing
  property or (ii) both the weak mixing property and the
  near-Markov property for open circuits.  There exist constants
  $c_{i}$ as follows.  Let $A > K_{51}, s_{+} \geq 0,
  s = w_{1}\sqrt{A} + s_{+}$ and $d^{\prime} \geq 0$.  Then
  \begin{align} \label{E:onebound}
    P(|\Int(&\Gamma_{0})| \geq A, \mathcal{W}(\partial
      \Co(\Gamma_{0})) \geq s,
      D_{(q,r)}^{\prime}(\Gamma_{0}) \geq d^{\prime}) \\
    &\leq \exp\left(-u(K_{52}(\log A)^{2/3},A) - \frac{1}{1520}s_{+}
      - \frac{1}{20}d^{\prime}\right). \notag
  \end{align}
\end{theorem}
\begin{proof}
Let $K_{53} \geq K_{46}$ (of Proposition \ref{P:withbnbound}) and 
\[
  q_{B} = \frac{1}{15}K_{53}\log B, \qquad r_{B} = K_{53}\log B,
\]
\[
  t_{+}(n,A^{\prime}) = \max(s - 19n 
  - w_{1}\sqrt{A^{\prime}},0),
\]
\[
  I_{1}(n) = \{A^{\prime} \in \mathbb{Z}^{+}: t_{+}(n,A^{\prime})
  \geq \frac{s_{+}}{2}\},
\]
\[
  I_{2}(n) = \{A^{\prime} \in \mathbb{Z}^{+}: t_{+}(n,A^{\prime})
  < \frac{s_{+}}{2}, w_{1}\sqrt{A^{\prime}}
  \leq w_{1}\sqrt{A} + 19n\},
\]
and
\[
  I_{3}(n) = \{A^{\prime} \in \mathbb{Z}^{+}: t_{+}(n,A^{\prime})
  < \frac{s_{+}}{2}, w_{1}\sqrt{A^{\prime}}
  > w_{1}\sqrt{A} + 19n\}.
\]
Then using Lemma \ref{L:hullsize},
\begin{align} \label{E:sumbound}
  P(&|\Int(\Gamma_{0})| \geq A, \mathcal{W}(\partial
    \Co(\Gamma_{0})) \geq s,
    D_{(q,r)}^{\prime}(\Gamma_{0}) \geq d^{\prime}) \\
  &\leq \sum_{B \geq A} \sum_{n \geq d^{\prime}}
    P(|\Int(\Gamma_{0})| = B,
    D_{(q,r)}^{\prime}(\Gamma_{0}) \in [n,n+1), \notag \\
  &\qquad \qquad \qquad
    \mathcal{W}(\partial \Co(\alpha_{\max,\Gamma_{0}})) \geq
    s - 19n) \notag \\
  &\leq \sum_{B \geq A} \sum_{n \geq d^{\prime}}
    \Bigg[\sum_{A^{\prime} \leq B,A^{\prime} \in I_{1}(n)} \sum_{k \geq 0}
    P\bigl(M_{0}(k,q_{B},r_{B},B,A^{\prime},n,w_{1}\sqrt{A^{\prime}}
    + t_{+}(n,A^{\prime}))\bigr) \notag \\
  &\qquad \qquad + \sum_{A^{\prime} \leq B,A^{\prime} \in I_{2}(n)}
    \sum_{k \geq 0} P(M_{0}(k,q_{B},r_{B},B,A^{\prime},n,
    w_{1}\sqrt{A^{\prime}})) \notag \\
  &\qquad \qquad + \sum_{A^{\prime} \leq B,A^{\prime} \in I_{3}(n)}
    \sum_{k \geq 0} P(M_{0}(k,q_{B},r_{B},B,A^{\prime},n,
    w_{1}\sqrt{A^{\prime}})) \Bigg]. \notag 
\end{align}
The events $M_{0}(k,q_{B},r_{B},B,A^{\prime},n,\cdot)$ 
are empty unless $n + 1 > (k + 1)\kt r/4$ (cf. (\ref{E:dbounds})); 
if $K_{53}$, and hence $r$, is large enough, this implies $k \leq n$,
so we may restrict the sums in (\ref{E:sumbound}) to such $k$.
Presuming $A$ is large enough, $u(K_{47}(\log B)^{2/3},B)$
is strictly positive for all $B \geq A$.
For $A^{\prime} \in I_{1}(n)$  we apply Proposition \ref{P:withbnbound}
to get
\begin{align} \label{E:AprimeI1}
  &\sum_{A^{\prime} \leq B,A^{\prime} \in I_{1}(n)}
    \sum_{k \leq n} P(M_{0}(k,q_{B},r_{B},B,A^{\prime},n,
    w_{1}\sqrt{A^{\prime}} + t_{+}(n,A^{\prime})) \\
  &\qquad \leq 
    (n+1)B\exp\left(- w_{1}\sqrt{B} + K_{47}B^{1/6}(\log B)^{2/3}
    - \frac{s_{+}}{120} - \frac{1}{10}n\right). \notag
\end{align}
Note that if $I_{2}(n)$ or $I_{3}(n)$ is nonempty we must have 
$s_{+} = s - w_{1}\sqrt{A} > 0$.
If $A^{\prime} \in I_{2}(n)$ we have 
\[  
  \frac{1}{2}(s - w_{1}\sqrt{A}) > s - 19n - w_{1}\sqrt{A^{\prime}}
  \geq s - w_{1}\sqrt{A} - 38n
\]
and hence $n \geq s_{+}/64$.  Therefore
\begin{align} \label{E:AprimeI2}
  &\sum_{A^{\prime} \leq B,A^{\prime} \in I_{2}(n)}
    \sum_{k \leq n} P(M_{0}(k,q_{B},r_{B},B,A^{\prime},n,
    w_{1}\sqrt{A^{\prime}})) \\
  &\qquad \leq (n+1)B\exp\left(- w_{1}\sqrt{B} + K_{47}B^{1/6}(\log B)^{2/3}
    - \frac{1}{10}n\right) \notag \\
  &\qquad \leq (n+1)B\exp\left(- w_{1}\sqrt{B} + K_{47}B^{1/6}(\log B)^{2/3}
    - \frac{1}{20}n - \frac{1}{1520}s_{+}\right). \notag
\end{align}
If $A^{\prime} \in I_{3}(n)$ we have 
\[
  s_{+} + w_{1}\sqrt{A} - w_{1}\sqrt{A^{\prime}} - 19n
  \leq t_{+}(n,A^{\prime}) \leq \frac{s_{+}}{2}
\]
so that
\[
  2(w_{1}\sqrt{A^{\prime}} - w_{1}\sqrt{A}) > 
  w_{1}\sqrt{A^{\prime}} - w_{1}\sqrt{A} + 19n 
  \geq \frac{s_{+}}{2}
\]
which implies
\[
  w_{1}\sqrt{B} \geq w_{1}\sqrt{A^{\prime}}
  > w_{1}\sqrt{A} + \frac{s_{+}}{4}.
\]
Therefore
\begin{align} \label{E:AprimeI3}
  &\sum_{A^{\prime} \leq B,A^{\prime} \in I_{3}(n)}
    \sum_{k \leq n} P(M_{0}(k,q_{B},r_{B},B,A^{\prime},n,
    w_{1}\sqrt{A^{\prime}})) \\
  &\qquad \leq (n+1)B\exp\left(- w_{1}\sqrt{B} + K_{47}B^{1/6}(\log B)^{2/3}
    - \frac{1}{10}n\right) \notag \\
  &\qquad \leq (n+1)B\exp\left(- \frac{1}{2}w_{1}\sqrt{B} - 
    \frac{1}{2}w_{1}\sqrt{A}
    - \frac{1}{8}s_{+}
    + K_{47}B^{1/6}(\log B)^{2/3}
    - \frac{1}{10}n\right). \notag
\end{align}
We can now use (\ref{E:AprimeI1}),(\ref{E:AprimeI2}) and
(\ref{E:AprimeI3}) to sum over $n$ and $B$ in
(\ref{E:sumbound}), obtaining
(\ref{E:onebound}).
\end{proof}

Part of our main result is an easy consequence of Theorem \ref{T:onebound}.

\begin{proof}[Proof of Theorem \ref{T:main},(\ref{E:ALRmax}) and
(\ref{E:Hausmax})]
From the definition of $w_{1}$ and Theorem \ref{T:onebound}, 
for any $c$, if $A$ is sufficiently large,
\begin{align} \label{E:ALRbound}
  P(&|\Int(\Gamma_{0})| \geq A, ALR(\Gamma_{0}) > 
    cl^{1/3}(\log l)^{2/3}) \\
  &\leq P(|\Int(\Gamma_{0})| \geq A, |\Co(\Gamma_{0})|
    \geq A + cw_{1}l^{4/3}(\log l)^{2/3})
    \notag \\
  &\leq P\left(|\Int(\Gamma_{0})| \geq A, 
    \mathcal{W}(\partial\Co(\Gamma_{0})) \geq w_{1}\sqrt{A} + 
    \frac{cw_{1}^{2}l^{1/3}(\log l)^{2/3}}{3}\right) \notag \\
  &\leq \exp(-u(K_{52}(\log A)^{2/3},A) - \epsilon_{19}cl^{1/3}(\log l)^{2/3}).
    \notag
\end{align}
If we take $c$ sufficiently large, this and Theorem \ref{T:lowerbound}
prove that (\ref{E:ALRmax}) holds with
conditional probability approaching 1 as
$A \to \infty$.

Next, from the quadratic nature of the Wulff variational minimum (see 
\cite{Al92}, \cite{DKS}), for any $a, b$, if $A$ is sufficiently large,
\begin{align} \label{E:dHbound}
  P(A &\leq |\Int(\Gamma_{0})| \leq A + aw_{1}l^{4/3}(\log l)^{2/3},
    \Delta_{A}(\partial\Co(\Gamma_{0})) > bl^{2/3}(\log l)^{1/3}) \\
  &\leq \sum_{B} P(|\Int(\Gamma_{0})| = B,
    \Delta_{B}(\partial\Co(\Gamma_{0})) > 
    \frac{b}{2}l^{2/3}(\log l)^{1/3}) \notag \\
  &\leq \sum_{B} P(|\Int(\Gamma_{0})| = B,
    \mathcal{W}(\partial\Co(\Gamma_{0})) \geq w_{1}\sqrt{A} + 
    \epsilon_{20}bl^{1/3}(\log l)^{2/3}) \notag \\
  &\leq \exp(-u(K_{52}(\log A)^{2/3},A) - \epsilon_{21}bl^{1/3}(\log l)^{2/3}),
    \notag
\end{align}
where the sums are over $A \leq B \leq A + aw_{1}l^{4/3}(\log l)^{2/3}$.
Now 
\[
  P(|\Int(\Gamma_{0})| > A + aw_{1}l^{4/3}(\log l)^{2/3})
\]
can be bounded as in (\ref{E:ALRbound}), so if
we take $a, b$ sufficiently large, (\ref{E:dHbound})
and Theorem \ref{T:lowerbound} prove that (\ref{E:Hausmax}) 
holds with conditional probability approaching 1 as $A \to \infty$.
\end{proof}

\section{Proof of the Single Droplet Theorem}
Let $\Phi_{N}$ denote the open dual circuit in $\Lambda_{N}$ enclosing
maximal area.  Let
\[
  T_{N} = \sum_{\gamma \in \mathfrak{C}_{N}} 
  \diam_{\tau}(\gamma), \quad
  T_{N}^{\prime} = \sum_{\gamma \in \mathfrak{C}_{N} 
  \backslash \{\Phi_{N}\}} \diam_{\tau}(\gamma).
\]
Let $G_{N}(k,A,A^{\prime},d,d^{\prime})$ denote the event that there
are exactly $k+1 \ (K \log N)$-large exterior open dual circuits in 
$\Lambda_{N}$, with
\[
  \sum_{\gamma \in \mathfrak{C}_{N}} |\Int(\gamma)| = A, \quad
  |\Int(\Phi_{N})| = A^{\prime}, \quad d \leq T_{N} < d+1,
  \quad d^{\prime} \leq T_{N} < d^{\prime} + 1.
\]
As we will see, by
mimicking the proof of (\ref{E:mkbound3}) it is easy to obtain
\begin{align} \label{E:singlebd}
  P_{N,w}&(G_{N}(k,A,A^{\prime},d^{\prime},t)) \\
  &\leq \exp\left(-\frac{4}{5}d^{\prime}\right)P_{N,w}
    (|\Int(\Gamma_{y})| = A^{\prime} \text{ for some } y 
    \in \Lambda_{N} \cap \ZZ),
    \notag
\end{align}
and then, summing as in Theorem \ref{T:onebound}, roughly
\begin{align}\label{E:secondcost}
  P_{N,w}&(\sum_{\gamma \in \mathfrak{C}_{N}} |\Int(\gamma)| 
    \geq A, T_{N}^{\prime} \geq d^{\prime}) \\
  &\leq \exp\left(-\frac{1}{2}d^{\prime}\right)P_{N,w}(
    |\Int(\Gamma_{y})| \geq A - v \text{ for some } 
    y \in \Lambda_{N} \cap \ZZ). \notag
\end{align}
with $v \leq \sqrt{A}$.  (Statement (\ref{E:secondcost}) is for 
motivation only---the actual statement we prove is (\ref{E:sumbound4}).)
Note that $|\mathfrak{C}_{N}| > 1$ implies $T_{N}^{\prime} \geq K\log N$;
hence to prove (\ref{E:single}) we 
would like a result somewhat
like (\ref{E:secondcost}) but with $A$ on the right side in place
of $A - v$.  To replace $A - v$ with $A$
we need to know something of how the probability on the
right side of (\ref{E:secondcost}) behaves as a function of 
$v$, which is obtainable from our next two results.
For $N > 0$ and $0 < A < B$,
we say that a lattice site $y$ is $(\Lambda_{N},A,B)$-\emph{compatible}
if there exists $z$ such that $y \in z + \sqrt{A}\mathcal{K}_{1}$ and
$z + \sqrt{B}\mathcal{K}_{1} \subset \Lambda_{N}$.

\begin{proposition} \label{P:changeA}
  Let $P$ be a percolation model on $\mathcal{B}(\ZZ)$ satisfying
  (\ref{E:assump}), the near-Markov property for open circuits,
  and the ratio weak mixing property.  There exist
  $K_{i}, \epsilon_{i}$ such that for $A \geq K_{54}$ and $\delta \geq 
  K_{55}\log A$ we have 
  \begin{equation} \label{E:changeA}
    P(|\Int(\Gamma_{y})| \geq A + \delta\sqrt{A})
      \geq e^{-\epsilon_{22}\delta}P(|\Int(\Gamma_{y})| \geq A).
  \end{equation}
\end{proposition}

From this proposition we will obtain its analog for $P_{N,w}$,
which is as follows.

\begin{proposition} \label{P:changeA2}
  Let $P$ be a percolation model on $\mathcal{B}(\ZZ)$ satisfying
  (\ref{E:assump}), the near-Markov property for open circuits,
  and the ratio weak mixing property.  There exist
  $K_{i}, \epsilon_{i}$ such that for $N \geq 1, K_{56} \leq
  A \leq c_{2}N^{2}, K_{57}\log A
  < \delta \leq \epsilon_{23}\sqrt{A}$ and $y$ $(\Lambda_{N},A/2,
  (1 + \epsilon_{24})A)$-compatible we have
  \begin{equation} \label{E:changeA2}
    P_{N,w}(|\Int(\Gamma_{y})| \geq A + \delta\sqrt{A})
      \geq e^{-\epsilon_{25}\delta}P_{N,w}(|\Int(\Gamma_{y})| 
      \geq A).
  \end{equation}
  Here $c_{2}$ is from Theorem \ref{T:single}.
\end{proposition}

These propositions will require some
preliminary results.

\begin{lemma} \label{L:AAprime}
  Suppose $\tau$ is positive.  There exists $\epsilon_{26}$ such that 
  if $q \geq 1, r \geq 15q, A > A^{\prime} > 0$ 
  and $\gamma$ is a dual circuit with 
  $|\Int(\gamma)| = A, |\Int(\alpha_{\max,\gamma})| = A^{\prime}$
  then
  \[
    D^{\prime}_{(q,r)}(\gamma) \geq \epsilon_{26}\sqrt{A - A^{\prime}}.
  \]
\end{lemma}
\begin{proof}
Let $\alpha_{1},..,\alpha_{k}$ be the 
non-maximal final $(q,r)$-descendants of
$\gamma$, and $A_{i} = |\Int(\alpha_{i})|$.  Then
\[
  D^{\prime}_{(q,r)}(\gamma) \geq \sum_{i \leq k} 
  \max\left(w_{1}\sqrt{A_{i}},\frac{\kt r}{3}\right) \geq
  \frac{1}{2}w_{1}\sqrt{\sum_{i \leq k} A_{i}} +
  \frac{1}{6}k \kt r
\]
while (cf. (\ref{E:Rbound}))
\[
  A - A^{\prime} \leq kr^{2} + \sum_{i \leq k} A_{i}.
\]
The lemma follows easily.
\end{proof}

\begin{remark} \label{R:wired}
The proof of Proposition \ref{P:withbnbound} shows that the 
bounds (\ref{E:mkbound2}) and (\ref{E:mkbound3}) are valid
conditionally on any enclosure event $E$, with the following
modification:  on the right side of (\ref{E:mkbound3})
one must replace $P(M_{0}(\cdots))$ with 
$\sup_{y \in \Lambda_{N} \cap \ZZ}P_{N,w}(M_{y}(\cdots) \mid E)$.  
Let $E_{N}$ denote
the enclosure event $\Open(\partial \Lambda_{N}) \cap 
[\partial \Lambda_{N} \lra \infty]$.  Assuming $P$ has the
near-Markov property for open circuits, probabilities under
$P_{N,w}$ and under $P(\cdot \mid E_{N})$ differ
by a factor of at most $1 + Ce^{-aN}$ for some $C,a$.  
Therefore Proposition \ref{P:withbnbound} is valid with 
(\ref{E:mkbound2}) and (\ref{E:mkbound3}) replaced by
\begin{equation} \label{E:mkbound2alt}
  P_{N,w}(M_{y}(k,q,r,A,A^{\prime},d^{\prime},t)) \leq
    \exp(- u(K_{47}r^{2/3},A)
    - \frac{1}{60}t_{+} - \frac{1}{10}d^{\prime})
\end{equation}
and
\begin{equation} \label{E:mkbound3alt}
  P_{N,w}(M_{y}(k,q,r,A,A^{\prime},d^{\prime},t)) 
    \leq \exp\left(-\frac{1}{2}d^{\prime}\right)
    P_{N,w}(\cup_{x \in B(y,K_{58}A)}
    M_{x}(0,q,r,A^{\prime},A^{\prime},0,t))
\end{equation}
for all $y \in \Lambda_{N} \cap \ZZ$, for some $K_{58}$.
It follows easily from (\ref{E:mkbound2alt})
that Theorem \ref{T:onebound}, and then (\ref{E:ALRmax})
and (\ref{E:Hausmax}) of
Theorem \ref{T:main}, extend similarly.  
Further, because the 
constraint $A \leq c_{2}N^{2}$ in Theorem
\ref{T:single} and Proposition \ref{P:changeA}
allows the appropriate size
of Wulff shape to fit inside $\Lambda_{N}$, and weak mixing 
adequately eliminates boundary effects, Theorem
\ref{T:lowerbound} is also valid for $P_{N,w}$ 
when $A \leq c_{2}N^{2}$.
\end{remark}

For $x, y \in (\ZZ)^{*}, r > 0$ and $G \subset \mathbb{R}^{2}$,
we say there is an \emph{r-near dual connection} from $x$ to $y$ 
in $G$ if for some $u, v \in (\ZZ)^{*}$ with $d(u,v) \leq r$,
there are open dual paths from $x$ to $u$ and from $y$ to $v$
in $G$.  Let $N(x,y,r,G)$ denote the event that such a 
connection exists.
The following result is from \cite{Al97pwr}.

\begin{lemma}\label{L:rnear}
  Let $P$ be a percolation model on $\mathcal{B}(\ZZ)$ satisfying
  (\ref{E:assump})
  and the ratio weak mixing property.  There exist $K_{i}$ such that if 
  $|x| > 1$ and $r \geq K_{59}\log |x|$ then
  \[
    P(N(0,x,r,\mathbb{R}^{2})) \leq e^{-\tau(x) + K_{60}r}.
  \]
\end{lemma}

The next lemma proves (\ref{E:MLRmax}) in Theorem
\ref{T:main}.

\begin{lemma} \label{L:inward}
  Let $P$ be a percolation model on $\mathcal{B}(\ZZ)$ satisfying
  (\ref{E:assump}), the near-Markov property for open circuits,
  and the ratio weak mixing property.  There exist
  $\epsilon_{i},K_{i}$ such that for $A > K_{61},
  l = \sqrt{A}$ and $\epsilon_{27}A \geq r \geq 15q \geq K_{62}\log A$,
  under the measure $P(\cdot \mid |\Int(\Gamma_{0})| \geq A)$, with 
  probability approaching 1 as $A \to \infty$ we have
  \begin{equation} \label{E:MLRbound}
    MLR(\Gamma_{0}) \leq K_{63}l^{2/3}(\log l)^{1/3}.
  \end{equation}
\end{lemma}

We may and henceforth do assume that $K_{63} \geq K_{4}$ of Theorem 
\ref{T:main}.

\begin{proof}[Proof of Lemma \ref{L:inward}]
The proof is partly a modification of that of 
Proposition \ref{P:nobnbound}, so we use the notation of that proof.
The basic idea is that a large inward deviation of 
$\Gamma_{0}^{[y_{k},y_{k+1}]}$ from
$\partial\Co(\Gamma_{0})$ for some $k$ reduces the factor
$P(w_{k} \lra x_{k})$ in (\ref{E:L0bound2}).

First observe that the proof of (\ref{E:separated}) actually
shows that 
\begin{equation} \label{E:separated2}
  d(\Gamma_{0}^{[y_{i},y_{i+1}]},\Gamma_{0}^{[w_{j},x_{j}]}) >
  q/2 \quad \text{for all } i \neq j \text{ with } j \in I.
\end{equation}
Let $Z$ be the site in $\Gamma_{0}$ most distant from
$\partial\Co(\Gamma_{0})$.
Let $K_{64} > \sqrt{4\pi^{-1} w_{1}K_{3}}$
be a constant to be specified, and
suppose that $\Gamma_{0}$ satisfies (\ref{E:ALRmax}) and 
(\ref{E:Hausmax}), but for some $k \leq m$ and some $z$ we have
$Z = z \in \Gamma_{0}^{[y_{k},y_{k+1}]}$ and
$d(z,\partial\Co(\{y_{0},\ldots,y_{m+1}\})) 
> 51K_{64}l^{2/3}(\log l)^{1/3}$.  By (\ref{E:outside3})
we have $z \in \Co(\{y_{0},\ldots,y_{m+1}\})$.
Let $R_{k}$ denote the region bounded by 
$\Gamma_{0}^{[y_{k},y_{k+1}]}$ and by the
segment of $\partial\Co(\Gamma_{0})$ from
$y_{k}$ to $y_{k+1}$.  Let $\ell_{k}$ denote the line through
$y_{k}$ and $y_{k+1}$.  There exists a site $z^{\prime} \in
\ell_{k}$ with the following property:  a line tangent to 
$\partial B_{\tau}(y_{k},\tau(z^{\prime}-y_{k}))$ 
at $z^{\prime}$ passes 
through $z$.  This $z^{\prime}$ is a projection of $z$ onto
$\ell_{k}$ and satisfies
\begin{equation} \label{E:taucomp}
  \tau(y_{k} - z) \geq \tau(y_{k} - z^{\prime}), \quad
  \tau(y_{k+1} - z) \geq \tau(y_{k+1} - z^{\prime}).
\end{equation}
Let
\[
  D = B(z,K_{64}l^{2/3}(\log l)^{1/3}), \qquad
  E = B(z,2K_{64}l^{2/3}(\log l)^{1/3}),
\] 
\[
  D^{\prime} = B(z^{\prime},K_{64}l^{2/3}(\log l)^{1/3}), \qquad
  E^{\prime} = B(z^{\prime},13K_{64}l^{2/3}(\log l)^{1/3}),
\]
\[
  F^{\prime} = B(z^{\prime},14K_{64}l^{2/3}(\log l)^{1/3}).
\]
Write $D_{\mathbb{Z}}, D_{\mathbb{Z}}^{\prime}$, etc. for the 
intersection of $D, D^{\prime}$, etc. with $(\ZZ)^{*}$.
We have several cases.

\emph{Case 1.}  $\Gamma_{0}^{[y_{k},y_{k+1}]}$
contains a $(K_{59} \log A)$-near dual connection from $y_{k}$
to $y_{k+1}$ outside $E$.  Here $K_{59}$ is from Lemma \ref{L:rnear}.
Then this event and the event
$[z \lra \partial_{in}D_{\mathbb{Z}}]$ occur at separation 
$K_{64}l^{2/3}(\log l)^{1/3}$, 
so by Lemmas \ref{L:ratiowm} and \ref{L:rnear},
\begin{align} \label{E:rnearE}
  P&([z \lra \partial_{in}D_{\mathbb{Z}}] 
    \cap N(y_{k},y_{k+1},K_{59}\log A,
    E^{c})) \\
  &\leq 2P(z \lra \partial_{in}E_{\mathbb{Z}})
    P(N(y_{k},y_{k+1},K_{59}\log A,E^{c}) \notag \\
  &\leq \exp\left(-\tau(y_{k+1}-y_{k}) 
    - \frac{1}{2}K_{64}l^{2/3}(\log l)^{1/3}\right). \notag
\end{align} 

\emph{Case 2.}  $\Gamma_{0}^{[y_{k},y_{k+1}]}$
contains no $(K_{59} \log A)$-near connection from $y_{k}$
to $y_{k+1}$ outside $E$, and $D^{\prime} \cap \ell_{k}
\subset \overline{y_{k}y_{k+1}}$.  In this case,
if $\Gamma_{0}^{[y_{k},y_{k+1}]}$ did not
intersect $D^{\prime}$, then at least half of $D^{\prime}$  
would be contained in $R_{k}$.  (This requires (\ref{E:outside3}),
which shows that we cannot have $D^{\prime} 
\subset \Int(\Gamma_{0})$.)
Then from (\ref{E:ALRmax}),
\[
  \frac{1}{2}\pi K_{64}^{2}l^{4/3}(\log l)^{2/3} 
  \leq |R_{k}| \leq |\partial\Co(\Gamma_{0})|ALR(\Gamma_{0})
  \leq 2w_{1}K_{3}l^{4/3}(\log l)^{2/3},
\]
contradicting the definition of $K_{64}$.  Thus 
$\Gamma_{0}^{[y_{k},y_{k+1}]}$ must intersect $D^{\prime}$.  

Let $h_{1}$ be the first site of $\Gamma_{0}^{[y_{k},v]}$ in
$\partial E_{\mathbb{Z}}$, and let $h_{2}$ be the last site of 
$\Gamma_{0}^{[v,y_{k+1}]}$ in
$\partial E_{\mathbb{Z}}$.  
Note that by the definition of
Case 2, we have
\[
  d(\Gamma_{0}^{[y_{k},h_{1}]},\Gamma_{0}^{[h_{2},y_{k+1}]})
  > K_{59}\log A.
\]
We have two subcases within Case 2.

\emph{Case 2a}.  $\Gamma_{0}^{[y_{k},h_{1}]}$ and
$\Gamma_{0}^{[h_{2},y_{k+1}]}$ contain
$(K_{59}\log A)$-near connections
from $y_{k}$ to $h_{1}$ and from $h_{2}$ to $y_{k+1}$, 
respectively, both outside $F^{\prime}$.  Then these 
near-connections
occur at separation $K_{59}\log A$ from each other and from
the event $[\partial_{in} D^{\prime}_{\mathbb{Z}} \lra 
\partial E^{\prime}_{\mathbb{Z}}]$,
so by Lemmas \ref{L:decouple} and \ref{L:rnear}, since
$\tau(h_{2} - h_{1}) \leq \kt \sqrt{2}(4K_{64}l^{2/3}
(\log l)^{1/3} + 2)$,
\begin{align} \label{E:rnearFpr}
  P(&N(y_{k},h_{1},K_{59}\log A,(F^{\prime})^{c}),
    N(h_{2},y_{k+1},K_{59}\log A,(F^{\prime})^{c}),
    \text{ and } [\partial_{in} D^{\prime}_{\mathbb{Z}}
    \lra \partial E^{\prime}_{\mathbb{Z}}] \\
  &\qquad \text{ occur at mutual separation } K_{59}\log A
    \text{ for some } h_{1}, h_{2} \in \partial E_{\mathbb{Z}}) \notag \\
  &\leq 4\sum_{h_{1}, h_{2} \in \partial E_{\mathbb{Z}}} 
    P(N(y_{k},h_{1},K_{59}\log A,(F^{\prime})^{c}) \notag \\
  &\qquad \qquad \cdot P(N(h_{2},y_{k+1},K_{59}\log A,(F^{\prime})^{c}))
    P(\partial_{in} D^{\prime}_{\mathbb{Z}} \lra 
    \partial E^{\prime}_{\mathbb{Z}}) \notag \\
  &\leq 4|\partial_{in} D^{\prime}_{\mathbb{Z}}|
    |\partial E^{\prime}_{\mathbb{Z}}|
    \sum_{h_{1}, h_{2} \in \partial E_{\mathbb{Z}}}
    \exp\Bigl(-\tau(h_{1} - y_{k}) - \tau(y_{k+1} - h_{2}) \notag \\
  &\qquad \qquad \qquad \qquad \qquad \qquad
    - 6\sqrt{2}K_{64}\kt l^{2/3}(\log l)^{1/3} 
    + K_{65}\log A\Bigr) \notag \\
  &\leq \exp(-\tau(y_{k+1} - y_{k}) - 
    K_{64}\kt l^{2/3}(\log l)^{1/3}). \notag
\end{align}
(Here we have actually used a trivial extension of Lemma 
\ref{L:decouple}, since the near-connection events occur not 
necessarily on the dual cluster of a single $x^{*}$ but rather on the 
union of two such clusters.)

\emph{Case 2b}.  $\Gamma_{0}^{[y_{k},h_{1}]}$ does not contain
a $(K_{59}\log A)$-near connection
from $y_{k}$ to $h_{1}$ outside $F^{\prime}$.  (The other 
alternative to Case 2a within Case 2, symmetric to this one, is that 
$\Gamma_{0}^{[h_{2},y_{k+1}]}$ does not contain
a $(K_{59}\log A)$-near connection
from $h_{2}$ to $y_{k+1}$ outside $F^{\prime}$; the proof
is similar.)  In this case the events $[y_{k} \lra 
\partial F_{\mathbb{Z}}^{\prime}$ outside $F^{\prime} \cup E],
[\partial F_{\mathbb{Z}}^{\prime} \lra \partial E_{\mathbb{Z}}$ 
outside $F^{\prime}
\cup E]$ and $[\partial E_{\mathbb{Z}} \lra y_{k+1}$ outside 
$E]$ occur at separation $K_{59}\log A$ from each other.  Further,
we have for $f, g \in \partial F_{\mathbb{Z}}^{\prime}$, 
using (\ref{E:taucomp}),
\begin{align} \label{E:taubd1}
  \tau(&f - y_{k}) + \tau(y_{k+1} - h_{2}) + \tau(h_{1} - g) \\
  &\geq \tau(z^{\prime} - y_{k}) + \tau(y_{k+1} - z) 
    - \sqrt{2}\kt(|z^{\prime} - f| + |z - h_{2}|) \notag \\
  &\qquad \qquad + \frac{\kt}{\sqrt{2}}(|z - z^{\prime}| 
    - |z - h_{1}| - |z^{\prime} - g|) \notag \\
  &\geq \tau(z^{\prime} - y_{k}) + \tau(y_{k+1} - z^{\prime})
    + 2K_{64}\kt l^{2/3}(\log l)^{1/3} \notag \\
  &\geq \tau(y_{k+1} - y_{k}) 
    + 2K_{64}\kt l^{2/3}(\log l)^{1/3}. \notag
\end{align}
Hence by Lemma \ref{L:decouple},
\begin{align}\label{E:nornearFpr}
  P(&[y_{k} \lra \partial F_{\mathbb{Z}}^{\prime} \text{ outside } 
    F^{\prime} \cup E], [\partial F_{\mathbb{Z}}^{\prime} \lra 
    \partial E_{\mathbb{Z}}
    \text{ outside } F^{\prime} \cup E] \text{ and } \\
  &\qquad [\partial E_{\mathbb{Z}} \lra y_{k+1} \text{ outside } 
    F^{\prime} \cup E] 
    \text{ occur at mutual separation } K_{59}\log A) \notag \\
  &\leq 4P(y_{k} \lra \partial F_{\mathbb{Z}}^{\prime})
    P(\partial F_{\mathbb{Z}}^{\prime} \lra \partial E_{\mathbb{Z}})
    P(\partial E_{\mathbb{Z}} \lra y_{k+1}) \notag \\
  &\leq 4\sum_{h_{1},h_{2} \in \partial F_{\mathbb{Z}}^{\prime}}
    \sum_{f, g \in \partial E_{\mathbb{Z}}} \exp(-
    \tau(f - y_{k}) - \tau(y_{k+1} - h_{2}) - \tau(h_{1} - g)) \notag \\
  &\leq \exp(-\tau(y_{k+1} - y_{k}) - 
    K_{64}\kt l^{2/3}(\log l)^{1/3}). \notag
\end{align}

\emph{Case 3}.  $\Gamma_{0}^{[y_{k},y_{k+1}]}$
contains no $(K_{59} \log A)$-near connection from $y_{k}$
to $y_{k+1}$ outside $E$, and $D^{\prime} \cap \ell_{k}
\not\subset \overline{y_{k}y_{k+1}}$.  As in Case 2b
there are two symmetric alternatives within this, and we
need only consider one, so we assume $d(z,y_{k+1}) \leq
d(z,y_{k})$.  Then the events $[y_{k} \lra \partial E_{\mathbb{Z}}]$ and
$[\partial E_{\mathbb{Z}} \lra y_{k+1}]$ occur at separation 
$K_{59} \log A$ and we have
\[
  \tau(z - y_{k}) \geq \tau(z^{\prime} - y_{k}) \geq
  \tau(y_{k+1} - y_{k})
  - K_{64}\sqrt{2}\kt l^{2/3}(\log l)^{1/3}.
\]
Hence
\begin{align} \label{E:farz}
  P(&[y_{k} \lra \partial E_{\mathbb{Z}}] \text{ and } 
    [\partial E_{\mathbb{Z}} \lra y_{k+1}] \text{ occur at separation }
    K_{59} \log A) \\
  &\leq 2P(y_{k} \lra \partial E)P(\partial E \lra y_{k+1})
    \notag \\
  &\leq 2\sum_{h_{1},h_{2} \in \partial E} \exp(-\tau(h_{1} - y_{k})
    - \tau(y_{k+1} - h_{2})) \notag \\
  &\leq 2\exp\left(-\tau(z - y_{k})  
    + \sqrt{2}\kt |z - h_{1}|) - \frac{\kt}{\sqrt{2}}
    (|y_{k+1} - z| - |z - h_{2}|\right) \notag \\
  &\leq 2\exp(-\tau(y_{k+1} - y_{k}) 
    - 30K_{64}\kt l^{2/3}(\log l)^{1/3}). \notag
\end{align}

Now let $J(y_{k},y_{k+1})$ denote the event that there is an
open dual path from $y_{k}$ to $y_{k+1}$ containing a site
$z$ with $d(z,\ell_{k}) \geq 51K_{64}l^{2/3}(\log l)^{1/3}$.
Combining the three cases we obtain 
\[
  P(J(y_{k},y_{k+1})) \leq 3\exp(-\tau(y_{k+1} - y_{k})
    - \frac{1}{2}K_{64}\kt l^{2/3}(\log l)^{1/3}),
\]
and then analogously to (\ref{E:L0bound2}),
\[
  P(A(w_{0},x_{0},\ldots,w_{m},x_{m}) \cap J(y_{k},y_{k+1}))
  \leq 2^{m}P(J(y_{k},y_{k+1}))\prod_{i \in I \backslash \{k\}}
  P(w_{i} \lra x_{i}).
\]
Analogously to (\ref{E:L0bound2}) -- (\ref{E:mprimebd}), this 
leads to
\begin{align} \label{E:bigMLRbd}
  P(&|\Int(\Gamma_{0})| \geq A, MLR(\Gamma_{0}) \geq 
    51K_{64}l^{2/3}(\log l)^{1/3}, \Gamma_{0}
    \text{ is } (q,r)-\text{bottleneck-free}) \\
  &\leq \exp(-w_{1}\sqrt{A}
    - \frac{1}{3}K_{64}\kt l^{2/3}(\log l)^{1/3}). \notag
\end{align}

We will need the following straightforward extension of (\ref{E:mkbound3}),
under the conditions of Proposition \ref{P:withbnbound}:
\begin{align} \label{E:55extn}
  P(M_{0}&(k,q,r,A,A^{\prime},d^{\prime},t) \cap
    [MLR(\alpha_{\max,\Gamma_{0}}) \geq 51K_{64}l^{2/3}(\log l)^{1/3}]) \\
  &\leq \exp\left(-\frac{1}{2}d^{\prime}\right)
    P\bigl(M_{0}(0,q,r,A^{\prime},A^{\prime},0,t) \cap
    [MLR(\Gamma_{0}) \geq 51K_{64}l^{2/3}(\log l)^{1/3}]\bigr). \notag
\end{align}

It is easy to see (cf. the proof of Lemma \ref{L:hullsize})
that
\begin{equation} \label{E:alphaMLR}
  MLR(\alpha_{\max,\Gamma_{0}}) \geq MLR(\Gamma_{0}) - 
  3\kt^{-1} D^{\prime}_{(q,r)}(\Gamma_{0}).
\end{equation}
Let $g(A) = (3\epsilon_{26})^{-2}K_{64}^{2}\kt^{2}l^{4/3}(\log l)^{2/3}$,
with $\epsilon_{26}$ from Lemma \ref{L:AAprime}.
From Theorem \ref{T:onebound}, 
Lemma \ref{L:AAprime}, (\ref{E:bigMLRbd}), (\ref{E:alphaMLR})
and (\ref{E:55extn}),
\begin{align} \label{E:M0bound}
  P&(|\Int(\Gamma_{0})| \geq A, MLR(\Gamma_{0}) \geq 
    52K_{64}l^{2/3}(\log l)^{1/3}) \\
  &\leq P(|\Int(\Gamma_{0})| \geq A, D_{(q,r)}^{\prime}(\Gamma_{0}) 
    \geq \frac{1}{3}K_{64}\kt 
    l^{2/3}(\log l)^{1/3}) \notag \\
  &\qquad + P(|\Int(\Gamma_{0})| \geq A, MLR(\alpha_{\max,\Gamma_{0}}) 
    \geq 51K_{64}l^{2/3}(\log l)^{1/3}, \notag \\
  &\qquad \qquad \qquad |\Int(\Gamma_{0})|
    - |\Int(\alpha_{\max,\Gamma_{0}})| < g(A)) \notag \\
  &\leq \exp(-w_{1}\sqrt{A} 
    - \frac{1}{60}K_{64}\kt l^{2/3}(\log l)^{1/3}) \notag \\
  &\qquad + \sum_{B \geq Q} \quad \sum_{d^{\prime} \geq 0} 
    \quad \sum_{B - g(A) < 
    A^{\prime} \leq B} \quad \sum_{k \leq d^{\prime}} \quad
    P\Bigl(M_{0}(k,q,r,B,A^{\prime},d^{\prime},w_{1}\sqrt{A^{\prime}})
    \notag \\
  &\qquad \qquad \qquad \qquad \qquad \qquad \qquad \qquad \cap
    [MLR(\alpha_{\max,\Gamma_{0}}) \geq 51K_{64}l^{2/3}(\log l)^{1/3}]
    \Bigr) \notag \displaybreak[0]\\
  &\leq \exp(-w_{1}\sqrt{A} 
    - \frac{1}{60}K_{64}\kt l^{2/3}(\log l)^{1/3}) \notag \\
  &\quad + \sum_{A^{\prime} \geq A - g(A)} \quad \sum_{A^{\prime} \leq
    B < A^{\prime} + g(A)} \quad \sum_{d^{\prime} \geq 0} \quad
    (d^{\prime} + 1)e^{-d^{\prime}/2}
    P\Bigl(M_{0}(0,q,r,A^{\prime},A^{\prime},0,
    w_{1}\sqrt{A^{\prime}}) \notag \\
  &\qquad \qquad \qquad \qquad \qquad \qquad \qquad \qquad \cap
    [MLR(\Gamma_{0}) \geq 51K_{64}l^{2/3}(\log l)^{1/3}]\Bigr) \notag \\
  &\leq \exp(-w_{1}\sqrt{A} 
    - \frac{1}{60}K_{64}\kt l^{2/3}(\log l)^{1/3}) \notag \\
  &\qquad + 10g(A) P\Bigl(|\Int(\Gamma_{0})| \geq 
    A - g(A), \Gamma_{0} \text{ is } 
    (q,r)-\text{bottleneck-free}, \notag \\
  &\qquad \qquad \qquad \qquad MLR(\Gamma_{0}) \geq 
    51K_{64}l^{2/3}(\log l)^{1/3}\Bigr) \notag \\
  &\leq \exp(-w_{1}\sqrt{A} 
    - \frac{1}{70}K_{64}\kt l^{2/3}(\log l)^{1/3}). \notag 
\end{align}
Here the restriction to $k \leq d^{\prime}$ is permissible as in the 
proof of Theorem \ref{T:onebound}.  Provided we take $K_{64}$
large enough,
with Theorem \ref{T:lowerbound}
this completes the proof.
\end{proof}

\begin{proof}[Proof of Proposition \ref{P:changeA}]
Let $l = \sqrt{A}$ and $r = 15q = K_{46}\log A$, where $K_{46}$ is from
Proposition \ref{P:withbnbound}.
Let $\Mid(\Gamma_{0})$ denote the Wulff shape of area 
$\tfrac{3}{4}|\Gamma_{0}|$ centered at the center of mass of 
$\Int(\Gamma_{0})$.
From translation invariance,
Theorems \ref{T:main} and \ref{T:lowerbound}, (\ref{E:DeltaA})
and (\ref{E:M0bound}) we have (recalling $K_{63} = 52K_{64}$)
\begin{align}
  &P(|\Int(\Gamma_{0})| \geq A, D_{(q,r)}^{\prime}(\Gamma_{0}) 
    \geq \frac{1}{3}K_{64}\kt 
    l^{2/3}(\log l)^{1/3}) \notag \\
  &+ P(|\Int(\Gamma_{0})| \geq A, \Delta_{A}(\Gamma_{0}) >
    2K_{63}l^{2/3}(\log l)^{1/3}) \notag \\
  &+ P(|\Int(\Gamma_{0})| \geq A, \Delta_{A}(\Gamma_{0}) \leq
    2K_{63}l^{2/3}(\log l)^{1/3}, 0 \notin 
    \Mid(\Gamma_{0})) \notag \\
  &\qquad \leq \frac{1}{2}P(|\Int(\Gamma_{0})| \geq A). \notag
\end{align}
With this fact, we can repeat the argument of 
(\ref{E:M0bound}), but excluding reference to 
$MLR(\cdot)$, to obtain (using again
$g(A) = (3\epsilon_{26})^{-2}K_{64}^{2}\kt^{2}l^{4/3}(\log l)^{2/3}$)
\begin{align} \label{E:M0bound3a}
  P(|&\Int(\Gamma_{0})| \geq A) \\
  &\leq \frac{1}{2}P(|\Int(\Gamma_{0})| \geq A)
    \notag \\
  &\qquad + P\bigl(|\Int(\alpha_{\max,\Gamma_{0}})| \geq 
    A - g(A), \notag \\
  &\qquad \qquad \qquad \Delta_{A}(\alpha_{\max,\Gamma_{0}})
    \leq 2K_{63}l^{2/3}(\log l)^{1/3}, 0 \in 
    \Mid(\Gamma_{0})\bigr) \notag \\
  &\leq \frac{1}{2}P(|\Int(\Gamma_{0})| \geq A)
    \notag \\
  &\qquad + 10g(A) P\bigl(|\Int(\Gamma_{0})| \geq 
    A - g(A), \Gamma_{0} \text{ is } 
    (q,r)-\text{bottleneck-free, } \notag \\
  &\qquad \qquad \qquad \qquad \qquad \Delta_{A}(\Gamma_{0})
    \leq 2K_{63}l^{2/3}(\log l)^{1/3}, 0 \in 
    \Mid(\Gamma_{0})\bigr), \notag
\end{align}
so that 
\begin{align} \label{E:M0bound3}
  P(|&\Int(\Gamma_{0})| \geq A) \\
  &\leq 20g(A) P\bigl(|\Int(\Gamma_{0})| \geq 
    A - g(A), \Gamma_{0} \text{ is } 
    (q,r)-\text{bottleneck-free, } \notag \\
  &\qquad \qquad \qquad \qquad \qquad \Delta_{A}(\Gamma_{0})
    \leq 2K_{63}l^{2/3}(\log l)^{1/3}, 0 \in 
    \Mid(\Gamma_{0})\bigr). \notag
\end{align}

The idea now is to split $\Gamma_{0}$ into two halves and 
approximate the probability on the right
side of (\ref{E:M0bound3}) by the 
product of the probabilities of the two halves.  With this
independence the two halves can in effect be pulled apart
from one another to increase the area enclosed by $\Gamma_{0}$
at only a small cost in increased boundary length.  To
accomplish this we first need some definitions.  Let $\rho$ be a 
path from $x_{2} = (a_{2},b_{2})$ to $x_{1} = (a_{1},b_{1})$
in the slab $S_{x_{1}x_{2}} = \{(x,y) \in 
\mathbb{R}^{2}: b_{2} \leq y \leq b_{1}\}$.  Let $J_{L}(\rho)$
and $J_{R}(\rho)$ denote the regions to the left and right,
respectively, of $\rho$ in $S_{x_{1}x_{2}}$.  The 
\emph{right-side area} determined by $\rho$ is 
\[
  \mu_{R}(\rho) = |J_{R}(\rho) \cap \Lambda_{N}|
  - |J_{R}(\overline{x_{1}x_{2}}) \cap \Lambda_{N}|,
\]
evaluated for $N$ large enough
that $\Lambda_{N}$ contains $\rho$.  (Note that for such $N$,
the right-side area does not vary with $N$.  Also, in our
definition the path $\rho$ must be oriented so that
$b_{1} \leq b_{2}$.)  The \emph{left-side area} $\mu_{L}(\rho)$
is defined similarly using the left side of $\rho$.  Let
$X_{1}$ and $X_{2}$ be the points of $\Gamma_{0}$ of 
maximum and minimum second coordinate, respectively, using
the leftmost if there is more than one.  Then
\[
  |\Int(\Gamma_{0})| = \mu_{L}(\Gamma_{0}^{[X_{2},X_{1}]})
    +  \mu_{R}(\hat{\Gamma_{0}}^{[X_{2},X_{1}]}),
\]
where $\hat{\Gamma_{0}}$ is $\Gamma_{0}$ traversed in the 
direction of negative orientation.  Let
\[
  B_{i} = B(X_{i},4r), \qquad i = 1,2,
\]
and let $U_{i}, V_{i}$ be the first and last lattice sites, 
respectively, of the segment of $\Gamma_{0} \cap B_{i}$
containing $X_{i}$.  Let $W_{1}$ be the first site in
$\Gamma_{0}^{[U_{1},X_{1}]}$ for which 
\[
  d(W_{1},\Gamma_{0}^{[X_{1},U_{2}]}) \leq q,
\]
and let $Z_{1}$ be the closest site to $W_{1}$ in 
$\Gamma_{0}^{[X_{1},X_{2}]}$.  $W_{2}$ and $Z_{2}$ are 
defined similarly with subscripts 1 and 2 interchanged.

Suppose now that $\Gamma_{0}$ is $(q,r)$-bottleneck-free
and satisfies
\begin{equation} \label{E:approxWulff}
  \Delta_{A}(\Gamma_{0}) \leq 2K_{63}l^{2/3}(\log l)^{1/3}
\end{equation}
for $K_{63}$ of Lemma \ref{L:inward}.
Then $B_{1}$ and $B_{2}$ are disjoint, and since 
\[
  \diam(\Gamma_{0}^{[U_{1},X_{1}]}) > r \quad \text{and} \quad
  \diam(\Gamma_{0}^{[X_{1},V_{1}]}) > r,
\]
the absense of bottlenecks implies 
\[
  d(\Gamma_{0}^{[X_{2},U_{1}]},\Gamma_{0}^{[X_{1},U_{2}]}) 
  > q \quad \text{and} \quad
  d(\Gamma_{0}^{[V_{2},X_{1}]},\Gamma_{0}^{[V_{1},X_{2}]}) 
  > q.
\]
It follows that
$W_{i} \neq U_{i}, Z_{i} \in \Gamma_{0}^{[X_{i},V_{i}]}$ 
($i = 1,2$) and
\[
  d(\Gamma_{0}^{[X_{2},W_{1}]},\Gamma_{0}^{[X_{1},W_{2}]}) 
  > q - 1 > \frac{q}{2}.
\]
When $\rho$ and $\sigma$ are paths such that the endpoint of 
$\rho$ is the initial point of $\sigma$, we let $(\rho,\sigma)$
denote the path obtained by concatenating $\sigma$ and $\rho$.
Then
\[
  |\mu_{L}(\Gamma_{0}^{[X_{2},X_{1}]}) - 
  \mu_{L}((\Gamma_{0}^{[X_{2},W_{1}]},\zeta_{W_{1}X_{1}}))| 
  \leq |B_{1}| = 16 \pi r^{2},
\]
since the paths differ only
inside $B_{1}$.  Again we may interchange 1, 2 and $L, R$.  Let $K_{66} > 0$
and let $D_{x} = B(x,K_{66}l^{2/3}(\log l)^{1/3})$
for $x \in \mathbb{R}^{2}$.  Presuming $K_{66}$
is large enough, (\ref{E:approxWulff}) implies that 
\[
  \Gamma_{y}^{[X_{1},X_{2}]} \backslash (D_{X_{1}} \cup D_{X_{2}})
  \subset J_{L}(\overline{X_{1}X_{2}}), \quad 
  \Gamma_{y}^{[X_{2},X_{1}]} \backslash (D_{X_{1}} \cup D_{X_{2}})
  \subset J_{R}(\overline{X_{1}X_{2}}).
\]
Hence using Lemma \ref{L:decouple},
\begin{align} \label{E:decomp}
  P&(|\Int(\Gamma_{0})| \geq 
    A - u, \Gamma_{0} \text{ is } 
    (q,r)-\text{bottleneck-free, } \\
  &\qquad \qquad \Delta_{A}(\Gamma_{0})
    \leq 2K_{63}l^{2/3}(\log l)^{1/3}, 0 \in \Mid(\Gamma_{0})) \notag \\
  &\leq \sum_{A-u \leq A_{L}+A_{R} \leq 2A} \quad \sum_{x_{1},x_{2} 
    \in B(y,\sqrt{A})} \quad
    \sum_{w_{1},z_{1} \in B(x_{1},4r)} \quad
    \sum_{w_{2},z_{2} \in B(x_{2},4r)} \notag \\
  &\qquad \qquad P(X_{1} = x_{1},
    X_{2} = x_{2}, W_{1} = w_{1}, W_{2} = w_{2},
    Z_{1} = z_{1}, Z_{2} = z_{2} \notag \\
  &\qquad \qquad \text{ and there exist paths } 
    \rho_{R} \text{ from } x_{2}
    \text{ to } w_{1} \text{ and } \rho_{L} \text{ from } w_{2}
    \text{ to } x_{1} \notag \\
  &\qquad \qquad \text{ satisfying } \mu_{L}((\rho_{R},
    \zeta_{w_{1}x_{1}})) \geq A_{R} - 16\pi r^{2}, 
     \notag \\
  &\qquad \qquad \mu_{R}((\zeta_{x_{2}w_{2}},
    \rho_{L})) \geq A_{L} - 16\pi r^{2},
    d(\rho_{L},\rho_{R}) \geq \frac{q}{2}, z_{1} \in 
    \rho_{L}, z_{2} \in \rho_{R}, \notag \\
  &\qquad \qquad \rho_{L} \backslash (D_{x_{1}} \cup D_{x_{2}})
    \subset J_{L}(\overline{x_{1}x_{2}}), \rho_{R} 
    \backslash (D_{x_{1}} \cup D_{x_{2}})
    \subset J_{R}(\overline{x_{1}x_{2}}), \notag \\
  &\qquad \qquad 0 \in J_{R}((\zeta_{x_{2}w_{2}},
    \rho_{L})) \cap J_{L}((\rho_{R},\zeta_{w_{1}x_{1}})) \backslash
    (D_{x_{1}} \cup D_{x_{2}})) \notag \\
  &\leq \sum_{A-u \leq A_{L}+A_{R} \leq 2A} \quad \sum_{x_{1},x_{2} 
    \in B(y,\sqrt{A})} \quad
    \sum_{w_{1},z_{1} \in B(x_{1},4r)} \quad
    \sum_{w_{2},z_{2} \in B(x_{2},4r)} \notag \\
  &\qquad \qquad 2P(x_{2} \lra w_{1} \text{ via an open dual path }
    \rho_{R} \notag \\
  &\qquad \qquad \qquad \text{ with } \mu_{L}((\rho_{R},
    \zeta_{w_{1}x_{1}})) \geq A_{R} - 16\pi r^{2}, z_{2} \in
    \rho_{R}, \notag \\
  &\qquad \qquad \qquad \rho_{R} 
    \backslash (D_{x_{1}} \cup D_{x_{2}})
    \subset J_{R}(\overline{x_{1}x_{2}}), \notag \\
  &\qquad \qquad \qquad 0 \in J_{L}((\rho_{R},
    \zeta_{w_{1}x_{1}})) \backslash (D_{x_{1}} \cup D_{x_{2}})) \notag \\
  &\qquad \qquad \cdot P(w_{2} \lra x_{1} 
    \text{ via an open dual path }
    \rho_{L} \notag \\
  &\qquad \qquad \qquad  \text{ with } \mu_{R}((\zeta_{x_{2}w_{2}},
    \rho_{L})) \geq A_{L} - 16\pi r^{2}, z_{1} \in
    \rho_{L}, \notag \\
  &\qquad \qquad \qquad \rho_{L} 
    \backslash (D_{x_{1}} \cup D_{x_{2}})
    \subset J_{L}(\overline{x_{1}x_{2}}), \notag \\
  &\qquad \qquad \qquad 0 \in J_{R}((\zeta_{x_{2}w_{2}},
    \rho_{L})) \backslash (D_{x_{1}} \cup D_{x_{2}})) \notag
\end{align}
Let us assume for convenience that $\delta$ is an integer
(if not, the necessary modifications are simple), and let
$x_{1}^{\prime}, w_{1}^{\prime}, z_{2}^{\prime}$ and $x_{2}^{\prime}$ 
be the lattice sites which
are $2\delta$ units to the right of $x_{1},w_{1},z_{2}$ and 
$x_{2}$, respectively.  We now ``pull apart'' the two halves
of $\Gamma_{0}$ by replacing each of these four sites by its 
right-shifted counterpart in the first
probability on the right side of (\ref{E:decomp}).  Specifically, by the FKG
property,
\begin{align} \label{E:FKG}
  P(&x_{2}^{\prime} \lra w_{1}^{\prime} \text{ via an open dual path }
    \rho_{R} \\
  &\qquad \qquad \text{ with } \mu_{L}((\rho_{R},
    \zeta_{w_{1}^{\prime}x_{1}^{\prime}})) \geq A_{R} 
    - 16\pi r^{2}, z_{2}^{\prime} \in \rho_{R}, \notag \\
  &\qquad \qquad \rho_{R} 
    \backslash (D_{x_{1}^{\prime}} \cup D_{x_{2}^{\prime}})
    \subset J_{R}(\overline{x_{1}^{\prime}x_{2}^{\prime}}), \notag \\
  &\qquad \qquad 0 \in J_{L}((\rho_{R},
    \zeta_{w_{1}^{\prime}x_{1}^{\prime}})) \backslash 
    (D_{x_{1}^{\prime}} \cup D_{x_{2}^{\prime}})) \notag \\
  &\quad \cdot P(w_{2} \lra x_{1} 
    \text{ via an open dual path }
    \rho_{L} \notag \\
  &\qquad \qquad \text{ with } \mu_{R}((\zeta_{x_{2}w_{2}},
    \rho_{L})) \geq A_{L} - 16\pi r^{2}, z_{1} \in
    \rho_{L}, \notag \\
  &\qquad \qquad \rho_{L} 
    \backslash (D_{x_{1}} \cup D_{x_{2}})
    \subset J_{L}(\overline{x_{1}x_{2}}), \notag \\
  &\qquad \qquad 0 \in J_{R}((\zeta_{x_{2}w_{2}},
    \rho_{L})) \backslash (D_{x_{1}} \cup D_{x_{2}})) \notag \\
  &\quad \cdot P(\Open(\zeta_{z_{1}w_{1}^{\prime}}))
    P(\Open(\zeta_{w_{2}z_{2}^{\prime}})) \notag \\
  &\leq P(|\Gamma_{0}| \geq A_{L} + A_{R} 
    + \frac{3}{2}\delta\sqrt{A} - 32\pi r^{2}
    - |D_{x_{1}}| - |D_{x_{2}}| - |D_{x_{1}^{\prime}}|
    - |D_{x_{2}^{\prime}}|) \notag \\
  &\leq P(|\Gamma_{0}| \geq A + \delta\sqrt{A}). \notag
\end{align}
From (\ref{E:M0bound3}), (\ref{E:decomp}), (\ref{E:FKG}) and
the bounded energy property we obtain
\begin{align} \label{E:combine}
  P(&|\Int(\Gamma_{0})| \geq A) \\
  &\leq K_{67}A^{4}r^{4}ue^{K_{68}\delta}
    P(|\Int(\Gamma_{0})| \geq A + \delta \sqrt{A}) \notag \\
  &\leq e^{K_{69}\delta}
    P(|\Int(\Gamma_{0})| \geq A + \delta \sqrt{A}), \notag 
\end{align}
completing the proof.
\end{proof}

\begin{proof}[Proof of Propostion \ref{P:changeA2}]
Let $B = A + \delta\sqrt{A}$.
Since $[|\Int(\Gamma_{y})| \geq A]$ is a decreasing event,
we have
\begin{equation} \label{E:PvsPlam}
  P(|\Int(\Gamma_{y})| \geq A) \geq 
  P_{N,w}(|\Int(\Gamma_{y})| \geq A),
\end{equation}
so it is enough to show
\begin{equation} \label{E:PvsPlam2}
  P_{N,w}(|\Int(\Gamma_{y})| \geq B)
  \geq e^{-\epsilon_{25}\delta}P(|\Int(\Gamma_{y})| \geq B).
\end{equation}
Let
\[
  r^{+}_{i} = \sqrt{B} 
  + 2K_{63}il^{2/3}(\log l)^{1/3}), \quad
  r^{-} = \sqrt{B} - 2K_{63}l^{2/3}(\log l)^{1/3}).
\]
Fix $y \in \ZZ$ and let $z \in \ZZ$ with 
$y \in z + r^{-}\mathcal{K}_{1}$ and 
$z + \sqrt{(1 + \epsilon_{24})A}\mathcal{K}_{1} 
\subset \Lambda_{N}$, where
$\epsilon_{24}$ is to be specified.  
Then $\sqrt{A/2} < 
r^{-} < r^{+}_{5} < \sqrt{B}$ and 
\begin{equation} \label{E:separation}
  d(z+r_{5}^{+}\mathcal{K}_{1},\partial\Lambda_{N}) 
  \geq \tfrac{1}{8}\epsilon_{24}\sqrt{A},
\end{equation}
provided $\epsilon_{23}$ and $\epsilon_{24}$ are small and
$A$ is large enough. 

Let $\tilde{\Gamma}_{y}^{i}$ denote the outermost open dual 
circuit surrounding $y$ in $z + r^{+}_{i}\mathcal{K}_{1}$
and define the event
\begin{align} 
  F_{i} = [|\Int(\tilde{\Gamma}_{y}^{i})| &\geq B, 
    d_{H}(\partial\Co(\tilde{\Gamma}_{y}^{i}),z + 
    \partial\sqrt{B}\mathcal{K}_{1}) \leq 
    K_{4}l^{2/3}(\log l)^{1/3}, \notag \\
  &MLR(\tilde{\Gamma}_{y}^{i}) \leq K_{63}l^{2/3}(\log l)^{1/3}]
    \notag
\end{align}
Note that $F_{i} \subset [\tilde{\Gamma}_{y}^{i} \subset
z + r^{+}_{1}\mathcal{K}_{1}]$.  Let
$E_{i}$ be the event that
there exist $u \notin z + r^{+}_{i+1}\mathcal{K}_{1}, v 
\notin z + r^{+}_{i}\mathcal{K}_{1}$ for which $d(u,v) \geq
\tfrac{1}{2}d(u,z+r^{+}_{i}\mathcal{K}_{1})$ and 
$u \lra v$ via 
an open dual path outside $z + r^{+}_{i}\mathcal{K}_{1}$.  Then
for some $\epsilon_{28}$,
\begin{equation} \label{E:Eibound}
  P_{N,w}(E_{i}) \leq 
  P(E_{i}) \leq \exp(-\epsilon_{28}l^{2/3}(\log l)^{1/3}),
\end{equation}
\[
  F_{5} 
  \subset [\tilde{\Gamma}_{y}^{5} \subset z + 
  r^{+}_{1}\mathcal{K}_{1}] \subset 
  [\tilde{\Gamma}_{y}^{5} = \tilde{\Gamma}_{y}^{3}]
\]
and 
\[
  F_{5} \cap [\tilde{\Gamma}_{y}^{5} \neq \Gamma_{y}]
  \subset [\Gamma_{y} \not\subset z + r^{+}_{5}\mathcal{K}_{1},
  \tilde{\Gamma}_{y}^{5} \subset z + 
  r^{+}_{1}\mathcal{K}_{1}] \subset E_{4}.
\]
Here the last inclusion follows from the fact that if
$\Gamma_{y} \not\subset z + r^{+}_{5}\mathcal{K}_{1}$ and
$\tilde{\Gamma}_{y}^{5} \subset z + 
r^{+}_{1}\mathcal{K}_{1}$ then $\Gamma_{y}$ must 
surround or intersect $\tilde{\Gamma}_{y}^{5}$.
The events $E_{4}$ and $F_{3}$ necessarily
occur at separation $2K_{63}l^{2/3}(\log l)^{1/3}$, so
by Lemma \ref{L:ratiowm} we have
\begin{equation} \label{E:mixing1}
  P_{N,w}(F_{5} \cap 
  [\tilde{\Gamma}_{y}^{5} \neq \Gamma_{y}]) 
  \leq P_{N,w}(E_{4} \cap F_{3})
  \leq 2P_{N,w}(E_{4})P_{N,w}(F_{3}). 
\end{equation}
We want to replace $P(F_{3})$
with $P(F_{5})$ on the right side of (\ref{E:mixing1}).  We have
\begin{equation} \label{E:contain}
  F_{3} \backslash F_{5} \subset F_{3} \cap [\tilde{\Gamma}_{y}^{5}
  \not\subset z + r^{+}_{3}\mathcal{K}_{1}] \subset
  F_{3} \cap E_{2}.
\end{equation}
Let $\omega_{1} \in F_{1}$ be a bond configuration on
$\mathcal{B}(z + r^{+}_{1}\mathcal{K}_{1})$.  Conditionally
on $\omega_{1}$, $F_{3}$ is an increasing event (since $F_{3}$ requires
$\tilde{\Gamma}_{y}^{3} \subset z + r^{+}_{1}\mathcal{K}_{1}$,
meaning $\tilde{\Gamma}_{y}^{3}$ is part of $\omega_{1}$)
and $E_{2}$ is a 
decreasing one, so using (\ref{E:Eibound}), (\ref{E:contain}) and 
Lemma \ref{L:ratiowm} again,
\begin{align} \label{E:F3F5}
  P_{N,w}&(F_{3} \backslash F_{5} \mid \omega_{1}) \\
  &\leq P_{N,w}(F_{3} \mid \omega_{1}) 
    P_{N,w}(E_{2} \mid \omega_{1}) \notag \\
  &\leq P_{N,w}(F_{3} \mid \omega_{1})
    P(E_{2} \mid \omega_{1}) \notag \\
  &\leq 2P_{N,w}(F_{3} \mid \omega_{1})
    P(E_{2}) \notag \\
  &\leq \frac{1}{2}P_{N,w}(F_{3} \mid \omega_{1}). \notag
\end{align}
Therefore $P_{N,w}(F_{3}) \leq 2P_{N,w}(F_{5})$.
With (\ref{E:mixing1}), (\ref{E:Eibound}), (\ref{E:separation}) and Lemma
\ref{L:ratiowm} this shows that
\begin{align} \label{E:PvsPlam3}
  P(&|\Int(\Gamma_{y})| \geq B, 
    d_{H}(\partial\Co(\Gamma_{y}),z + 
    \partial\sqrt{B}\mathcal{K}_{1}) \leq 
    K_{4}l^{2/3}(\log l)^{1/3}, \\
  &\qquad \qquad MLR(\Gamma_{y}) \leq 
    K_{63}l^{2/3}(\log l)^{1/3}) \notag \\
  &\leq P(F_{5}) \notag \\
  &\leq 2P_{N,w}(F_{5}) \notag \\
  &\leq 4P_{N,w}(F_{5} \cap [\Gamma_{y} = 
    \tilde{\Gamma}_{y}^{5}]) \notag \\
  &\leq 4P_{N,w}(|\Int(\Gamma_{y})| \geq B, 
    d_{H}(\partial\Co(\Gamma_{y}),z + 
    \partial\sqrt{B}\mathcal{K}_{1}) \leq 
    K_{4}l^{2/3}(\log l)^{1/3}, \notag \\
  &\qquad \qquad \qquad MLR(\Gamma_{y}) \leq 
    K_{63}l^{2/3}(\log l)^{1/3}). \notag 
\end{align}

By Theorem \ref{T:main}, Lemma \ref{L:inward} and translation
invariance, there exists a site $y^{\prime}$ such that
\begin{align} \label{E:choosez}
  \frac{1}{2B}&P(|\Int(\Gamma_{y})| \geq B) \\
  &= \frac{1}{2B}P(|\Int(\Gamma_{y^{\prime}})| \geq B) \notag \\
  &\leq \frac{1}{B}P(|\Int(\Gamma_{y^{\prime}})| \geq B,
    \Delta_{B}(\partial\Co(\Gamma_{y^{\prime}})) \leq 
    K_{4}l^{2/3}(\log l)^{1/3}, \notag \\
  &\qquad \qquad MLR(\Gamma_{y^{\prime}}) \leq 
    K_{63}l^{2/3}(\log l)^{1/3}) \notag \\
  &\leq P(|\Int(\Gamma_{y^{\prime}})| \geq B, 
    d_{H}(\partial\Co(\Gamma_{y^{\prime}}),z + 
    \partial\sqrt{B}\mathcal{K}_{1}) \leq 
    K_{4}l^{2/3}(\log l)^{1/3}, \notag \\
  &\qquad \qquad MLR(\Gamma_{y^{\prime}}) \leq 
    K_{63}l^{2/3}(\log l)^{1/3}). \notag 
\end{align}
But the last event implies that $\Gamma_{y^{\prime}}$
surrounds $z + r^{-}\mathcal{K}_{1}$, which contains $y$.  Hence
the last probability
in (\ref{E:choosez}) is bounded by the first probability in
(\ref{E:PvsPlam3}), so that
\begin{equation} \label{E:PvsPlam4}
  \frac{1}{2B}P(|\Int(\Gamma_{y})| \geq B)
  \leq 4P_{N,w}(|\Int(\Gamma_{y})| \geq B). 
\end{equation}
Since $\delta \geq K_{57}\log A$, this 
completes the proof of (\ref{E:PvsPlam2}).
\end{proof}

The next lemma includes the
analog of Lemma \ref{L:inward} for $P_{N,w}$.

\begin{lemma} \label{L:inward2}
  Let $P$ be a percolation model on $\mathcal{B}(\ZZ)$ satisfying
  (\ref{E:assump}), the near-Markov property for open circuits,
  and the ratio weak mixing property.  There exist
  $\epsilon_{i},K_{i}$ such that for $N \geq 1,
  K_{70}(\log N)^{2} \leq A \leq c_{2}N^{2}$ (with $c_{2}$ from
  Theorem \ref{T:single}) and $l = \sqrt{A}$,
  we have
  \begin{align} \label{E:MLRbound3}
    P_{N,w}&(|\Int(\Gamma_{y})| \geq A \text{ and }
      MLR(\Gamma_{y}) > K_{71}l^{2/3}(\log l)^{1/3} \text{ for some }
      y \in \Lambda_{N} \cap \ZZ) \\
    &\leq \exp(-\epsilon_{29}l^{2/3}(\log l)^{1/3})P_{N,w}
      (|\Int(\Gamma_{y})| \geq A \text{ for some }
      y \in \Lambda_{N} \cap \ZZ), \notag
  \end{align}
  \begin{align} \label{E:Hausmax4}
    P_{N,w}&(|\Int(\Gamma_{y})| \geq A \text{ and }
      \Delta_{A}(\Gamma_{y}) > K_{72}l^{2/3}(\log l)^{1/3} \text{ for some }
      y \in \Lambda_{N} \cap \ZZ) \\
    &\leq \exp(-\epsilon_{30}l^{2/3}(\log l)^{1/3})P_{N,w}
      (|\Int(\Gamma_{y})| \geq A \text{ for some }
      y \in \Lambda_{N} \cap \ZZ), \notag
  \end{align}
  \begin{align} \label{E:ALRbound3}
    P_{N,w}&(|\Int(\Gamma_{y})| \geq A \text{ and }
      ALR(\Gamma_{y}) > K_{73}l^{1/3}(\log l)^{2/3} \text{ for some }
      y \in \Lambda_{N} \cap \ZZ) \\
    &\leq \exp(-\epsilon_{31}l^{1/3}(\log l)^{2/3})P_{N,w}
      (|\Int(\Gamma_{y})| \geq A \text{ for some }
      y \in \Lambda_{N} \cap \ZZ), \notag
  \end{align}
  and for $\epsilon_{32}A \geq r \geq 15q \geq K_{74}\log A$ and 
  $\kt r/3 \leq d^{\prime} \leq \sqrt{A}$,
  \begin{align} \label{E:Dprimebd}
    P_{N,w}&(|\Int(\Gamma_{y})| \geq A,
      D_{(q,r)}^{\prime}(\Gamma_{y}) > d^{\prime}
      \text{ and } \\
    &\qquad \qquad |\Int(\Gamma_{y})| - |\Int(\alpha_{\max,\Gamma_{y}})|
      < \epsilon_{33}d^{\prime}\sqrt{A} \text{ for some } y 
      \in \Lambda_{N} \cap \ZZ) \notag \\
    &\leq \exp(-\epsilon_{34}d^{\prime})P_{N,w}
      (|\Int(\Gamma_{y})| \geq A \text{ for some }
      y \in \Lambda_{N} \cap \ZZ). \notag
  \end{align}
\end{lemma}
\begin{proof}
We begin with (\ref{E:MLRbound3}).
The proof of Lemma
\ref{L:inward} is valid for $P_{N,w}$ through (\ref{E:M0bound})
(see Remark \ref{R:wired}), which gives that for $c$ sufficiently large,
\begin{align} \label{E:M0boundsum}
  P_{N,w}&(|\Int(\Gamma_{y})| \geq A, MLR(\Gamma_{y}) \geq 
    52cl^{2/3}(\log l)^{1/3} \text{ for some } y \in \Lambda_{N} \cap \ZZ) \\
  &\leq |\Lambda_{N} \cap \ZZ|\exp\left(-w_{1}\sqrt{A} 
    - \frac{1}{80}c\kt l^{2/3}(\log l)^{1/3}\right). \notag 
\end{align}
If $l^{2/3}(\log l)^{1/3} \geq \log N$, we can complete the proof as
we did in Lemma \ref{L:inward}.  Thus suppose
\begin{equation} \label{E:lvsN}
  l^{2/3}(\log l)^{1/3} < \log N.
\end{equation}
Since $A \leq c_{2}N^{2}$, provided we choose $\epsilon_{35}$ small enough
the set $J_{A} = \{y \in \Lambda_{N} \cap \ZZ: y \text{ is } (\Lambda_{N},
A/2,(1 + \epsilon_{35})A)-\text{compatible}\}$ 
satisfies $|J_{A}| \geq \epsilon_{36}
N^{2}$.  Let $E_{y}$ denote the event that 
$A \leq |\Int(\Gamma_{y})| \leq 2A$ and 
$\Delta_{A}(\partial\Co(\Gamma_{y})) \leq
K_{4} l^{1/3}(\log l)^{2/3}$.  Let $F_{y}$ denote the event 
that $\Gamma_{y}$ is the 
unique exterior dual circuit in $\Lambda_{N}$ satisfying both
$A \leq |\Int(\Gamma_{x})| \leq 2A$ for some $x \in J_{A}$
and $\Delta_{A}(\partial\Co(\Gamma_{x})) \leq K_{4}
l^{1/3}(\log l)^{2/3}$.  From Theorem \ref{T:main}, Theorem
\ref{T:lowerbound} and Remark \ref{R:wired}, we have for $y \in J_{A}$
\begin{equation} \label{E:badshape}
  P_{N,w}(E_{y}^{c} \mid |\Int(\Gamma_{y})| \geq A) 
  \leq \frac{1}{4}.
\end{equation}
Also, from the near-Markov property for open circuits, the FKG property,
and Theorem \ref{T:onebound}, provided $K$ is large,
\begin{align} \label{E:notunique}
  P_{N,w}&(F_{y}^{c} \mid E_{y}) \\
  &= \sum_{\nu} 
    P_{N,w}(F_{y}^{c} \mid E_{y} \cap [\Theta(\Gamma_{y})
    = \nu])P_{N,w}(\Theta(\Gamma_{y})
    = \nu \mid E_{y}) \notag \\
  &\leq \sum_{\nu} 
    2P_{N,w}(|\Int(\Gamma_{x})| \geq A \text{ for some }
    x \in J_{A} \cap \Ext(\nu)) P_{N,w}(\Theta(\Gamma_{y})
    = \nu \mid E_{y}) \notag \\
  &\leq 2P_{N,w}(|\Int(\Gamma_{x})| \geq A \text{ for some }
    x \in J_{A}) \notag \\
  &\leq \frac{1}{4}, \notag
\end{align}
so by Theorem \ref{T:lowerbound} and Remark \ref{R:wired}, for 
$y \in J_{A}$, using the $(\Lambda_{N},
A/2,(1 + \epsilon_{35})A)$-compatibility of $y$,
\[
  P_{N,w}(F_{y}) \geq 
  \frac{1}{2}P_{N,w}(|\Int(\Gamma_{y})| \geq A) \geq
  \exp(-w_{1}\sqrt{A} - K_{76}l^{1/3}(\log l)^{2/3}).
\]
Combining these facts we obtain
\begin{align} \label{E:sumlower}
  2A&P_{N,w}(|\Int(\Gamma_{y})| \geq A \text{ for some }
    y \in J_{A}) \\
  &\geq 2AP_{N,w}(\cup_{y \in J_{A}}F_{y}) \notag \\
  &\geq \sum_{y \in J_{A}} P_{N,w}(F_{y}) \notag \\
  &\geq \epsilon_{36}N^{2}\exp(-w_{1}\sqrt{A} 
    - K_{76}l^{1/3}(\log l)^{2/3}). \notag 
\end{align}
Here the second inequality uses the fact that $F_{y}$ can occur 
in each bond configuration for 
at most $2A$ sites $y \in J_{A}$.
With (\ref{E:lvsN}) and (\ref{E:M0boundsum}), provided we take $c$
large enough, (\ref{E:sumlower}) gives
\begin{align} 
  \exp&(-l^{2/3}(\log l)^{1/3})
  P_{N,w}(|\Int(\Gamma_{y})| \geq A \text{ for some }
    y \in \Lambda_{N} \cap \ZZ) \notag \\
  &\geq |\Lambda_{N} \cap \ZZ|\exp(-w_{1}\sqrt{A} 
    - K_{77}l^{2/3}(\log l)^{1/3}) \notag \\
  &\geq P_{N,w}(|\Int(\Gamma_{y})| \geq A, MLR(\Gamma_{y}) 
    \geq 52cl^{2/3}(\log l)^{1/3} \text{ for some }
    y \in \Lambda_{N} \cap \ZZ), \notag 
\end{align}
which completes the proof of (\ref{E:MLRbound3}).

The proofs of (\ref{E:Hausmax4}) and (\ref{E:ALRbound3})
are essentially the same, except that for (\ref{E:Hausmax4}),
in place of (\ref{E:M0boundsum}) we have the following.  
From Remark \ref{R:wired} and the proof of
(\ref{E:Hausmax}) of Theorem \ref{T:main}, for some $K_{78}$ and for
$c$ sufficiently large:
\begin{align} 
  P_{N,w}&(|\Int(\Gamma_{y})| \geq A, 
    \Delta_{A}(\partial\Co(\Gamma_{y})) 
    \geq cl^{1/3}(\log l)^{2/3} \text{ for some } y \in 
    \Lambda_{N} \cap \ZZ) \\
  &\leq |\Lambda_{N} \cap \ZZ|\exp(-w_{1}\sqrt{A} 
    - K_{78}c\kt l^{1/3}(\log l)^{2/3}). \notag 
\end{align}
With (\ref{E:DeltaA}) and (\ref{E:MLRbound3}) this proves
(\ref{E:Hausmax4}).
The modification is similar for (\ref{E:ALRbound3}).

For (\ref{E:Dprimebd}) we must do more, because $d^{\prime}$ 
may be much smaller than the order of the error term
$l^{1/3}(\log l)^{2/3}$ in Theorem \ref{T:lowerbound}.  Let $a =
\epsilon_{33}d^{\prime}\sqrt{A}$, with $\epsilon_{33}$ still to be specified.
From Theorem \ref{T:onebound}, Remark \ref{R:wired} and
(\ref{E:sumlower}) we have
\begin{align} \label{E:verybig}
  P_{N,w}&(|\Int(\Gamma_{y})| \geq 2A,
    D_{(q,r)}^{\prime}(\Gamma_{y}) \geq d^{\prime} \text{ for some }  
    y \in \Lambda_{N} \cap \ZZ) \\
  &\leq |\Lambda_{N} \cap \ZZ|\exp\left(-\frac{1}{20}d^{\prime} 
    - u(K_{52}(\log 2A)^{2/3},2A)\right) \notag \\
  &\leq \exp\left(-\frac{1}{20}d^{\prime}\right)
    P_{N,w}(|\Int(\Gamma_{y})| \geq A \text{ for some }  
    y \in \Lambda_{N} \cap \ZZ). \notag
\end{align}
From Proposition \ref{P:withbnbound} and Remark \ref{R:wired},
similarly to (\ref{E:M0bound}) we have
\begin{align} \label{E:regsize}
  P_{N,w}&\bigl(A \leq |\Int(\Gamma_{y})| < 2A,
    D_{(q,r)}^{\prime}(\Gamma_{y}) > d^{\prime}
    \text{ and } \\
  &\qquad \qquad |\Int(\Gamma_{y})| - |\Int(\alpha_{\max,\Gamma_{y}})|
    < a \text{ for some } y \in \Lambda_{N} \cap \ZZ \bigr) \notag \\
  &\leq \sum_{y \in \Lambda_{N} \cap \ZZ} \sum_{A \leq B < 2A} 
    \sum_{n \geq d^{\prime}} 
    \sum_{B - a < A^{\prime} \leq B} \sum_{k \leq n} 
    P_{N,w}(M_{y}(k,q,r,B,A^{\prime},n,w_{1}\sqrt{A^{\prime}}))
    \notag \\
  &\leq \sum_{y \in \Lambda_{N} \cap \ZZ} 
    \sum_{A - a \leq A^{\prime} < 2A} 
    \sum_{A^{\prime} \leq B < A^{\prime} + a} \sum_{n \geq d^{\prime}} 
    \notag \\
  &\qquad \qquad (n + 1)e^{-n/2}P_{N,w}
    (\cup_{x \in B(y,K_{58}A)}
    M_{x}(0,q,r,A^{\prime},A^{\prime},0,
    w_{1}\sqrt{A^{\prime}})) \notag \\
  &\leq \sum_{y \in \Lambda_{N} \cap \ZZ} 
    K_{79}a(d^{\prime} + 1)e^{-d^{\prime}/2}
    \notag \\
  &\qquad \qquad \cdot
    P_{N,w}(A - a \leq |\Int(\Gamma_{x})| < 2A \text{ for some } 
    x \in B(y,K_{58}A)) \notag \\
  &\leq K_{80}a(d^{\prime} + 1)e^{-d^{\prime}/2}A^{2}
    \sum_{x \in \Lambda_{N} \cap \ZZ}
    P_{N,w}(A - a \leq |\Int(\Gamma_{x})| < 2A). \notag 
\end{align}
Let $\tilde{F}_{y}$ denote the event that $A - a \leq |\Int(\Gamma_{y})|
< 2A$ and $\Gamma_{y}$ is the unique 
exterior dual circuit in $\Lambda_{N}$ satisfying  
$|\Int(\Gamma_{y})| \geq A - a$.  Similarly to (\ref{E:notunique})
and (\ref{E:sumlower}) we have
\begin{align} \label{E:unique2}
  &\sum_{x \in \Lambda_{N} \cap \ZZ} P_{N,w}(A - a \leq 
    |\Int(\Gamma_{x})| < 2A) \\
  &\quad \leq \sum_{x \in \Lambda_{N} \cap \ZZ} 
    2P_{N,w}(\tilde{F}_{x})
    \notag \\
  &\quad \leq 4AP_{N,w}(A - a \leq |\Int(\Gamma_{x})| <2A 
    \text{ for some } x \in \Lambda_{N} \cap \ZZ). \notag
\end{align}
From (\ref{E:Hausmax4}) we have 
\begin{align} \label{E:goodcirc}
  P_{N,w}&(|\Int(\Gamma_{x})| \geq A - a 
    \text{ for some } x \in \Lambda_{N} \cap \ZZ) \\
  &\leq 2\sum_{x \in \Lambda_{N} \cap \ZZ}
    P_{N,w}(|\Int(\Gamma_{x})| \geq A - a,
    \Delta_{A}(\Gamma_{x}) 
    \leq K_{72}l^{2/3}(\log l)^{1/3}). \notag
\end{align}
But the last event, which says roughly that $\Gamma_{x}$ approximates
the appropriate Wulff shape, implies that $\Gamma_{x}$ surrounds 
some $(\Lambda_{N},(A-a)/2,(1+\epsilon_{24})(A-a))$-compatible site 
$y$. (Take $y$ to be the closest site to 0 in 
$\Int(\Gamma_{x})$.)  Here $\epsilon_{24}$ is from 
Proposition \ref{P:changeA2}.  Thus provided $\epsilon_{33}$
is small, (\ref{E:goodcirc}) and Proposition
\ref{P:changeA2} give
\begin{align} \label{E:AatoA}
  P_{N,w}&(|\Int(\Gamma_{x})| \geq A - a 
    \text{ for some } x \in \Lambda_{N} \cap \ZZ) \\
  &\leq 2e^{\epsilon_{25}\epsilon_{33}d^{\prime}}
    \sum_{x \in \Lambda_{N} \cap \ZZ \cap (\ZZ)^{*}}
    P_{N,w}(|\Int(\Gamma_{y})| \geq A \text{ for some }
    y \in B(x,\sqrt{A}) \cap J_{A-a}) \notag \\
  &\leq 2K_{81}Ae^{\epsilon_{25}\epsilon_{33}d^{\prime}}\sum_{y \in J_{A-a}}
    P_{N,w}(|\Int(\Gamma_{y})| \geq A). \notag 
\end{align}
We can repeat the argument 
of (\ref{E:badshape}) -- (\ref{E:sumlower}) (excluding the last
inequality of (\ref{E:sumlower})) once more to get
\begin{equation} \label{E:sumout}
  \sum_{y \in J_{A-a}}
    P_{N,w}(|\Int(\Gamma_{y})| \geq A)
  \leq 4AP_{N,w}(|\Int(\Gamma_{y})| \geq A 
    \text{ for some } y \in \Lambda_{N} \cap \ZZ).
\end{equation}
Provided $\epsilon_{33}$ is small, together, (\ref{E:verybig}),
(\ref{E:regsize}), (\ref{E:unique2}), (\ref{E:AatoA}) and 
(\ref{E:sumout}) prove the lemma.
\end{proof}

\begin{proof}[Proof of Theorem \ref{T:single}]
We begin by obtaining an analog of (\ref{E:mkbound3}), by 
mimicking its proof.  We omit some details becuase of
the similarity.  First fix $N$ and let $K_{82}$ be a constant
to be specified.
Our induction hypothesis is that for
every $j < k, K_{82} \leq A^{\prime} \leq A, d \geq 
d^{\prime} \geq 0$
and enclosure event $E$, we have 
\begin{align} \label{E:indhyp1}
  P_{N,w}&(G_{N}(j,A,A^{\prime},d,d^{\prime}) 
    \mid E) \\
  &\leq e^{-\frac{4}{5}d^{\prime}}P_{N,w}
    (|\Int(\Gamma_{y})| = A^{\prime} \text{ for some }
    y \in \Lambda_{N} \cap \ZZ \mid E) \notag
\end{align}
and
\begin{equation} \label{E:indhyp2}
  P_{N,w}(G_{N}(j,A,A^{\prime},d,d^{\prime}) \mid E)
  \leq \exp\left(-\frac{9}{10}d + \frac{1}{20}jK\log N\right)
\end{equation}
and
\begin{equation} \label{E:indhyp3}
  P_{N,w}(G_{N}(j,A,A^{\prime},d,d^{\prime}) \mid E)
  \leq \exp\left(-u(K_{52}(\log N)^{2/3},A) - \frac{1}{20}d^{\prime}\right).
\end{equation}
We need only consider parameters satisfying
\begin{equation} \label{E:params}
  d + 1 \geq (j+1)K\log N, \qquad d^{\prime} 
  + 1 \geq jK\log N, \quad 
  d^{\prime} + 1 \geq \frac{1}{6}w_{1}\sqrt{A - A^{\prime}},
\end{equation}
the last inequality following from (\ref{E:Ddiam}) and subadditivity of
the square root.
For $j = 0$, we need only consider
$A = A^{\prime}, d^{\prime} = 0$ and (\ref{E:indhyp1}) is
immediate from the definitions, while (\ref{E:indhyp2}) 
follows easily from the first inequality in (\ref{E:biggamma})
and the fact that $E$ is increasing, and (\ref{E:indhyp3})
follows from Theorem \ref{T:onebound}.
Hence we may assume $k \geq 1$ and fix 
$A,A^{\prime},d,d^{\prime},E$.  

Let $\omega \in G_{N}(k,A,A^{\prime},d,d^{\prime})$.  There may be open 
dual paths $\alpha$ in $\omega$ each connecting 
$\Phi_{N}$ to another open dual 
circuit $\gamma \in \mathfrak{C}_{N}$.  Since $\Phi_{N}$ is exterior,
for each such $\gamma$ there is a unique bond $\langle xy \rangle$ with
$x \in \Phi_{N}, y \in \Ext(\Phi_{N})$ which is part of such an $\alpha$.
We denote the set of all such bonds by $\mathcal{A}_{N}$.  Thus
$|\mathcal{A}_{N}| \leq |\mathfrak{C}_{N}| - 1$.  Let 
\[
  \tilde{G}_{N}(k,A,A^{\prime},d,d^{\prime}) = 
  G_{N}(k,A,A^{\prime},d,d^{\prime}) \cap [\mathcal{A}_{N} = \phi].
\]
When this event occurs, the only $(K\log N)$-large open dual circuit
in $\Int(\Theta(\Phi_{N}))$ is $\Phi_{N}$.
From the above and the bounded energy property we have
\[
  P_{N,w}(G_{N}(k,A,A^{\prime},d,d^{\prime}) \mid E) \leq
  |\mathcal{B}(\Lambda_{N})|^{k}p_{0}^{-k}
  P_{N,w}(\tilde{G}_{N}(k,A,A^{\prime},d,d^{\prime}) \mid E),
\]
where $p_{0}$ is from (\ref{E:boundener}).
Given a circuit $\nu \subset \Lambda_{N}$ define the event
\begin{align} \label{E:onlyone}
  F(\nu,d_{1}) = &[|\Int(\Gamma_{y})| = A^{\prime}, d_{1} \leq 
    \diam_{\tau}(\Gamma_{y}) < d_{1} + 1, \Theta(\Gamma_{y}) = \nu
    \text{ and } \\
  &\quad \Gamma_{y} \text{ is the only } 
    (K\log N)-\text{large open dual circuit in } \Int(\nu), \notag \\
  &\quad \text{ for some } y \in \Int(\nu)]. \notag 
\end{align}
Note that $F(\nu,d_{1}) \in \mathcal{G}_{\mathcal{B}(\Int(\nu) \cup \nu)}$.
It follows easily from the near-Markov property and the first inequality
in (\ref{E:biggamma}) that provided $K$ is 
large, given $F(\nu,d_{1}) \cap E$ with high 
probability there are no $(K\log N)$-large open dual circuits outside
$\nu$; more precisely,
\[
  P_{N,w}\bigl(G_{N}(0,A^{\prime},A^{\prime},d_{1},0)
    \cap [\Theta(\Phi_{N}) = \nu] \mid E \cap F(\nu,d_{1})
    \bigr) \geq \frac{1}{2}.
\]
Hence as in the proof of Proposition
\ref{P:withbnbound}, conditioning on $\Theta(\Phi_{N})$ and using
the near-Markov property and (\ref{E:indhyp2}) for $j = k - 1$ gives
\begin{align} \label{E:GNbound1}
  P_{N,w}&(G_{N}(k,A,A^{\prime},d,d^{\prime}) \mid E) \\
  &\leq |\mathcal{B}(\Lambda_{N})|^{k}p^{-k}
    P_{N,w}(\tilde{G}_{N}(k,A,A^{\prime},d,d^{\prime}) \mid E)
    \notag \\
  &\leq K_{83}p^{-k}N^{2k}
    \sum_{\nu} P_{N,w}(\tilde{G}_{N}(k,A,A^{\prime},d,d^{\prime}) 
    \cap [\Theta(\Phi_{N}) = \nu] \mid E) \notag \\
  &\leq 2K_{83}p^{-k}N^{2k} \sum_{d_{1}} 
    \quad \sum_{\nu} P_{N,w}(F(\nu,d_{1}) \mid E) \notag \\
  &\qquad \cdot 
    P_{N,w}(\tilde{G}_{N}(k,A,A^{\prime},d,d^{\prime}) 
    \cap [\Theta(\Phi_{N}) = \nu] \mid E \cap F(\nu,d_{1})) \notag \\
  &\leq 8K_{83}p^{-k}N^{2k} \sum_{d_{1},d_{2}^{\prime}} 
    \quad \sum_{0 < A_{1} \leq A-A^{\prime}} \quad
    \sum_{\nu} \notag \\
  &\qquad P_{N,w}(G_{N}(0,A^{\prime},
    A^{\prime},d_{1},0) \cap
    [\Theta(\Phi_{N}) = \nu] \mid E) \notag \\
  &\qquad \cdot 
    P_{N,w}(G_{N}(k-1,A-A^{\prime},A_{1},d^{\prime},d_{2}^{\prime}) 
    \mid E \cap \Open(\nu) \cap [\nu \lra \infty]) \notag \\
  &\leq 8K_{83}p^{-k}N^{2k+6} \sum_{d_{1}} \sum_{\nu}
    P_{N,w}(G_{N}(0,A^{\prime},A^{\prime},d_{1},0) \cap
    [\Theta(\Phi_{N}) = \nu] \mid E) \notag \\
  &\qquad \qquad \qquad \cdot \exp\left(-\frac{9}{10}d^{\prime} 
    + \frac{1}{20}(k-1)K\log N\right) \notag \\
  &\leq 8K_{83}p^{-k}N^{2k+6} \exp\left(-\frac{9}{10}d^{\prime} 
    + \frac{1}{20}(k-1)K\log N\right) \notag \\
  &\qquad \cdot P_{N,w}
    (|\Int(\Gamma_{y})| = A^{\prime}, d - d^{\prime} - 1 \leq 
    \diam_{\tau}(\Gamma_{y}) \notag \\
  &\qquad \qquad \qquad < d - d^{\prime} + 1
    \text{ for some }
    y \in \Lambda_{N} \cap \ZZ \mid E) \notag \\
  &\leq \exp\left(-\frac{4}{5}d^{\prime}\right) P_{N,w}
    (|\Int(\Gamma_{y})| = A^{\prime}, d - d^{\prime} - 1 \leq 
    \diam_{\tau}(\Gamma_{y}) \notag \\
  &\qquad \qquad \qquad \qquad \qquad < d - d^{\prime} + 1
    \text{ for some } y \in \Lambda_{N} \cap \ZZ \mid E). \notag
\end{align}
Here the sums include $d_{1} \in \{d-d^{\prime}-1,d-d^{\prime}\}$
and $0 \leq d_{2}^{\prime} < d^{\prime}$.  Now (\ref{E:GNbound1}) 
establishes (\ref{E:indhyp1}) for $j = k$; we next establish
(\ref{E:indhyp2}).  We have by (\ref{E:conupr})
\begin{align}
  P_{N,w}&(d - d^{\prime} - 1 \leq 
    \diam_{\tau}(\Gamma_{y}) < d - d^{\prime} + 1
    \text{ for some } y \in \Lambda_{N} \cap \ZZ \mid E) \notag \\
  &\leq |\Lambda_{N} \cap \ZZ|^{2}e^{-(d - d^{\prime} - 1)}.
    \notag 
\end{align}
Hence provided $K$ is large
enough, the right side of the next-to-last inequality in 
(\ref{E:GNbound1}) is bounded by 
\[
  \exp\left(-\frac{9}{10}d + \frac{1}{20}kK\log N\right),
\]
which proves (\ref{E:indhyp2}) for $j = k$.  Turning to 
(\ref{E:indhyp3}), from (\ref{E:indhyp2}),
(\ref{E:indhyp3}) for $j = k - 1$, and Theorem \ref{T:main},
the right side of the fourth inequality in
(\ref{E:GNbound1}) is bounded by
\begin{align} \label{E:GNbound3}
  &8K_{84}p^{-k}N^{2k+4} \exp(-w_{1}\sqrt{A^{\prime}} + 
    K_{52}(\log N)^{2/3}(A^{\prime})^{1/6})) \\
  &\qquad \qquad \cdot \exp\biggl(- \max\biggl[w_{1}\sqrt{A-A^{\prime}} - 
    K_{52}(\log N)^{2/3}(A-A^{\prime})^{1/6}), \notag \\
  &\qquad \qquad \qquad \qquad \qquad \frac{9}{10}
    d^{\prime} - \frac{1}{20}(k-1)K\log N\biggr]\biggr). \notag
\end{align}
We now have two cases.

\emph{Case 1}. $A^{\prime} \geq \tfrac{1}{3}A$.
From (\ref{E:params}) we have
\[
  d^{\prime} \geq \frac{1}{2}kK\log N + (K\log N)^{2/3}(\frac{1}{12}
  w_{1}\sqrt{A^{\prime}})^{1/3}.
\]
We can use this, (\ref{E:sqroot})
with $\theta = 1/\sqrt{2}$, (\ref{E:indhyp3}) with
$j = k-1$,
and the fact that the maximum exceeds any
convex combination to conclude that, provided $K$ is sufficiently
large, (\ref{E:GNbound3}) is bounded above by
\begin{align} \label{E:GNbound4}
  &8K_{84}p^{-k}N^{2k+4} \exp\biggl(-w_{1}\sqrt{A^{\prime}} - 
    \frac{3}{4}w_{1}\sqrt{A - A^{\prime}} - \frac{9}{40}d^{\prime}
    + K_{52}(\log N)^{2/3}(A^{\prime})^{1/6} \\
  &\qquad \qquad \qquad \qquad \qquad + \frac{3}{4}K_{52}(\log N)^{2/3}
    (A-A^{\prime})^{1/6} + \frac{1}{80}(k-1)K\log N\biggr) \notag \\
  &\leq \exp\biggl(-u(K_{52}(\log N)^{2/3},A) 
    - \frac{1}{20}d^{\prime}\biggr). \notag
\end{align}

\emph{Case 2}.  $A^{\prime} < \tfrac{1}{3}A$.  In
this case it is easily checked that we must
have $k \geq 2$ and $d + 1 \geq 3w_{1}\sqrt{A/3}$ for 
$G_{N}(k,A,A^{\prime},d,d^{\prime})$ to be nonempty.  But then
\[
  \frac{9}{10}d - \frac{1}{20}kK\log N \geq \frac{4}{5}(d+1)
  \geq w_{1}\sqrt{A} + \frac{1}{20}d^{\prime},
\]
so (\ref{E:indhyp3}) follows from (\ref{E:indhyp2}).

The proof of (\ref{E:indhyp3}) for $j = k$ is now complete.
Summing 
(\ref{E:indhyp3}) as in (\ref{E:sumbound}),
\begin{align} \label{E:sumbound2}
  P_{N,w}&(\sum_{\gamma \in \mathfrak{C}_{N}}
    |\Int(\gamma)| \geq A, T_{N}^{\prime} \geq d^{\prime}) \\
  &\leq \sum_{A \leq B \leq |\Lambda_{N}|} \quad 
    \sum_{d \leq \diam_{\tau}(\Lambda_{N})} \quad
    \sum_{d^{\prime} \leq n \leq d} \quad \sum_{A^{\prime} \leq B} 
    \quad \sum_{k \leq n} P_{N,w}(G_{N}(k,B,A^{\prime},d,n)) 
    \notag \\
  &\leq K_{85}N^{4}\sum_{B \geq A} \quad \sum_{n \geq d^{\prime}}
    \exp\left(-\frac{1}{20}n - u(K_{52}(\log N)^{2/3},B)\right) \notag \\
  &\leq \exp\left(-\frac{1}{20}d^{\prime} - w_{1}\sqrt{A} 
    + K_{86}(\log N)^{2/3}l^{1/3}\right). \notag
\end{align}
If $K_{87}$ and $K_{7}$ are large enough, 
then (\ref{E:sumbound2}), Theorem
\ref{T:lowerbound} and Remark \ref{R:wired} show that
\begin{align} \label{E:bigTNpr}
  P_{N,w}&\left(\sum_{\gamma \in \mathfrak{C}_{N}}
    |\Int(\gamma)| \geq A, T_{N}^{\prime} \geq 
    K_{87}l^{1/3}(\log N)^{2/3}\right) \\
  &\leq \exp\left(-\frac{1}{40}K_{87}l^{1/3}(\log N)^{2/3}\right)
    P_{N,w}(|\Int(\Gamma_{0})| \geq A) \notag
\end{align}
and
\[
  P_{N,w}\left(\sum_{\gamma \in \mathfrak{C}_{N}}
    |\Int(\gamma)| \geq 2A\right) \leq e^{-\epsilon_{37}l}
  P_{N,w}(|\Int(\Gamma_{0})| \geq A).
\]
Therefore
\begin{align} \label{E:maxTN}
  P_{N,w}&\left(\sum_{\gamma \in \mathfrak{C}_{N}}
    |\Int(\gamma)| \geq A, |\mathfrak{C}_{N}| > 1\right) \\
  &\leq P_{N,w}\left(\sum_{\gamma \in \mathfrak{C}_{N}}
    |\Int(\gamma)| \geq A, T_{N}^{\prime} \geq K\log N\right) \notag 
    \displaybreak[0]\\
  &\leq P_{N,w}\left(\sum_{\gamma \in \mathfrak{C}_{N}}
    |\Int(\gamma)| \geq A, T_{N}^{\prime} \geq 
    K_{87}l^{1/3}(\log N)^{2/3}\right) \notag \\
  &\qquad + P_{N,w}\left(\sum_{\gamma \in \mathfrak{C}_{N}}
    |\Int(\gamma)| \geq 2A\right) \notag \\
  &\qquad + P_{N,w}\left(A \leq \sum_{\gamma \in \mathfrak{C}_{N}}
    |\Int(\gamma)| < 2A, K\log N \leq T_{N}^{\prime} <
    K_{87}l^{1/3}(\log N)^{2/3}\right) \notag \\
  &\leq 2\exp\left(-\frac{1}{40}K_{87}l^{1/3}(\log N)^{2/3}\right)
    P_{N,w}\left(\sum_{\gamma \in \mathfrak{C}_{N}}
    |\Int(\gamma)| \geq A\right) \notag \\
  &\qquad + P_{N,w}\left(A \leq \sum_{\gamma \in \mathfrak{C}_{N}}
    |\Int(\gamma)| < 2A, K\log N \leq T_{N}^{\prime} <
    K_{87}l^{1/3}(\log N)^{2/3}\right). \notag
\end{align}
To prove (\ref{E:single}), then, we need to bound the last 
probability in (\ref{E:maxTN}). We will sum as in 
(\ref{E:sumbound2}), but this time using (\ref{E:indhyp1}) instead of
(\ref{E:indhyp3}).  Let $a = 72w_{1}^{-2}K_{87}^{2}
l^{2/3}(\log l)^{4/3}$.
By (\ref{E:params}), we need only consider
\[
  A^{\prime} \geq B - 36w_{1}^{-2}(d^{\prime}+1)^{2} \geq B - a.
\]
Therefore using (\ref{E:indhyp2}),
\begin{align} \label{E:sumbound4}
  &P_{N,w}\left(A \leq \sum_{\gamma \in \mathfrak{C}_{N}}
    |\Int(\gamma)| < 2A, K\log N \leq T_{N}^{\prime} <
    K_{87}l^{1/3}(\log N)^{2/3}\right) \\
  &\quad \leq \sum_{A \leq B < 2A} \quad \sum_{B - a < 
    A^{\prime} \leq B} \quad \sum_{d \leq \diam_{\tau}(\Lambda_{N})} 
    \quad \sum_{d^{\prime} \geq  
    K\log N} \quad \sum_{k \leq d^{\prime}} 
    P_{\Lambda_{N}}(G_{N}(k,B,A^{\prime},d,d^{\prime}))
    \notag \\
  &\quad \leq \sum_{A^{\prime} \geq A - a} \quad \sum_{A^{\prime} \leq B < 
    A^{\prime} + a} \quad \sum_{d^{\prime} \geq K\log N} \notag \\
  &\qquad \qquad 4\kt N(d^{\prime}+1)e^{-\frac{4}{5}d^{\prime}}
    P_{N,w}
    (|\Int(\Gamma_{y})| = A^{\prime} \text{ for some }
    y \in \Lambda_{N} \cap \ZZ) \notag \\
  &\quad \leq e^{-\frac{1}{2}K\log N} 
    P_{N,w}
    (|\Int(\Gamma_{y})| \geq A - a \text{ for some }
    y \in \Lambda_{N} \cap \ZZ). \notag
\end{align}
By Lemma \ref{L:inward2} and (\ref{E:DeltaA}),
\begin{align} \label{E:goodcase}
  &P_{N,w}
    (|\Int(\Gamma_{y})| \geq A - a \text{ for some }
    y \in \Lambda_{N} \cap \ZZ) \\
  &\qquad \leq 2P_{N,w}
    (|\Int(\Gamma_{y})| \geq A - a, \notag \\
  &\qquad \qquad \qquad \qquad \qquad \Delta_{A-a}(\Gamma_{y}) <
    K_{63}l^{2/3}(\log l)^{1/3} \text{ for some }
    y \in \Lambda_{N} \cap \ZZ). \notag
\end{align}
With (\ref{E:sumbound4})
and Proposition \ref{P:changeA2} (taking $\delta = K_{57}\log A$)
this shows that, if $K$ is sufficiently large,
\begin{align} \label{E:conclude}
  P_{N,w}&(A \leq \sum_{\gamma \in \mathfrak{C}_{N}}
    |\Int(\gamma)| < 2A, K\log N \leq T_{N}^{\prime} <
    K_{87}l^{1/3}(\log N)^{2/3}) \\
  &\leq 2e^{-\frac{1}{2}K\log N} P_{N,w}
    \bigl(|\Int(\Gamma_{y^{\prime}})| \geq A - a \text{ for some } 
    \notag \\
  &\qquad  \qquad \qquad \qquad (\Lambda_{N},(A-a)/2,
    (1+\epsilon_{24})(A-a))-\text{compatible }
    y^{\prime}\bigr) \notag \\
  &\leq 2e^{-\frac{1}{2}K\log N + \epsilon_{25}\delta} 
    \sum_{y^{\prime} \in \Lambda_{N} \cap \ZZ} P_{N,w}
    (|\Int(\Gamma_{y^{\prime}})| \geq A)\notag \\
  &\leq e^{-\frac{1}{4}K\log N} P_{N,w}
    \left(\sum_{\gamma \in \mathfrak{C}_{N}} |\Int(\gamma)| \geq A\right).
    \notag
\end{align}
With (\ref{E:maxTN}) this proves (\ref{E:single}).  Then
(\ref{E:single}) and Lemma \ref{L:inward2} prove 
(\ref{E:ALRmax2})--(\ref{E:Hausmax2}).

It remains to establish (\ref{E:bnfree}).  Let $\epsilon_{38} > 0$ to
be specified, 
let $n_{0} = \min\{n: 2^{n}\kt r/3 > \epsilon_{38}\sqrt{A}\}$,
and let $b_{n} = \epsilon_{26}^{-2}2^{2n}(\kt r/3)^{2}$.
Then provided $\epsilon_{38}$ is small enough, we have
\[
  b_{n} < \epsilon_{33}2^{n-1}\frac{\kt r}{3}\sqrt{A}
  \quad \text{for all } n \leq n_{0},
\]
with $\epsilon_{33}$ from Lemma \ref{L:inward2}.
We have
\begin{align} \label{E:bneckbd}
  P_{N,w}&(|\Int(\Gamma_{y})| \geq A \text{ and }
    \Gamma_{y} \text{ is not } (q,r)-\text{bottleneck-free for some } 
    y \in \Lambda_{N} \cap \ZZ) \\
  &\leq P_{N,w}(|\Int(\Gamma_{y})| \geq A \text{ and }
    D^{\prime}_{(q,r)}(\Gamma_{y}) \geq \frac{\kt r}{3} 
    \text{ for some } y \in \Lambda_{N} \cap \ZZ) \notag \\
  &\leq P_{N,w}(|\Int(\Gamma_{y})| \geq A \text{ and }
    D^{\prime}_{(q,r)}(\Gamma_{y}) \geq \epsilon_{38}\sqrt{A}
    \text{ for some } y \in \Lambda_{N} \cap \ZZ) \notag \\
  &\qquad + \sum_{n = 1}^{n_{0}}
    P_{N,w}\Bigl(|\Int(\Gamma_{y})| \geq A, 
    2^{n-1}\frac{\kt r}{3} \leq
    D^{\prime}_{(q,r)}(\Gamma_{y}) < 2^{n}\frac{\kt r}{3}, \notag \\
  &\qquad \qquad \qquad \qquad
    |\Int(\Gamma_{y})| - |\Int(\alpha_{\max,\Gamma_{y}})| < b_{n}
    \text{ for some } y \in \Lambda_{N} \cap \ZZ\Bigr) \notag 
    \displaybreak[0]\\
  &\leq |\Lambda_{N} \cap \ZZ|\exp\left(-\frac{1}{20}
    \epsilon_{38}\sqrt{A} - 
    u(K_{52}(\log A)^{2/3},A)\right) \notag \\
  &\qquad + \sum_{n = 1}^{n_{0}} \exp\left(-\epsilon_{34}2^{n-1}
    \frac{\kt r}{3}\right) P_{N,w}(|\Int(\Gamma_{y})| \geq A
    \text{ for some } y \in \Lambda_{N} \cap \ZZ) \notag \\
  &\leq \exp\left(-\frac{1}{40}\epsilon_{38}\sqrt{A}\right)
    P_{N,w}(|\Int(\Gamma_{0})| \geq A) \notag \\
  &\qquad + K_{88}\exp\left(-\epsilon_{34}\frac{\kt r}{3}\right)
    P_{N,w}(|\Int(\Gamma_{y})| > A
    \text{ for some } y \in \Lambda_{N} \cap \ZZ), \notag \\
  &\leq K_{89}\exp(-\epsilon_{39}r)
    P_{N,w}\left(\sum_{\gamma \in \mathfrak{C}_{N}}
    |\Int(\gamma)| > A\right). \notag
\end{align}
Here the second inequality uses Lemma \ref{L:AAprime}, the third 
inequality uses Theorem \ref{T:onebound} and Lemma \ref{L:inward2}, 
and the fourth
inequality follows from Theorem \ref{T:lowerbound}.
With (\ref{E:single}) this proves (\ref{E:bnfree}).
\end{proof}

\begin{proof}[Proof of Theorem \ref{T:main}, (\ref{E:MLRmax}) and
(\ref{E:bnkfree})]
For (\ref{E:MLRmax}) we need only observe that 
Lemma \ref{L:inward2} is valid for $P$
in place of $P_{N,w}$ and only $y = 0$ considered.  After these same
changes, (\ref{E:bnkfree}) follows from (\ref{E:bneckbd}) with 
the last inequality omitted.
\end{proof}

\section{Conditioning on the Exact Area} \label{S:exact}

In Theorems \ref{T:main}---\ref{T:single}, 
one could as well consider conditioning on
$|\Int(\Gamma_{0})| = A$ in place of $|\Int(\Gamma_{0})| \geq A$.  It 
is straightforward to alter the existing
proof to prove these theorems under this 
conditioning once one has a lower bound like Theorem \ref{T:lowerbound}
for $P(|\Int(\Gamma_{0})| = A)$.  For this we need first some definitions 
and lemmas.  In the interest of space we will not give full details in the 
proofs; these details involve many of the same technicalities we have
encountered earlier.
Consider distinct points $x,y \in \mathbb{R}^{2}$ with $|y - x| \geq 
4\sqrt{2}$.  We let
$U(x,y)$ denote the open slab between the tangent line to 
$\partial B_{\tau}(x,\tau(y-x))$ at $y$ and the parallel line through $x$; we
call $U(x,y)$ the \emph{natural slab} of $x$ and $y$.  (Note the tangent
line is not necessarily unique; if it is not we make some arbitrary choice.)
We have $U(x,y) = U(y,x)$.
It follows from the definition 
of $U(x,y)$ that if $u$ and $v$ are on opposite sides of 
$U(x,y)$, then 
\begin{equation} \label{E:taubound}
  \tau(v - u) \geq \tau(y - x).
\end{equation}
The portion of $U(x,y)$ which is strictly to the right of the line from
$x$ to $y$ is called the \emph{natural half-slab} of $x$ and $y$ and
denoted $U_{R}(x,y)$.  For $x, y \in U(u,v)$ we let 
$\tilde{U}_{uv}(x,y)$ denote the open slab with sides parallel to those of
$U(u,v)$ with $x$ in one edge of $\tilde{U}_{uv}(x,y)$ and $y$ in the 
other edge.  We let $\tilde{U}_{uv}^{R}(x,y) = \tilde{U}_{uv}(x,y) \cap
U_{R}(u,v)$.

One or two of the dual sites adjacent to $x$ are
in $U(x,y)$; we let $x^{\prime}$ denote such a site, making an
arbitrary choice if there are two.  We define $y^{\prime}$ analogously
as a dual site in $U(x,y)$ adjacent to $y$.

By a \emph{face of the dual lattice} we mean a square
$z + [-\tfrac{1}{2},\tfrac{1}{2}]^{2}$ with $z \in \ZZ$.  Let $V(x,y)$
denote the interior of the union of all faces of the dual lattice whose
interiors are contained in
$\tilde{U}_{xy}(x^{\prime},y^{\prime})$, and let 
$T_{x}(x,y)$ and $T_{y}(x,y)$ denote
the components of $V(x,y)^{c}$ containing $x$ and $y$, respectively.  
(These components are distinct since $|x - y| \geq 4\sqrt{2}$.)  Then $V(x,y)
\subset \tilde{U}_{xy}(x^{\prime},y^{\prime}) \subset U(x,y)$.

The following lemma says roughly that typical open dual paths from $x$ to $y$ 
are connected to the boundary of the natural slab only near $x$ and $y$. 

\begin{lemma} \label{L:extraneous}
  Let $P$ be a percolation model on $\mathcal{B}(\ZZ)$ satisfying 
  (\ref{E:assump}).  There exist constants $K_{90}, K_{91},
  \epsilon_{40}$ as follows.  
  For all $x, y \in (\ZZ)^{*}$ and $J \geq K_{90} \log |x - y|$, 
  \[
    P(x \lra z \text{ for some } z \in \partial V(x,y) \backslash 
    (x + \Lambda_{J}) \mid x \lra y)
    \leq K_{91}e^{-\epsilon_{40}J}.
  \]
\end{lemma}
\begin{proof}  As mentioned above, we omit some details.  Suppose $x \lra y$ 
and $x \lra z$ for some $z \in \partial T_{x}(x,y) \backslash (x + \Lambda_{J})$,
via open dual paths.  If $K_{92}$ is large, one can trivially dispose of the case
in which $x \lra B(x,K_{92}|y - x|)^{c}$, so we hence forth tacitly consider only
connections occuring inside $B(x,K_{92}|y - x|)$; in particular this means
$|z - x| \leq K_{92}|y - x|$.
There are then two cases:  either there is 
an $(\epsilon_{41}J)$-near connection from $x$ to $y$ 
in $B(z,J/10)^{c}$, or
there is not; here $\epsilon_{41}$ is to be 
specified.  In the first
case,  there is also an open dual path from $z$ to 
$\partial (B(z,J/20) \cap (\ZZ)^{*})$,
so we can apply (\ref{E:conlwr}) and Lemmas \ref{L:ratiowm} and \ref{L:rnear} 
(assuming $K_{90}$ is large) to obtain
\begin{align} \label{E:Jnear}
  P&(N(x,y,\epsilon_{41}J,B(z,J/10)^{c}) \cap 
    [z \lra \partial (B(z,J/20) \cap (\ZZ)^{*})]) \\
  &\leq K_{93}Je^{-\tau(y-x) + K_{60}\epsilon_{41}J}
    e^{-\kt J/40} \notag \\
  &\leq K_{94}e^{-\epsilon_{42}J}P(x \lra y), \notag
\end{align}
provided $\epsilon_{41}$ is chosen small enough and $K_{90}$ large enough.
In the second case, there exists dual sites $v,w$ just outside 
$B(z,J/10)$
and open dual paths $x \lra v$ and $w \lra y$ occuring at separation 
$\epsilon_{41}J$.  It follows easily from (\ref{E:taubound}) and the fact that
$z$ is close to $\partial U(x,y)$ that 
\[
  \tau(y - w) \geq \tau(y - x) - \tfrac{1}{5}\kt J.
\]
Also, since $|z - x| \geq J$, 
\begin{equation} \label{E:taucomp2}
  \tau(v - x) \geq \tau(z - x) - \tfrac{1}{5}\kt J
  \geq \tfrac{1}{2}\kt J. 
\end{equation}
Therefore by Lemma \ref{L:decouple} and (\ref{E:conlwr}),
\begin{align} \label{E:noJnear}
  P&(N(x,y,\epsilon_{41}J,B(z,J/10)^{c})^{c} \cap 
    [x \lra y] \cap [x \lra z]) \\
  &\leq \sum_{v,w} 2P(x \lra v)P(w \lra y) \notag \\
  &\leq K_{95}J^{2}e^{-\tau(y - x) - \tfrac{3}{10}\kt J} \notag \\
  &\leq K_{96}e^{-\epsilon_{43}J}P(x \lra y). \notag
\end{align}
Now (\ref{E:Jnear}) and (\ref{E:noJnear}), summed over $z$ with
$|z - x| \leq K_{92}|y - x|$, prove the lemma.
\end{proof}

For dual sites $x$ and $y$, we say that $x \lra y$ 
\emph{cylindrically} if there
is an open dual path $\gamma$ from $x$ to $y$ in $U(x,y)$ and every open 
dual path from $\gamma$ to $U(x,y)^{c}$ passes through $x$ or $y$.

For $\mathcal{D}$ a subgraph of the dual lattice (or just a set of dual
bonds, which we may view as such a subgraph), and $A \subset \mathbb{R}^{2}$,
we define the \emph{bond boundary of} $\mathcal{D}$ \emph{in} $A$ to be the set
of bonds contained in $A$ having exactly one endpoint in $\mathcal{D}$.
(As always, we view bonds as open intervals contained in the plane.)

\begin{lemma} \label{L:slabconn}
  Let $P$ be a percolation model on $\mathcal{B}(\ZZ)$ satisfying 
  (\ref{E:assump}).  There exist constants $K_{97}, \epsilon_{44}$ as follows.  
  For all $x, y \in (\ZZ)^{*}$, 
  \[
    P(x \lra y \text{ cylindrically})
    \geq \epsilon_{44}|y - x|^{-K_{97}}P(x \lra y).
  \]
\end{lemma}
\begin{proof}
Fix $x, y$ and let $J = K_{90}\log |y - x|$, with $K_{90}$ from Lemma
\ref{L:extraneous}. Let $\mathcal{J}_{x}(x,y)$ denote the bond boundary of 
$\partial T_{x}(x,y) \cap \mathcal{B}(x + \Lambda_{J})$ in $\mathbb{R}^{2}$,
and define $\mathcal{J}_{y}(x,y)$ analogously.  For $e \in 
\mathcal{J}_{x}(x,y)$ let $A_{e}$ denote the event that all dual bonds
in $\mathcal{B}(\partial T_{x}(x,y) \cap (x + \Lambda_{J}))$ 
are open, $e$ 
and $\langle xx^{\prime} \rangle$ are open 
and all other bonds in $\mathcal{J}_{x}(x,y)$ are closed; define
$A_{e}$ analogously for $e \in \mathcal{J}_{y}(x,y)$.
Define the event
\[
  B = [x \lra z \text{ for some } z \in \partial V(x,y) \backslash 
        (x + \Lambda_{J})]
\]
(cf. Lemma \ref{L:extraneous}.)  Given a configuration $\omega \in 
[x \lra y] \cap B^{c}$, there necessarily exists an open dual path from 
$\partial T_{x}(x,y) \cap (x + \Lambda_{J})$ to 
$\partial T_{y}(x,y) \cap (y + \Lambda_{J})$ in $V(x,y)$ which contains 
only one bond in $\mathcal{J}_{x}(x,y)$ and one bond in 
$\mathcal{J}_{y}(x,y)$; we denote these two bonds by $b_{xy}(x,\omega)$ 
and $b_{xy}(y,\omega)$, respectively, making an arbitrary choice if more 
than one choice is possible.  If the configuration $\omega$ is in
$[x \lra y] \cap B^{c} \cap [b_{xy}(x) = e] \cap [b_{xy}(y) = f]$ for 
some $e, f$, then we can modify at most $16J$ bonds (those in 
$\partial T_{x} \cap \mathcal{B}(\partial T_{x}(x,y) \cap (x + 
\Lambda_{J}))$) to
obtain a  configuration in $A_{e} \cap A_{f}$; in the resulting configuration 
we have $x \lra y
\text{ cylindrically}$. Therefore from the bounded energy property we have
\[
  P(x \lra y \text{ cylindrically} \mid [x \lra y] \cap B^{c})
  \geq e^{-K_{98}J}
\]
This and Lemma \ref{L:extraneous} prove the lemma.
\end{proof}

It is easy to check that Lemmas \ref{L:extraneous} and \ref{L:slabconn}
are valid if we restrict to connections in the halfspace $H_{xy}$.  More
precisely, under the assumptions of Lemma \ref{L:slabconn} we have
\begin{equation} \label{E:halfslab}
    P(x \lra y \text{ cylindrically in } U_{R}(x,y))
    \geq \epsilon_{44}|y - x|^{-K_{97}}P(x \lra y).
\end{equation}
Additionally, we can extend the idea of cylindrical connections
as follows:  for $u, v \in \mathbb{R}^{2}$ and $x, y$ dual sites in 
$U(u,v)$, we say that $x \lra y \ (u,v)$-\emph{cylindrically} if there is
an open dual path from $x$ to $y$ in $\tilde{U}_{uv}(x,y)$ and the dual cluster
of $x$ and $y$ intersects $\partial \tilde{U}_{uv}(x,y)$ only at $x$ and $y$.
Provided $x, y \in U_{R}(u,v)$, the proof of Lemma \ref{L:slabconn} shows that
\begin{equation} \label{E:halfslab2}
    P(x \lra y \ (u,v)\text{-cylindrically in } \tilde{U}_{uv}^{R}(x,y))
    \geq \epsilon_{44}|y - x|^{-K_{97}}P(x \lra y \text{ in } H_{xy}).
\end{equation}

\begin{theorem} \label{T:lowerbound2}
  Let $P$ be a percolation model on $\mathcal{B}(\ZZ)$ satisfying
  (\ref{E:assump}), the near-Markov property for open circuits,
  positivity of $\tau$
  and the ratio weak mixing property.  There exist
  $K_{i}$ such that for $A > K_{99}$ and $l = \sqrt{A}$,
  \[
    P(|\Int(\Gamma_{0})| = A) \geq \exp(-w_{1}\sqrt{A} - K_{100}l^{1/3}
    (\log l)^{2/3}).
  \]
\end{theorem}
\begin{proof}
Fix $A$ large and let $l = \sqrt{A}$.  Let $a_{1}$ denote
the vertical coordinate of the point where $\partial\mathcal{K}_{1}$
meets the positive vertical axis.  Let $s = l^{2/3}(\log l)^{1/3}$ 
and $\delta = K_{30}s^{2}/l$, with $K_{30}$ as in the proof of Theorem
\ref{T:lowerbound}.
Let $\alpha = \partial (l + \delta)\mathcal{K}_{1}$ and let 
$(z^{\prime}_{0},..,z^{\prime}_{n},z^{\prime}_{0})$ 
be the $s$-hull skeleton of $\alpha$.  
It is an easy exercise in geometry to   
see that the natural half-slabs $U_{R}(z^{\prime}_{i},
z^{\prime}_{i+1}), i = 0,..,n$,
are disjoint.  (Our labeling as usual is cyclical: 
$z^{\prime}_{n+1} = z^{\prime}_{0}$.)
For some $K_{101}$ to be specified,
let us call a pair $(z^{\prime}_{i},z^{\prime}_{i+1})$ from the skeleton
\emph{very short} if $|z^{\prime}_{i+1} - z^{\prime}_{i}| 
\leq 2\sqrt{2}$,
\emph{short} if $2\sqrt{2} < |z^{\prime}_{i+1} - z^{\prime}_{i}| 
\leq 2K_{101} \log l$ and \emph{long} if
$|z^{\prime}_{i+1} - z^{\prime}_{i}| > 2K_{101} \log l$.  In what follows, very
short pairs can be handled quite trivially but tediously, so for convenience
we will assume there are no very short pairs.
For long pairs we
define $x_{i}^{\prime}$ and $y_{i+1}^{\prime}$ 
to be the points on the line segment
$\overline{z^{\prime}_{i}z^{\prime}_{i+1}}$ at distance 
$K_{101}\log l$ from $z^{\prime}_{i}$ and
from $z^{\prime}_{i+1}$, respectively, and let $x_{i}, y_{i+1}$ be dual sites in
$U_{R}(z^{\prime}_{i},z^{\prime}_{i+1})$ within distance $\sqrt{2}$ of 
$x_{i}^{\prime}$ and $y_{i+1}^{\prime}$, respectively. 
For short pairs we let $x_{i} = y_{i+1}$ be a dual site in
$U_{R}(z^{\prime}_{i},z^{\prime}_{i+1})$ within distance $\sqrt{2}$ 
of the midpoint of 
$\overline{z^{\prime}_{i}z^{\prime}_{i+1}}$.
With minor modification of the definition of the
$s$-hull skeleton, 
we may assume the set $\{z^{\prime}_{0},..,z^{\prime}_{n}\}$ has 
lattice symmetry, that is, for each
$z^{\prime}_{i}$ the reflection of $z^{\prime}_{i}$ 
across the horizontal or vertical axis 
is another $z^{\prime}_{j}$, and analogously for the 
sites $x_{i}$ and $y_{i}$.  For each $i$
we let $\phi_{i}$ denote a dual lattice path of  minimal length from $y_{i}$ to
$x_{i}$ outside
$\Co(\{z^{\prime}_{0},..,z^{\prime}_{n}\}) \cup 
\tilde{U}_{z^{\prime}_{i-1}z^{\prime}_{i}}(x_{i-1},y_{i})
\cup \tilde{U}_{z^{\prime}_{i}z^{\prime}_{i+1}}(x_{i},y_{i+1})$.  
We call such a $\phi_{i}$ a \emph{short link}.  Let
$\mathcal{C}_{i}$ denote the bond boundary of
$\phi_{i}$ in $(\Co(\{z^{\prime}_{0},..,z^{\prime}_{n}\}) \cup 
\tilde{U}_{z^{\prime}_{i-1}z^{\prime}_{i}}(x_{i-1},y_{i})
\cup \tilde{U}_{z^{\prime}_{i}z^{\prime}_{i+1}}(x_{i},y_{i+1}))^{c}$.

Let $\lambda_{a}$ be the  vertical line through $(a,0)$.  
Let $H_{L}(x)$ and $H_{R}(x)$
denote the open half planes to the left and right, respectively, of 
the vertical line through $x$.
Let $H_{U}(x)$ and $H_{B}(x)$ denote the open half planes above and below
the horizontal line through $x$, respectively.  
(In general we use the convention that subscripts 
$L, R, U, B$ refer to left, right, upper and lower halfspaces, respectively, with
combinations, such as $LU$, referring to quadrants.)
Let $S(x,y)$ denote the open slab between the vertical lines through $x$ and $y$.
Let $N$ be the integer part of $a_{1}l/2$, 
$M$ the integer part of $l^{2/3}(\log l)^{1/3}$
and $D$ the integer part of $l^{1/3}(\log l)^{2/3}$.
Let 
\[
  u_{RU}^{\prime} = H_{R}(0) \cap H_{U}(0) \cap \alpha
  \cap \lambda_{M+\frac{1}{2}}, \quad v_{RU}^{\prime} = 
  H_{R}(0) \cap H_{U}(0) \cap \alpha 
  \cap \lambda_{M+D+\frac{1}{2}},
\]
\[
  w_{RU}^{\prime} = H_{R}(0) \cap H_{U}(0) \cap \alpha 
  \cap \lambda_{N+\frac{1}{2}}, \quad x_{RU}^{\prime} = 
  H_{R}(0) \cap H_{U}(0) \cap \alpha 
  \cap \lambda_{N+D+\frac{1}{2}}.
\]
(Note each of these intersections is a single point.)  We call these 4 
points \emph{determining points}. Lattice symmetry
yields corresponding determining points with appropriate subscripts in the other
three quadrants.  We may assume that $u_{RU}^{\prime}$ is 
one of the sites $z^{\prime}_{i}$ of the 
$s$-hull skeleton of $\alpha$ (if not, we add $u_{RU}^{\prime}$ to the
skeleton), and  analogously for the other sites just defined.
Let $u_{RU}$ 
be the second closest dual site above $u_{RU}^{\prime}$ in 
$\lambda_{M+\frac{1}{2}}$, and analogously for $v_{RU},
w_{RU}, x_{RU}$.  If $u_{RU}^{\prime} = z^{\prime}_{i}$, 
for some $i$, we define $z_{i}$ to 
be $u_{RU}$, and again analogously for the other determining points.
Loosely, the idea is to remove from $\Gamma_{0}$ its intersection with each of 
the width-$D$ vertical slabs $S(x_{LU},w_{LU}), S(v_{LU},u_{LU}), S(u_{RU},
v_{RU})$, $S(w_{RU},x_{RU})$, then raise or lower the segments of $\Gamma_{0}$
between these slabs to adjust the area as desired, then reconnect these 
segments to make a new circuit enclosing area $A$.  To do this we must 
first ensure that 
$\Gamma_{0}$ intersects each vertical line bounding any of these four slabs only
twice.

We refer to the 4 width-$D$ vertical slabs above as \emph{removal slabs}.  We
call the 5 regions $H_{L}(x_{LU}),S(w_{LU},v_{LU}),S(u_{LU},u_{RU}),
S(v_{RU},w_{RU}),H_{R}(x_{RU})$ (whose closures together
form the complement of the 4 removal slabs) \emph{retention regions}.
By a \emph{retained segment} we mean a connected
component of the intersection of $\alpha$ with a retention 
region.  Each retained
segment has the form $\alpha^{(z^{\prime}_{j},z^{\prime}_{k})}$ 
for some $j, k$; we call 
$z^{\prime}_{j}$ an \emph{initial determining point} and 
$z^{\prime}_{k}$ a \emph{final determining
point}, and call $(j,k)$ a \emph{retention pair}.
We let $J^{ret}$ denote the set of all 8 retention pairs.
For each initial determining point $z_{j}^{\prime}$, in the
boundary of some retained region $F$, we let $\psi_{j}$ be 
a dual lattice path from $z_{j}$ to $x_{j}$ in $F \backslash 
\tilde{U}_{z_{j}^{\prime}z_{j+1}^{\prime}}(x_{j},y_{j+1})$, of minimal
length, and let $\mathcal{D}_{j}$ be the bond boundary of 
$\psi_{j}$ in $\overline{F} \backslash 
  \tilde{U}_{z_{j}^{\prime}z_{j+1}^{\prime}}(x_{j},y_{j+1}) \}$.
For each final determining point $z_{k}^{\prime}$ we let 
$\psi_{k}$ be 
a dual lattice path from $y_{k}$ to $z_{k}$ in $F \backslash 
\tilde{U}_{z_{k-1}^{\prime}z_{k}^{\prime}}(x_{k-1},y_{k})$, of minimal
length, and define $\mathcal{D}_{k}$ analogously to
$\mathcal{D}_{j}$.
We refer to $\psi_{j}$ and $\psi_{k}$
as the \emph{endpaths} of the retention pair $(j,k)$.  

For each retention pair $(j,k)$ let
$I_{jk} = \{i:z^{\prime}_{i}, z^{\prime}_{i+1} \in 
\alpha^{(z^{\prime}_{j},z^{\prime}_{k})} \}$ and let
$Q_{jk}$ denote the event that (i) for 
each $i \in I_{jk} \cup \{j\}$, we have $x_{i} \lra y_{i+1}\ 
(z^{\prime}_{i},z^{\prime}_{i+1})$-cylindrically in 
$\tilde{U}^{R}_{z^{\prime}_{i}z^{\prime}_{i+1}}(x_{i},y_{i+1})$, (ii) 
for each $i \in I_{jk}$ we have $\phi_{i}$ open and all bonds in 
$\mathcal{C}_{i}$ closed,
and (iii) both endpaths of $(j,k)$
are open and all bonds in $\mathcal{D}_{j} \cup
\mathcal{D}_{k}$ are closed. These 3 component events are denoted 
$Q_{(i)}(j,k)$, 
$Q_{(ii)}(j,k)$
and $Q_{(iii)}(j,k)$.  For a configuration in
$Q_{jk}$, the paths $x_{i} \lra y_{i+1}$
together with the short links $\phi_{i}$ and the two endpaths form an open
dual path from $z_{j}$ to $z_{k}$ outside $\Co(\{z^{\prime}_{0},..,
z^{\prime}_{n}\})$, and there is no open dual connection from this path
to any point of the retention region boundary except $z_{j}$ and $z_{k}$.
By Lemma \ref{L:ratiowm},
(\ref{E:halfslab2}) and Theorem \ref{T:halfspace},
provided $K_{101}$ is large we have
\begin{align} \label{E:Qibound}
  P&(Q_{(i)}(j,k)) \\
  &\geq \left(\tfrac{1}{2}\right)^{|I_{jk}|} \prod_{i \in I_{jk} \cup \{j\}}
    P(x_{i} \lra y_{i+1}\ 
    (z^{\prime}_{i},z^{\prime}_{i+1})-\text{cylindrically in }
    \tilde{U}^{R}_{z^{\prime}_{i}z^{\prime}_{i+1}}(x_{i},y_{i+1})) \notag \\
  &\geq \left(\frac{\epsilon_{45}}{l}\right)^{|I_{jk}|+1}
    \exp\left(-\sum_{i \in I_{jk} \cup \{j\}} \tau(y_{i+1} - x_{i})\right) 
    \notag \\
  &\geq \exp\left(-\sum_{i \in I_{jk} \cup \{j\}} \tau(z_{i+1} - z_{i}) -
    K_{102}|I_{jk}| \log l\right). \notag
\end{align}
From the bounded energy property, 
\[
  P\bigl(Q_{(ii)}(j,k) \cap Q_{(iii)}(j,k) \mid Q_{(i)}(j,k)\bigr) \geq
    \exp(-K_{103}|I_{jk}| \log l)
\]
which with (\ref{E:Qibound}) yields
\begin{equation} \label{E:allQbound}
  P(Q_{jk}) \geq
  \exp\left(-\sum_{i \in I_{jk} \cup \{j\}} \tau(z_{i+1} - z_{i}) -
    K_{104}|I_{jk}| \log l \right).
\end{equation}

For a configuration $\omega \in Q_{jk}$, 
and for $F$ the retention region with $z^{\prime}_{j},z^{\prime}_{k} \in 
\partial F$, we can
associate an area $R_{jk}(\omega)$ as follows.
There is a unique outermost 
open dual path $\Xi_{jk}(\omega)$ from $z_{j}$
to $z_{k}$ in $F$.  If $F$ is a halfspace ($H_{L}(x_{LU})$ or 
$H_{R}(x_{RU}))$, then 
$R_{jk}(\omega)$ is the area of the region between
$\Xi_{jk}(\omega)$ and
$\overline{z_{j}z_{k}}$.  If $F$ is a slab and $z^{\prime}_{j}, 
z^{\prime}_{k} \in H_{U}(0)$, then 
$R_{jk}(\omega)$ is the area of the region 
in $F \cap H_{U}(0)$ below $\Xi_{jk}(\omega)$.
If $F$ is a slab and $z^{\prime}_{j}, 
z^{\prime}_{k} \in H_{B}(0)$, then 
$R_{jk}(\omega)$ is the area of the region 
in $F \cap H_{B}(0)$ above $\Xi_{jk}(\omega)$.
We define corresponding nonrandom areas $R_{jk}^{0}$
similarly but using $\alpha^{(z^{\prime}_{j}z^{\prime}_{k})}$ in place of 
$\Xi_{jk}(\omega)$.  Then 
$R_{jk}(\omega) \geq 
R_{jk}^{0}$.

It is not hard to see that for some $K_{105}$ we have
\begin{align} \label{E:maxarea}
  P&\bigl(R_{jk} \geq R_{jk}^{0}
    + K_{105}l^{4/3}(\log l)^{2/3} \text{ for some} \\
  &\qquad \qquad \text{retention pair } (j,k) 
    \mid \cap_{(j,k) \in J^{ret}} \ Q_{jk} \bigr) 
    \leq \frac{1}{2}. \notag
\end{align}
In fact, if this were false we could obtain extra area 
$K_{105}l^{4/3}(\log l)^{2/3}$ almost ``for free'' in Theorem \ref{T:lowerbound};
more precisely, one could replace $A$ with $A + K_{105}l^{4/3}(\log l)^{2/3}$ 
on the left side in the conclusion of that theorem.  But, assuming $K_{105}$ 
is large, this would contradict Theorem \ref{T:onebound}.  
It follows from (\ref{E:maxarea}) that for each 
retention pair $(z^{\prime}_{j},z^{\prime}_{k})$ there exists 
$a_{jk} \in [R_{jk}^{0},
R_{jk}^{0} + K_{105}l^{4/3}(\log l)^{2/3})$
such that 
\begin{align} \label{E:fixedarea}
  P&\left(\cap_{(j,k) \in J^{ret}} [R_{jk} = a_{jk}] \mid
    \cap_{(j,k) \in J^{ret}} \ Q_{jk} \right) \\
  &\geq \frac{1}{2}(K_{105}l^{4/3}(\log l)^{2/3})^{-8} \notag.
\end{align}
Let
\[
  R_{1} = \sum_{(j,k) \in J^{ret}} 
    a_{jk},
\]
where the sum is over the 8 retention pairs.

We call $(k,j)$ a \emph{removal pair} if $\alpha^{(z^{\prime}_{k},
z^{\prime}_{j})}$ is a connected component of the intersection of $\alpha$
with some removal slab.
Let $J^{rem}$ 
denote the set of all 8 removal pairs.
For each removal
pair $(k,j)$ and corresponding removal 
slab $F$, let $\chi_{kj}$
be a dual path from $z_{k}$ to $z_{j}$ in $\overline{F}
\backslash \Co(\{z^{\prime}_{0},..,z^{\prime}_{n}\})$, and let 
$G_{kj}$ denote the event that all bonds in
$\chi_{kj}$ are open, while all bonds in the bond boundary of 
$\chi_{kj}$ in $(\psi_{j} \cup 
\psi_{k})^{c}$ are closed.
We call $\chi_{kj}$
a \emph{long link}.  There are 2 long links in each removal slab, one
each in the upper and lower half planes.
Let $R_{2}$ be the total area in the 4 removal slabs, between 
the upper and lower long link in each slab.  Assuming the long links are
chosen to have length of order $D$ (say, at most $4D$),
we have from the bounded energy 
property that
\begin{align} \label{E:Gbound2}
  P&\bigl(\cap_{(k,j) \in J^{rem}} \ G_{z^{\prime}_{k},z^{\prime}_{j}}
    \mid (\cap_{(j,k) \in J^{ret}} \ Q_{jk}) \cap
    (\cap_{(j,k) \in J^{ret}} [R_{jk} = a_{jk}]) \bigr) \\
  &\geq e^{-K_{106}D} \notag \\
  &\geq e^{-K_{107}l^{1/3}(\log l)^{2/3}}. \notag
\end{align}

Define the event
\[
  E = \bigl(\cap_{(j,k) \in J^{ret}} \ Q_{jk} \bigr) \cap
  \bigl(\cap_{(j,k) \in J^{ret}} [R_{jk} = a_{jk}] \bigr)
  \cap \bigl(\cap_{(k,j) \in J^{rem}} \ G_{kj} \bigr).
\]
For a configuration $\omega \in E$, there is an open 
dual circuit surrounding $\Co(\{z^{\prime}_{0},..,z^{\prime}_{n}\})$ 
satsifying the constraint that it
include all of the short links $\phi_{i}$, long links $\chi_{kj}$
and endpaths $\psi_{j}$.  There is a unique outermost such
circuit subject to this constraint, obtained by taking the outermost 
$(z^{\prime}_{i},z^{\prime}_{i+1})$-cylindrical connection from $x_{i}$ 
to $y_{i+1}$ for each $i$; we denote this circuit $\Gamma_{1}(\omega)$.
Because of the cylindrical nature of these connections and the closed state of
the bond boundaries of the short and long links, we have $\Gamma_{1}(\omega)
= \Gamma_{0}(\omega)$, unless $\Gamma_{0}(\omega)$ and $\Gamma_{1}(\omega)$ are
disjoint with $\Gamma_{0}(\omega)$ surrounding $\Gamma_{1}(\omega)$ and no
open dual path connecting $\Gamma_{0}(\omega)$ to $\Gamma_{1}(\omega)$.  It
therefore follows from the near-Markov property that
\[
  P(\Gamma_{0} \neq \Gamma_{1} \mid E) \leq e^{-\epsilon_{46}l}
  \leq \frac{1}{2}.
\]
By (\ref{E:maxvert2}),
\[
  \sum_{(j,k) \in J^{ret}} |I_{jk}| \leq K_{108}\frac{l}{s}
  \leq K_{109}l^{1/3}(\log l)^{-1/3}.
\]
Using these facts with (\ref{E:fixedarea}), (\ref{E:Gbound2}), 
Lemma \ref{L:ratiowm},
(\ref{E:allQbound}) and (\ref{E:lengthhull}) (which is still valid here), we obtain
\begin{align} \label{E:exactsize}
  P&(|\Gamma_{0}| = R_{1} + R_{2}) \\
  &\geq P(E \cap [\Gamma_{0} = \Gamma_{1}]) \notag \\
  &\geq \frac{1}{2}P(E) \notag \\
  &\geq \frac{1}{2}(K_{105}l^{4/3}(\log l)^{2/3})^{-8}
    e^{-K_{107}l^{1/3}(\log l)^{2/3}}
    P(\cap_{(j,k) \in J^{ret}} \ Q_{jk})
    \notag \\
  &\geq e^{-K_{110}l^{1/3}(\log l)^{2/3}} \prod_{(j,k) \in J^{ret}}
    P(Q_{jk}) \notag \\
  &\geq \exp\left(-\sum_{i=0}^{n} \tau(z_{i+1} - z_{i}) -
    K_{111}l^{1/3}(\log l)^{2/3} \right) \notag \\
  &\geq \exp\left(-w_{1}\sqrt{A} - K_{112}l^{1/3}(\log l)^{2/3} \right). \notag
\end{align}

Let $\theta \omega$ be the upward shift of a configuration $\omega$ by 1 unit,
and for an event $G$ let $\theta^{m}G = \{\omega: \theta^{-m}\omega \in G\}$.  
Let $J_{B}^{ret}$ be the set consisting of the 3 retention pairs corresponding
to segments of $\alpha$ in the
lower half plane.  Given a constant $K_{113}$, for $m \leq 
K_{113}l^{1/3}(\log l)^{2/3}$ we can replace $Q_{jk}$
with $\theta^{m}Q_{jk}$ for all 
$(j,k) \in J_{B}^{ret}$ throughout the argument leading
to (\ref{E:exactsize}) at the expense of only a possible increase in $K_{111}$,
provided we alter the 4 long links $\chi_{kj}$
in the outer 2 removal slabs to connect to the appropriate shifted sites
$z_{i} + (0,m)$ instead of to $z_{i}$.  (The possible increase in $K_{111}$ 
reflects a possible reduction in the probabilities of the events 
$G_{kj}$, resulting in an increase of
$K_{107}$ in (\ref{E:Gbound2}).)
We can readily keep the area $R_{2}$
fixed when we so alter the long links.  We thereby obtain
\[
  P(|\Gamma_{0}| = R_{1} + R_{2} - (2N + 1)m) \geq 
  \exp(-w_{1}\sqrt{A} - K_{114}l^{1/3}(\log l)^{2/3}).
\]
Provided $K_{113}$ is large, we can choose $m \in \mathbb{Z}$ so that
\[
  |R_{1} + R_{2} - (2N + 1)m - A| \leq N.
\]
We can then repeat this entire argument, but shift upward (by some amount $q
\leq K_{113}l^{1/3}(\log l)^{2/3}$) only the event
$Q_{jk}$ for the central of the 3 retention pairs
in $J_{B}^{ret}$.  This gives
\begin{equation} \label{E:nearsize}
  P(|\Gamma_{0}| = R_{1} + R_{2} - (2N + 1)m - (2M + 1)q) \geq 
  \exp(-w_{1}\sqrt{A} - K_{115}l^{1/3}(\log l)^{2/3}).
\end{equation}
We can choose $k$ so that 
\[
  |R_{1} + R_{2} - (2N + 1)m - (2M + 1)q - A| \leq M 
  \leq l^{2/3}(\log l)^{1/3}.
\]
But it is easy to see that one can alter the long links to change $R_{2}$
by any amount up to $l^{2/3}(\log l)^{1/3}$, so that
\[
  R_{1} + R_{2} - (2N + 1)m - (2M + 1)q = A,
\]
at the expense of only a possible increase to 
$K_{115}$ in (\ref{E:nearsize}).  With (\ref{E:nearsize})
this completes the proof.
\end{proof}

Now that we can use Theorem \ref{T:lowerbound2} in place of 
Theorem \ref{T:lowerbound}, all proofs leading to Theorems \ref{T:main}---
\ref{T:single} remain valid under conditioning on $|\Int(\Gamma_{0})| = A$
in place of $|\Int(\Gamma_{0})| \geq A$.  This establishes the following
results.

\begin{theorem} \label{T:newmain}
Under the assumptions in Theorem \ref{T:main} or 
Theorem \ref{T:FKcase}, under the measure
$P(\cdot \mid |\Int(\Gamma_{0})| = A)$ the conclusions
(\ref{E:ALRmax}) -- (\ref{E:bnkfree}) hold with probability approaching
1 as $A \to \infty$.
\end{theorem}

\begin{theorem} \label{T:newsingle}
Under the assumptions in Theorem \ref{T:single}, under the measures
$P_{N,w}(\cdot \mid \sum_{\gamma \in \mathfrak{C}_{N}} 
|\Int(\Gamma_{0})| = A)$ the
conclusions (\ref{E:single}) -- (\ref{E:bnfree}) 
hold with probability approaching 1 as $A \to \infty$.
\end{theorem}


\begin{thebibliography}{99}

\bibitem{Al92} Alexander, K.S., \emph{Stability of the Wulff
minimum and fluctuations in shape for large finite clusters in
two-dimensional percolation}, Probab. Theory Rel. Fields
\textbf{91} (1992), 507-532.

\bibitem{Al97} Alexander, K.S., \emph{Approximation of
subadditive functions and rates of convergence in limiting
shape results}, Ann. Probab. \textbf{25} (1997), 30-55.

\bibitem{Al97mix} Alexander, K.S., \emph{On weak mixing
in lattice models}, Probab. Theory Rel. Fields \textbf{110}
(1998), 441-471.

\bibitem{Al97pwr}  Alexander, K.S., \emph{Power-law
corrections to exponential decay of connectivities
and correlations in lattice models}, preprint (1997).

\bibitem{Al97arc} Alexander, K.S., \emph{The asymmetric random
cluster model and comparison of Ising and Potts models},
preprint (1997).

\bibitem{ACC} Alexander, K.S., Chayes, J.T. and Chayes, L.,
\emph{The Wulff construction and asymptotics of the finite cluster
distribution for two dimensional Bernoulli percolation}, Commun.
Math. Phys. \textbf{131} (1990), 1-50.

\bibitem{AVSZ} Avron, J.E., van Beijeren, H., Schulman, L.S. and Zia,
R.K.P., \emph{Roughening transition, surface
tension and equilibrium droplet shapes in a two-dimensional Ising system},
J. Phys. A \textbf{15} (1982), L81-L86.

\bibitem{BDJ} Baik, J., Deift, P. and Johansson, K., 
\emph{On the distribution of the length of the longest increasing subsequence of
random permutations}, J. Amer. Math. Soc. \textbf{12} (1999), 1119-1178.

\bibitem{BuK} Burton, R. and Keane, M., \emph{Density and uniqueness
in percolation}, Commun. Math. Phys. \textbf{121} (1989), 501-505.

\bibitem{DH} Dobrushin, R.L. and Hryniv, O., \emph{Fluctuations 
of the phase boundary
in the $2D$ Ising ferromagnet}, Commun. Math. Phys. \textbf{189} (1997),
395-445.

\bibitem{DKS} Dobrushin, R.L., Koteck\'y, R. and Shlosman, S., 
\emph{Wulff construction. A global shape from local
interaction},  Translations of Mathematical Monographs, 104, American 
Mathematical Society, Providence (1992).

\bibitem{ES} Edwards, R.G. and Sokal, A.D., \emph{Generalization of the
Fortuin-Kasteleyn-Swendsen-Wang representation and Monte Carlo algorithm},
Phys. Rev. D \textbf{38} (1988), 2009-2012.

\bibitem{FK} Fortuin, C.M. and Kasteleyn, P.W., \emph{On
the random cluster model. I. Introduction and relation to other models},
Physica \textbf{57} (1972), 536-564.

\bibitem{Gr95} Grimmett, G.R., \emph{The stochastic
random-cluster process and uniqueness of random-cluster measures},
Ann. Probab. \textbf{23} (1995), 1461-1510.

\bibitem{Hr} Hryniv, O., \emph{On local behaviour of the phase 
separation line in the $2D$ Ising model}, Probab. Theory Rel. Fields
\textbf{110} (1998), 91-107.

\bibitem{IS97} Ioffe, D. and Schonmann, R.H., \emph{Dobrushin-Koteck\'y-Shlosman 
theorem up to the critical temperature}, Commun. Math. Phys. \textbf{199}
(1998), 117-167.

\bibitem{KS91} Krug, J. and Spohn, H., \emph{Kinetic roughening of 
growing interfaces}, in \emph{Solids Far from Equilibrium:  Growth, Morphology
and Defects} (C. Godr\`eche, ed.) 479-582, Cambridge University Press,
Cambridge (1991).

\bibitem{LMR} Laanait, L., Messager, A. and Ruiz, J., \emph{Phase
coexistence and surface tensions for the Potts model}, Commun. Math.
Phys. \textbf{105} (1986), 527-545.

\bibitem{LNP} Licea, C., Newman, C.M. and Piza, M.S.T., \emph{Superdiffusivity in
first-passage percolation}, Probab. Theory Rel. Fields \textbf{106} (1996), 
559-591.

\bibitem{MW} McCoy, B.M. and Wu, T.T., \emph{The Two-Dimensional Ising Model},
Harvard University Press, Cambridge, USA, 1973.

\bibitem{NP}  Newman, C.M. and Piza, M.S.T., \emph{Divergence of shape fluctuations
in two dimensions}, Ann. Probab. \textbf{23} (1995), 977-1005.

\bibitem{Pi} Piza, M.S.T., \emph{Directed polymers in a random
environment:  Some results on fluctuations}, J. Statist. Phys.
\textbf{89} (1997), 581-603.

\bibitem{Ta1}  Taylor, J.E., \emph{Existence and structure of solutions to
a class of nonelliptic variational problems}, Symp. Math. \textbf{14}
(1974), 499-508.

\bibitem{Ta2} Taylor, J.E., \emph{Unique structure of solutions to a class of
nonelliptic variational problems}, Proc. Symp. Pure Math. \textbf{27} (1975),
419-427.

\end{thebibliography}
\end{document}